# MINIMIZING AREA AMONG LAGRANGIAN SURFACES: THE MAPPING PROBLEM

R. SCHOEN AND J. WOLFSON

## 1. Introduction

This paper introduces a geometrically constrained variational problem for the area functional. We consider the area restricted to the lagrangian surfaces of a Kähler surface, or, more generally, a symplectic 4-manifold with suitable metric, and study its critical points and in particular its minimizers. We apply this study to the problem of finding canonical representatives of the lagrangian homology (that part of the homology generated by lagrangian cycles). We show that the lagrangian homology of a Kähler surface (or of a symplectic 4-manifold) is generated by lagrangian surfaces that are branched immersions except at finitely many singular points (Theorem 8.1). We precisely describe the structure of these singular points. In particular, these singular points are represented by lagrangian cones with an associated local Maslov index. Only those cones of Maslov index 1 or $-1$ may be area minimizing (Propositions 7.3 and 7.5). Our surfaces are minimizers of area among suitable classes of lagrangians. Their mean curvature satisfies a first-order system of partial differential equations of "Hodge-type" that we describe precisely below.

We show that an integral homology class may be represented by a lagrangian cycle (a cycle whose two-simplices are given by $C^1$ lagrangian maps) if and only the class is annihilated by the symplectic form $\omega$ (Proposition 2.1). On the other hand it is known that a homology class may be represented by an immersed lagrangian surface if and only if it is annihilated both by the symplectic form and the first Chern form. For a class which is annihilated by $\omega$, we are able to construct an area minimizing representative among lagrangian competitors and to show that it has a finite number of singularities whose local Maslov classes add to $\frac{1}{2}c_1(\alpha)$. In particular, nonflat singularities do arise in minimizing lagrangian surfaces (the note added in proof in [SW] incorrectly states the contrary). In case we have $c_1(\alpha) = 0$, it follows that there are an even number of singular points whose total sum of Maslov indices is 0. We do not know in general whether the minimizer

*Date*: August 30, 2018.

The first author was partially supported by NSF grant DMS-9803341. The second author was partially supported by NSF grant DMS-9802487. Parts of this paper were written while the second author was visiting the American Institute of Mathematics in Palo Alto, Ca. Both authors are grateful to AIM for their generous support of this project.

is free of singularities (other than branch points) in this case. We do show however that if a minimizing $S^2$ represents a class $\beta \in \pi_2(N)$, and if any integer ($> 1$) multiple of this $S^2$ also minimizes area in its homotopy class, then this $S^2$ is free of singularities (other than branch points) (Theorem 8.4).

A situation of great geometric interest in this area is the construction of special lagrangian cycles in Calabi-Yau manifolds (see [SYZ] for the connection with mirror symmetry). The variational method we employ here provides an approach to this construction. In particular, if an immersed lagrangian submanifold $\Sigma$ of a Kähler-Einstein manifold is stationary for volume, it is automatically minimal (special lagrangian in the Calabi-Yau case) (Lemma 8.2). This may not be true if $\Sigma$ is singular. In particular, in our two dimensional setting, the surface will be minimal provided there are no nonflat singular points. Note that in the Kähler-Einstein case, it is true that if $\omega(\alpha) = 0$ then $c_1(\alpha) = 0$, and hence the sum of the Maslov indices of the nonflat singular points is zero. We show that the nonflat cones are not minimizing if one allows nonorientable comparison surfaces (Proposition 7.8). This suggests that a lagrangian minimizer in a $\mathbb{Z}_2$ homology class will be minimal, and that the $\mathbb{Z}_2$ homology of any Kähler-Einstein surface (it is easy to show that every $\mathbb{Z}_2$ homology class may be represented by a lagrangian cycle) can be represented by minimal lagrangian cycles. We will study the $\mathbb{Z}_2$ homology in a forthcoming paper.

The technical work in this paper involves constructing area minimizing lagrangian maps from a surface into a symplectic 4 manifold. We follow the variational approach of first constructing such minimizers in the Sobolev class $W^{1,2}$, and then discussing the regularity of the minimizers. Comparison results of Gromov [G] and Allcock [Al] are used to show the Hölder continuity of such maps (Theorem 2.8). In contrast to the situation for the classical Plateau problem, the higher regularity is a difficult issue. This is because in the classical case the maps are harmonic, while in our case, the minimization of energy among lagrangian maps does not seem to imply any useful regularity. To get such regularity we derive a monotonicity formula (Proposition 3.2), and use it to obtain a partial regularity theorem and the existence of tangent cones (Theorem 4.10). The monotonicity formula is of a rather subtle nature. The setting for this formula is that of contact stationary surfaces in the Heisenberg group. A difficulty which we encounter here in proving the regularity is the issue of controlling both the map and the surface, and as such it involves a combination of classical and Geometric Measure Theory ideas. We believe the methods developed here will be useful in other problems involving the minimization of area among surfaces with geometric constraints (e.g. embeddedness).

We now introduce the technical setting for this problem and describe the PDE which we will be studying. Let $N$ be a symplectic $2n$-manifold with symplectic form $\omega$. We will equip $(N, \omega)$ with what is sometimes called a compatible metric $g$ and almost complex structure $J$. That is, the triple $(\omega, g, J)$ is required to satisfy, for $v, w \in TN$: $\omega(v, Jw) = g(v, w)$, $\omega(Jv, Jw) = \omega(v, w)$. An important example of this structure is the case that $(N, g, \omega, J)$ is a Kähler manifold with Kähler metric $g$, Kähler form $\omega$ and complex structure $J$. A vector field $X$ on $N$ that preserves the symplectic form, i.e., that satisfies $\mathcal{L}_X \omega = 0$ where $\mathcal{L}$ denotes the Lie derivative, is called a *symplectic vector field*. The diffeomorphisms it generates are called *symplectic diffeomorphisms*. Since $\mathcal{L}_X \omega = d(X \lrcorner \omega)$, it follows that $X$ is a symplectic vector field if and only if $X \lrcorner \omega$ is closed. If $X \lrcorner \omega$ is exact, i.e. $X \lrcorner \omega = dh$, for a function $h$ on $N$, $X$ is called a *hamiltonian vector field* for the *hamiltonian* $h$. The diffeomorphisms it generates are called *hamiltonian diffeomorphisms*. Suppose that $\Sigma$ is a real $n$-dimensional manifold (compact or with boundary). We say that a smooth ($C^1$) map $\ell : \Sigma \to (N, \omega)$ is *lagrangian* if $\ell^* \omega = 0$. If $\ell$ is an immersion then using the compatible metric the lagrangian condition can be formulated as follows: $\ell(\Sigma)$ is lagrangian if at each point $x \in \ell(\Sigma)$ and for each tangent vector $v \in T_x(\ell(\Sigma))$ then $\langle Jv, w \rangle = 0$ for all $w \in T_x(\ell(\Sigma))$ (i.e. $JT_x(\ell(\Sigma)) \perp T_x(\ell(\Sigma))$). We remark that if $\Sigma \subset \mathbb{R}^{2n}$ is given as a graph over the standard real $n$-plane; i.e. $\Sigma$ is given as $y = f(x)$ where $f$ is a map from $\mathbb{R}^n$ to $\mathbb{R}^n$, then $\Sigma$ is lagrangian if and only if there is a function $u(x)$ defined locally such that $f = \nabla u$. This is the simplest way to construct lagrangian submanifolds of $\mathbb{R}^{2n}$.

Lagrangian submanifolds enjoy the following remarkable property. Let $\Sigma$ be an immersed lagrangian submanifold and suppose that $(N, \omega)$ is a Kähler manifold. The mean curvature vector $H$ of $\Sigma$ in $N$ is a section of the normal bundle on $\Sigma$. Define a 1-form on $\Sigma$ by: $\sigma_H = H \lrcorner \omega$. Then a geometric computation (see the Appendix) gives: $d\sigma_H = \text{Ric}|_\Sigma$ where Ric denotes the Ricci 2-form on $N$. In particular, suppose $N$ is Kähler-Einstein so that the metric satisfies: $\text{Ric} = \text{R}\omega$ where R is the scalar curvature. Then since $\Sigma$ is lagrangian $d\sigma_H = 0$ and therefore $H$ is an infinitesimal symplectic motion.

This observation suggests that it is natural to consider variational problems for volume with a lagrangian constraint. Recall the an immersed submanifold is called *stationary* if the volume is critical for arbitrary smooth compactly supported variations. It is well-known that an immersed submanifold is stationary if and only if it is minimal, i.e. $H = 0$. We will call an immersed lagrangian submanifold $\Sigma$ *lagrangian stationary* if the volume is critical for arbitrary smooth variations through lagrangian submanifolds. We will call $\Sigma$ *hamiltonian stationary* if the volume is critical for compactly supported hamiltonian variations. If $X$ is a hamiltonian vector field then $X \lrcorner \omega = dh$ or, equivalently $X = J\nabla h$. It follows that if $\Sigma$ is an immersed lagrangian submanifold that is hamiltonian stationary then,

$$0 = \int_\Sigma \langle H, X \rangle = \int_\Sigma \langle \sigma_H, dh \rangle = \int_\Sigma \langle \delta \sigma_H, h \rangle.$$

Since this holds for any smooth function $h$ with compact support we conclude that $\delta \sigma_H = 0$. Combining this with the computation of $d\sigma_H$, we see that if

$\Sigma$ is an immersed hamiltonian stationary lagrangian submanifold then,

$$\begin{aligned} d\sigma_H &= \text{Ric}_{|\Sigma}, \\ \delta\sigma_H &= 0. \end{aligned} \quad (1.1)$$

Thus, if $N$ is Kähler-Einstein then $\sigma_H$ is a harmonic 1-form.

The equations (1.1) are third order in $\Sigma$. To understand the nature of these equations, we consider the special case in which $\Sigma \subset \mathbb{R}^{2n}$ is given by the graph of the gradient of a function $u(x)$. That is, $\Sigma$ is given by $\{(x, y = \nabla u(x)) : x \in \mathbb{R}^n\}$ for a scalar function $u(x)$. In this case, the volume of $\Sigma$ is given by $\int \sqrt{\det(g_{ij})}\, dx$ where $g_{ij} = \delta_{ij} + \sum_k u_{ik} u_{kj}$. Doing a variation of the form $u + tv$ with $v$ of compact support one derives the Euler Lagrange equation

$$\sum_{k=1}^{n} \frac{\partial}{\partial x_k}(\Delta_\mu(u_k)) = 0$$

where $\Delta_\mu$ is the Laplace operator with respect to the induced metric $\mu$. In particular we see that this is a quasilinear fourth order elliptic equation whose linearization at $u(x) = 0$ is the biharmonic equation. We remark that this set of hamiltonian stationary equations has been studied from the point of view of integrable systems by Hélein and Romon [HR].

## 2. Preliminaries

We will suppose, *for the remainder of the paper*, that $N$ denotes a symplectic 4-manifold with symplectic form $\omega$, compatible metric $g$ and almost complex structure $J$.

An integral lagrangian cycle in $N$ is an integral cycle with the 2-simplices given by $C^1$ lagrangian maps. We say a homology class $\alpha \in H_2(N, \mathbb{Z})$ is a *lagrangian homology class* if it can be represented by an integral lagrangian cycle.

**Proposition 2.1.** *A homology class $\alpha \in H_2(N, \mathbb{Z})$ is a lagrangian homology class if and only if $[\omega](\alpha) = 0$.*

*Proof.* Suppose that $\alpha \in H_2(N, \mathbb{Z})$ satisfies $[\omega](\alpha) = 0$. Represent $\alpha$ by an embedded surface $\Sigma \subset N$ and let $T$ denote a tubular neighborhood of $\Sigma$ in $N$. Since $\int_\Sigma \omega = 0$, $\omega$ is exact on $T$ and there is a one-form $\eta$ such that on $T$, $\omega = d\eta$. Note that $\eta$ is well-defined up to addition of a closed one-form. Choose a triangulation of $\Sigma$ that is $\varepsilon$-fine, in the sense that every one-simplex has length less than $\varepsilon$. Let $\Delta$ denote a two-simplex of the triangulation. Then

$$\int_\Delta \omega = \int_{\partial\Delta} \eta,$$

and the second integral is independent of the choice of $\eta$. Let $s$ denote a one-simplex in $\partial\Delta$. Then

$$|\int_s \eta| < C\varepsilon,$$

where the constant C depends only on $\eta$. Perturbing $s$ in $T$ keeping its endpoints fixed we can construct a one-simplex $\tilde{s}$ in $T$ with $\int_{\tilde{s}} \eta = 0$. Doing this for each one-simplex in the triangulation we have a one-complex in $T$ with the same vertices as the original triangulation and with the property that the integral of $\eta$ around the perturbation of $\partial \Delta$ is zero, for any two-simplex $\Delta$. It follows that the perturbation of $\partial \Delta$ spans a lagrangian two-simplex. The result follows. □

Note that the requirement $[\omega](\alpha) = 0$ imposes a rationality condition on $[\omega]$. This requirement while weaker than $[\omega] \in H^2(N, \mathbb{Q})$ nevertheless implies that for "generic" symplectic structure $\omega$ there are no nontrivial lagrangian classes.

**Proposition 2.2.** *Let $\alpha \in H_2(N, \mathbb{Z})$ be a lagrangian homology class. Then $\alpha$ can be represented by a piecewise $C^1$ lagrangian map $\ell : \Sigma \to N$ which is an immersion except in a neighborhood of finitely many points $\{x_1, \ldots, x_k\} \subset \Sigma$. At each point $x_j, j = 1, \ldots, k$ there is a well-defined local Maslov index $m_{x_j}$. The total Maslov index satisfies:*

$$\sum_{j=1}^{k} m_{x_j} = \frac{1}{2} c_1(N)(\alpha).$$

*Proof.* By Proposition 2.1 $\alpha$ can be represented by a lagrangian cycle. The cycle can be chosen so that along the one-simplices the tangent planes of adjacent two-simplices are transverse. We smooth this cycle along the one-simplices outside neighborhoods of the vertices, as follows. Choose a one-simplex $s$. We can suppose that $s$ lies in a Darboux ball. After a hamiltonian diffeomorphism we can suppose that $s$ is a segment of a coordinate axis. $s$ is the boundary of two smooth lagrangian simplices $L_1$ and $L_2$. Subdividing $s$, if necessary, we can find a lagrangian plane $P$ containing $s$ such that near $s$, $L_1$ is a graph over the half plane $H_1$ and $L_2$ is a graph over the half plane $H_2$, where $P = H_1 \cup H_2$. In particular, for $i = 1, 2$, in a neighborhood of $s$, $L_i$ is the graph of $\nabla f_i$, where $f_i$ is a smooth function on $H_i$. Let $V$ be a neighborhood of $s$ in $P$. Let $f$ be a smooth function on $P$ that coincides with $f_1$ on $H_1 \setminus V$ and with $f_2$ on $H_2 \setminus V$. The required smoothing is effected by glueing in the graph of $\nabla f$. Note that if $s$ must be subdivided to find $P$ then the smoothings can be constructed to agree at the subdivision points.

Let $x_j$ be a singular point. Let $D_r(x_j) \subset \Sigma$ be a disc of radius $r$, center $x_j$. Choose $r$ so that $D_r(x_j)$ contains no singular points except $x_j$, so that $\ell(D_r(x_j))$ lies in a Darboux neighborhood $U$ and so that $\ell$ is an immersion on $\partial D_r(x_j)$. Let $D_j$ be an oriented immersed disc in $U$ with $\partial D_j = \partial D_r(x_j)$ such that the induced orientations agree. The lagrangian immersion $\ell$ determines a trivialization of $TN$ along $\Sigma \setminus \{x_1, \ldots, x_k\}$ and in particular along $\partial D_r(x_j)$. Trivialize $TN$ on $D_j$ and compare the two trivializations on $\partial D_j$. The result is an element of $\pi_1(U(2)) \simeq \mathbb{Z}$. Choosing a generator of $\pi_1(U(2))$ associates an even integer, $2m_{x_j}$, to the point $x_j$. This integer is well-defined independent of the choices of $r$ and $D_j$, provided $D_j \subset U$ and $\ell$ is an immersion along $\partial D_r(x_j)$. We will say that $m_{x_j}$ is the *local Maslov index* of $x_j$

If, at each singular point $x_j$, we replace $\ell(D_r(x_j))$ with $D_j$ we construct a surface $\tilde{\Sigma}$, no longer lagrangian, that still represents $\alpha$. Join $x_1$ to $x_2$, $x_2$ to $x_3$, ..., and $x_{k-1}$ to $x_k$ by simple curves on $\Sigma \setminus \{x_1, \ldots, x_k\}$. Cutting along these curves we can realize $\tilde{\Sigma}$ as the union of $\Sigma \setminus \cup_j^k D_r(x_j)$ and $\cup_j^k D_j$ glued along an $S^1$. To compute $c_1(N)(\tilde{\Sigma})$ we compare the trivializations of $TN$ on the boundaries of $\Sigma \setminus \cup_j^k D_r(x_j)$ and $\cup_j^k D_j$. By the previous computation it follows that:

$$c_1(N)(\alpha) = c_1(N)(\tilde{\Sigma}) = 2 \sum_{j=1}^k m_{x_j}.$$

$\square$

The following proposition, due to Gromov and Lees [Le], gives necessary and sufficient conditions for a smooth map to be homotopic to a lagrangian immersion.

**Proposition 2.3.** *(Gromov, Lees) A smooth map $\phi : \Sigma \to N$ is homotopic to a lagrangian immersion if and only if:*

1. $\phi^*[\omega] = 0$.

2. $\phi^* c_1(N) = 0$,

*where $[\omega]$ is the cohomology class of $\omega$ and $c_1(N)$ is the first Chern class of $(N, J)$.*

It follows that:

**Corollary 2.4.** *An integral homology class $\alpha$ that satisfies $[\omega](\alpha) = 0$ and $c_1(N)(\alpha) = 0$ can be represented by an immersed lagrangian surface.*

Let $\Sigma$ be a compact surface with Riemannian metric and with volume form $d\mu$. Let $\mathbb{R}^n$ denote euclidean $n$-space. We denote the Hilbert space of square integrable maps $\Sigma \to \mathbb{R}^n$ by $L^2(\Sigma, \mathbb{R}^n)$ and the Hilbert space of maps $\Sigma \to \mathbb{R}^n$ with square integrable first derivatives by $W^{1,2}(\Sigma, \mathbb{R}^n)$. The energy of a map $f \in W^{1,2}(\Sigma, \mathbb{R}^n)$ is:

$$E(f) = \frac{1}{2} \int_\Sigma \langle df, df \rangle dv.$$

By the Nash imbedding theorem we can suppose that $N$ is isometrically embedded in $\mathbb{R}^n$, for some $n$. Define:

$$W^{1,2}(\Sigma, N) = \{f \in W^{1,2}(\Sigma, \mathbb{R}^n) : f(x) \in N \text{ a.e. } x \in \Sigma\}.$$

Recall that closed bounded subsets of $W^{1,2}(\Sigma, N)$ are weakly closed in $W^{1,2}(\Sigma, N)$.

We need the following lemma about $W^{1,2}$ maps.

**Lemma 2.5.** *There exists $\varepsilon_0 > 0$ depending only on $N$ such that if $\ell \in W^{1,2}(S^2, N)$ satisfies $E(\ell) \leq \varepsilon_0$, then $\ell$ is zero in homology in the sense that $\int_{S^2} \ell^*(\gamma) = 0$ for any smooth closed 2-form $\gamma$ on $N$.*

*Proof.* We follow [ScU], and consider smoothing the map $\ell$. Let us think of $S^2$ as the standard unit sphere in $\mathbb{R}^3$. Let $t \in S^2$, and consider a euclidean ball of radius $h \leq 1/2$ about $t$. Let $\zeta(r)$ be a smooth nonincreasing function of $r \geq 0$ with support in $[0,1]$ and $\zeta(r) = 1$ for $r \leq 1/2$. Consider the averaged map

$$\ell^{(h)}(t) = (|S^2 \cap B_h(t)|)^{-1} \int_{S^2 \cap B_h(t)} \ell \, da.$$

The Poincaré inequality then implies

$$\int_{S^2 \cap B_h(t)} |\ell - \ell^{(h)}|^2 \, da \leq ch^2 \int_{S^2 \cap B_h(t)} |\nabla \ell|^2 \, da.$$

It follows that the mean square distance from $\ell^{(h)}(t)$ to $\ell$ is small in the ball of radius $h$. In particular the averaged map $\ell^{(h)}$ has image which is uniformly close to $N$ independent of $h$. We may thus project this map into $N$ using the nearest point projection map $\Pi$ which is smooth in a neighborhood of $N$. We let $\ell_h = \Pi \circ \ell^{(h)}$, and we have that the family of maps $\ell_h$ for $0 < h \leq 1/2$ is a family of smooth maps depending smoothly on $h$, and from standard results we have that the limit as $h \to 0$ exists in the $W^{1,2}$ norm and is equal to $\ell$. Since we assumed that the energy of $\ell$ is small, it follows that the image of $\ell_{1/2}$ lies in a small coordinate ball in $N$, and hence we have that $\int_{S^2} \ell_{1/2}^* \gamma = 0$. Now from the homotopy formula since $\gamma$ is a closed form we have $\frac{d\ell_h^* \gamma}{dh} = d\ell_h^*(\frac{\partial \ell}{\partial h} \lrcorner \gamma)$. It follows that $\int_{S^2} \ell_h^* \gamma$ is constant, and therefore 0 for all $h \in (0, 1/2]$. But since the maps $\ell_h$ converge strongly in $W^{1,2}$ to $\ell$ as $h$ tends to 0, we see that $\int_{S^2} \ell^* \gamma = 0$ as claimed. $\square$

The notion of a lagrangian map can be extended to $W^{1,2}(\Sigma, N)$ as follows: We say a map $\ell \in W^{1,2}(\Sigma, N)$ is *weakly lagrangian* if:

$$\ell^* \omega = 0$$

almost everywhere on $\Sigma$. Note that by Lebesgue point theory, we can, equivalently, define a map $\ell \in W^{1,2}(\Sigma, N)$ to be weakly lagrangian if:

$$\int_D \ell^* \omega = 0,$$

for every disc $D \subset \Sigma$.

**Proposition 2.6.** *The set of weakly lagrangian maps in $W^{1,2}(\Sigma, N)$ satisfying a uniform energy bound is closed in the weak topology.*

*Proof.* Let $\ell_i$ be a sequence in $W^{1,2}(\Sigma, N)$ with $E(\ell_i) \leq c$ which converges weakly to $\ell \in W^{1,2}(\Sigma, N)$. We must show that $\ell$ is lagrangian. Note that any subsequence of $\ell_i$ also converges weakly to $\ell$, so we may choose subsequences at will. We first consider the energy measures $1/2|\nabla \ell_i|^2 d\mu$, and using the weak convergence theorem for measures we may assume that these measures converge weakly as measures to a limit measure $\nu$; that is, we assume that for any open set $\Omega$ in $\Sigma$ we have $1/2 \int_\Omega |\nabla \ell_i|^2 d\mu \to \int_\Omega d\nu$. For any $\varepsilon_1 > 0$, there is at most a finite set $S$ for which $\nu(\{Q\}) \geq \varepsilon_1$ for $Q \in S$. If $\ell$ is not weakly lagrangian, then we can find a Lebesgue point $P \in \Sigma \setminus S$ such that $\ell^*\omega(P) \neq 0$. Thus it follows that for all $r$ sufficiently small we have $\int_{D_r(P)} \ell^*\omega \neq 0$. On the other hand, because $P$ is not in $S$, we have $E_{D_r(P)}(\ell_i) \leq 2\varepsilon_1$ for $i$ sufficiently large. We fix such a radius $r$, and observe that we may choose a radius $\rho \in [r/2, r]$ and a subsequence again denoted $\ell_i$ so that $\ell_{i|\partial D_\rho}$ has energy bounded by $2r^{-1}\varepsilon_1$. This implies that $\ell_{i|\partial D_\rho}$ is continuous, rectifiable and has length bounded by a fixed constant times $\varepsilon_1$. It follows that a subsequence converges uniformly to $\ell_{|\partial D_\rho}$. Since the length of $\ell_i(\partial D_\rho)$ is small we can suppose that the images of $\partial D_\rho$ all lie in a Darboux ball $B \subset N$. For each $i$, $\ell_i(\partial D_\rho)$ extends to a (not necessarily lagrangian) map $\hat{\ell}_i : D_i \to B$ with small energy. For $\varepsilon_1$ sufficiently small (depending only on $N$) we may apply Lemma 2.5 to show that the cycle $\ell_i(D_\rho) \cup \hat{\ell}_i(D_i)$ (which can be represented as a map with small energy from $S^2$ to $N$) is trivial in homology. In $B$ write $\omega = d\eta$, then we have

$$0 = \int_{D_\rho} \ell_i^*\omega + \int_{D_i} \hat{\ell}_i^*\omega = \int_{D_i} \hat{\ell}_i^*\omega = \int_{\partial D_\rho} \ell_i^*\eta,$$

where the last equality follows by Stokes' theorem (and $\hat{\ell}_i = \ell_i$ on $\partial D_\rho = \partial D_i$). Note that the condition that the energy of $\ell_i|_{\partial D_\rho}$ is bounded implies that there is a subsequence again denoted $\ell_i$ which converges weakly in $W^{1,2}(\partial D_\rho)$ to $\ell$. It follows that $\int_{\partial D_\rho} \zeta \ell_i^*\eta$ converges to $\int_{\partial D_\rho} \zeta \ell^*\eta$ for any $L^2(\partial D_\rho)$ function $\zeta$. (To see this observe that $\ell_i^*\eta$ is a linear combination of first derivatives of $\ell_i$ with coefficients which are smooth functions of $\ell_i$, and these coefficient functions converge in $L^2$ norm.) Choosing $\zeta = 1$, we conclude,

$$0 = \int_{\partial D_\rho} \ell^*\eta.$$

Repeating the above argument for $\ell$, we may find a small energy map $\hat{\ell} : D \to B$ which agrees with $\ell$ on $\partial D = \partial D_\rho$. As above it follows that

$$0 = \int_{D_\rho} \ell^*\omega + \int_D \hat{\ell}^*\omega = \int_{D_\rho} \ell^*\omega + \int_{\partial D_\rho} \ell^*\eta,$$

where the second equality follows by Stokes Theorem (and $\hat{\ell} = \ell$ on $\partial D_\rho = \partial D$). Therefore we have $\int_{D_\rho} \ell^*\omega = 0$. This contradiction shows that $\ell$ is weakly lagrangian. □

A classical invariant of lagrangian maps that plays a role in our considerations is the *period*. We specialize to the case $N = \mathbb{R}^4$ (with the euclidean metric $g$ and the standard symplectic structure $\omega = \sum_i^2 dx_i \wedge dy_i$) to describe this invariant. Let $\ell : \Sigma \to \mathbb{R}^4$ be a $C^1$ lagrangian map. and $\gamma : S^1 \to \Sigma$ be a closed curve. Choose a 1-form $\eta$ on $\mathbb{R}^4$ such that $d\eta = \omega$. We define the *period* of $\gamma$ to be: $\text{period}(\gamma) = \int_\gamma \ell^* \eta$. This number is well-defined independent of the choice of primitive $\eta$ and depends only on the homology class of $\gamma$ in $H_1(\Sigma, \mathbb{Z})$. If the periods of all closed curves on $\Sigma$ vanish we say that the lagrangian map is *exact*. Equivalently, all periods vanish if and only if $\ell^* \eta = d\varphi$ for a function $\varphi$ on $\Sigma$. There is another description of exact lagrangian maps. Consider $\mathbb{R}^5$ with coordinates $(x, y, \varphi)$ and projection $p : \mathbb{R}^5 \to \mathbb{R}^4, (x, y, \varphi) \mapsto (x, y)$. Let $\alpha = d\varphi - (\sum_i x_i dy_i - y_i dx_i)$ be the contact 1-form on $\mathbb{R}^5$. The contact distribution is the distribution of hyperplanes defined by $\alpha = 0$. A $C^1$ map of a surface $\lambda : \Sigma \to \mathbb{R}^5$ is called *legendrian* if $\lambda^*(\alpha) = 0$ or, equivalently, if $\lambda(\Sigma)$ is everywhere tangent to the contact distribution. A lagrangian map $\Sigma \to \mathbb{R}^4$ is exact if and only if it admits a legendrian lift $\Sigma \to \mathbb{R}^5$. It is easy to construct non-exact lagrangian maps into $\mathbb{R}^4$. However every $C^1$ lagrangian map $\ell : \Sigma \to N$ is *locally exact* in the sense that $\Sigma$ can be covered by discs $\{D_\lambda\}$ so that the image of each $D_\lambda$ under $\ell$ lies in some Darboux ball. The notion of local exactness plays an essential role in our regularity theory. In fact, it is possible to show that any weakly lagrangian map in $W^{1,2}(\Sigma, N)$ is locally exact in a suitable sense. However for our purposes this will not be necessary.

The following comparison result is due to Gromov ([G]) and Allcock ([Al]). We will use it to prove a Hölder continuity result and (in a later section) a strong compactness theorem for minimizing lagrangian maps.

**Proposition 2.7.** *(Gromov, Allcock) Given a $W^{1,2}$ map $\varphi$ from the unit circle $C$ into $\mathbb{R}^4$ satisfying the period condition $\int_C \varphi^*(xdy) = 0$, there exists a lagrangian map $\ell$ in $W^{1,2}(D, \mathbb{R}^4)$ with $\ell = \varphi$ on $C$, and*

$$A(\ell(D)) \leq cL(\ell(C))^2$$

*for a fixed constant c.*

We may now use the idea of Morrey ([M1]) to establish the Hölder continuity of minimizing lagrangian maps from a disk. We will consider minimizing lagrangian maps which are weakly conformal in the sense that $\|\ell_*(\partial/\partial x)\| = \|\ell_*(\partial/\partial y)\|$ and $\langle \ell_*(\partial/\partial x), \ell_*(\partial/\partial y)\rangle = 0$ a.e.. Such maps will be constructed later in this paper.

**Theorem 2.8.** *Let $\ell \in W^{1,2}(D, N)$ be a weakly conformal, minimizing lagrangian map. There exists $\alpha \in (0, 1)$ and $c$ depending only on $N$ such that (after redefinition on a set of measure 0), $\ell$ is Hölder continuous on $D_{1/2}$ and satisfies the bound*

$$d(\ell(P), \ell(Q)) \leq c(A(\ell(D)))^{1/2}|P - Q|^\alpha$$

*for $P, Q \in D_{1/2}$.*

*Proof.* Consider a point $P \in D_{1/2}$ and a radius $r < 1/2$. We want to show that $E_{D_r(P)}(\ell) \leq c(L(\ell(\partial D_r(P))))^2$ for a fixed constant $c$ depending only on $N$. To see this, observe that if the right hand side is greater than a fixed $\varepsilon_0 > 0$, the inequality is trivial, so it suffices to assume that $L = L(\ell(\partial D_r(P)))$ is small. Thus we may assume that the curve $\ell(\partial D_r(P))$ lies in a Darboux chart in which the metric is nearly euclidean. In this case it follows from Proposition 2.7 that there exists a map $\ell_0$ in $W^{1,2}(D_r(P), N)$ with $\ell_0 = \ell$ on $\partial D_r(P)$ and $A(\ell_0) \leq cL^2$. Since $\ell$ is area minimizing and weakly conformal we have $E_{D_r(P)}(\ell) \leq A(\ell_0) \leq cL^2$ as required. Now by the Schwarz inequality we have $L^2 \leq 2\pi r \frac{d}{dr} E_{D_r(P)}(\ell)$, and we have shown that $E_{D_r(P)}(\ell) \leq cr \frac{d}{dr} E_{D_r(P)}(\ell)$. Integrating this differential inequality we get the decay estimate $E_{D_r(P)}(\ell) \leq cr^{2\alpha} E(D)$ for a fixed $\alpha \in (0, 1)$. By Morrey's Lemma ([GT]) this implies the desired Hölder estimate. □

**Corollary 2.9.** *If $\ell \in W^{1,2}(\Sigma, N)$ is a weakly conformal, minimizing lagrangian map then it is locally exact.*

## 3. The First Variation of Area

Let $\Sigma$ be a compact Riemann surface, and $\ell \in W^{1,2}(\Sigma, M)$ be a weakly conformal, locally exact, lagrangian map. We will assume throughout this section that there is a disk $D$ in $\Sigma$ such that $\ell(D)$ is contained in a fixed Darboux chart $O$ in $N$. Note that by Theorem 2.8, this will be true on sufficiently small disks if the map is minimizing. We will also assume that $\ell$ is *stationary* in the sense that the first variation of area is zero for allowable deformations of $\ell$ through lagrangian maps. We first define the variations which we will consider. We suppose that $x, y$ are local Darboux coordinates in $O$, so that the symplectic form is in the standard form $\omega = dx \wedge dy$. There is an obvious class of variations given in terms of compactly supported hamiltonian functions $h(x, y)$, namely the flow of the vector field $X_h = h_x \frac{\partial}{\partial y} - h_y \frac{\partial}{\partial x}$. If $F_t$ denotes this flow, then $\ell_t = F_t \circ \ell$ is an allowable variation. It will be important to extend this class of variations to allow variations which include dilations, and other collapsing deformations. In order to do this we use the local exactness. This enables us to define a function $\varphi$ on $\ell^{-1}(O)$ which satisfies $d\varphi = \ell^*\eta$ where $\eta = (xdy - ydx)$. We wish to consider hamiltonians which depend on $\varphi$ as well as $x, y$. In order to do this, we may form the 5 dimensional space with coordinates $x, y, \varphi$. In this space we have the contact 1-form $\alpha = d\varphi - \eta$, and we may consider the contact transformations. These are generated by vector fields, again denoted $X = X_h$, which depend on a hamiltonian $h(x, y, \varphi)$. They are given by

$$X_h = h_x \frac{\partial}{\partial y} - h_y \frac{\partial}{\partial x} - h_\varphi (x \frac{\partial}{\partial x} + y \frac{\partial}{\partial y}) + (-2h + (x h_x + y h_y)) \frac{\partial}{\partial \varphi}.$$

One can check that the above vector field satisfies

$$\mathcal{L}_{X_h} \alpha = -2 h_\varphi \alpha,$$

and hence its flow preserves the contact distribution. To construct a variation generated by such a vector field, we may take its flow $F_t$ in $\mathbb{R}^5$. This flow preserves the contact distribution defined by $\alpha$. If we denote by $\tilde{\ell}$ the local lift of $\ell$ given by $\tilde{\ell} = (\ell, \varphi)$, then we may construct a variation by setting $\ell_t = \Pi \circ F_t \circ \tilde{\ell}$ where $\Pi(x, y, \varphi) = (x, y)$ is the natural projection map. The first variation of the area for such a variation is given by

$$\frac{d}{dt} A(\ell_t)|_{t=0} = \int_\Sigma \operatorname{div}_\Sigma(X_h) \, da = 0 \tag{3.1}$$

where we have $\operatorname{div}_\Sigma(X) = \sum_{i=1}^2 \langle \nabla_{e_i} X, e_i \rangle$ with $e_1, e_2$ being an orthonormal basis for $T\ell(\Sigma)$ at a point and $\nabla$ being the lift under $\Pi$ of the Levi-Civita connection from $\mathbb{R}^4$ to $\mathbb{R}^5$. We have also lifted the metric $g$ by pullback under $\Pi$ to a (degenerate) metric on $\mathbb{R}^5$. In order to use (3.1) effectively, we will first need to do some calculations for the Euclidean metric, and then we will use the fact that $g$ can be made nearly Euclidean near any given point by a suitable choice of coordinates.

We let $\nabla^0$ denote the Euclidean connection, and let $\operatorname{div}_0(X)$ denote the corresponding operator with respect to the Euclidean metric (denoted as the dot product) and connection; thus we have $\operatorname{div}_0(X) = \sum_{i=1}^2 \nabla^0_{e_i} X \cdot e_i$ where $e_1, e_2$ is a Euclidean orthonormal basis. In particular, we see that $\operatorname{div}_0(X)$ does not depend on the $\frac{\partial}{\partial \varphi}$ component of $X$ since the inner product does not involve $\varphi$. On the other hand, the $e_j$ in general have a $\frac{\partial}{\partial \varphi}$ component, so the $\varphi$ dependence of the first 4 components of $X$ is significant. We define an important function $s$ by $s = \frac{1}{2}(x^2 + y^2)$. Since we are interested in the first variation, and it only depends on the first 4 components of a contact vector field, we will use the notation $X \approx Y$ to mean that $X$ and $Y$ have the same first 4 components. The contact vector fields associated with $s, \varphi$ then satisfy

$$X_s \approx x \frac{\partial}{\partial y} - y \frac{\partial}{\partial x} \tag{3.2}$$

and

$$X_\varphi \approx -x \frac{\partial}{\partial x} - y \frac{\partial}{\partial y} \approx J_0 X_s \tag{3.3}$$

where $J_0$ is the standard complex structure in the $x, y$ coordinates. It turns out to be very useful to do a change of variables using the fact that $s$ is positive. We introduce $t + i\theta = \log(s + i\varphi)$ where we choose $\theta \in [-\frac{\pi}{2}, \frac{\pi}{2}]$. In other words, we set $\tilde{s} = \sqrt{s^2 + \varphi^2}$, $t = \log(\tilde{s})$, and $\theta = \tan^{-1}(\varphi/s)$.

The function $\tilde{s}$ will play the role of the square of the distance function. In particular, if we set $r_0 = \sqrt{2}(s^2 + \varphi^2)^{1/4}$ so that $\tilde{s} = \frac{1}{2}r_0^2$, then the function $r_0$ will play the role of the distance to the origin. We define balls $B_\sigma(0)$ to be the points $p$ with $r_0(p) < \sigma$. (We note that the function $r_0$ arises naturally in the function theory of the sub-laplacian on the Heisenberg group; in fact, the fundamental solution of this operator with pole at 0 is a function of $r_0$.) In order to derive monotonicity identities, we consider functions $\eta(t, \theta)$ with compact support, and use the vector field $X_\eta$ in the first variation formula (3.1). One can check the following calculations

$$\mathrm{div}_0(X_\theta) = -2|\nabla^T \theta|^2$$

$$\mathrm{div}_0(X_t) = -2\frac{\sin\theta}{\tilde{s}} - 2\nabla^T \theta \cdot \nabla^T t.$$

Since $X_\eta \approx \eta_t X_t + \eta_\theta X_\theta$, we can compute

$$\begin{aligned}
\mathrm{div}_0 X_{\eta(t,\theta)} &= \eta_{\theta t}(|\nabla^T \theta|^2 - |\nabla^T t|^2) - 2\eta_t \tilde{s}^{-1} \sin\theta - 2\eta_\theta |\nabla^T \theta|^2 \\
&\quad + (\eta_{tt} - \eta_{\theta\theta} - 2\eta_t)\nabla^T \theta \cdot \nabla^T t.
\end{aligned}$$

Using the fact that $|\nabla^T t|^2 + |\nabla^T \theta|^2 = 2\tilde{s}^{-1} \cos\theta$, we may rewrite this

$$\begin{aligned}
\mathrm{div}_0 X_{\eta(t,\theta)} &= (2\eta_{\theta t} - 2\eta_\theta)|\nabla^T \theta|^2 - 2\eta_{\theta t}\tilde{s}^{-1} \cos\theta - 2\eta_t \tilde{s}^{-1} \sin\theta \\
&\quad + (\eta_{tt} - \eta_{\theta\theta} - 2\eta_t)\nabla^T \theta \cdot \nabla^T t.
\end{aligned}$$

We now assume that we have chosen Darboux coordinates centered at a given point $P$ such that the metric $g$ agrees with the Euclidean metric at $P$. Since the connection terms are bounded and the metrics agree at $P$, we may express the first variation condition (3.1) as

$$\int_\Sigma \mathrm{div}_0 X\, da_0 = \mathcal{R}(X) \tag{3.4}$$

where the remainder term satisfies

$$|\mathcal{R}(X)| \leq c \int_\Sigma (\sqrt{\tilde{s}}|\nabla X| + |X|)\, da. \tag{3.5}$$

We will choose a function $\zeta(t)$ below with the properties that it is a smooth nonincreasing function which is identically zero for $t \geq \log(1/2)$ and identically equal to $1 - 2\lambda e^t$ for $t \leq -c$ for positive constants $c$ and $\lambda$ (chosen below). We then define $\psi$ by $\psi(t) = -1/2 e^{-t} \zeta'(t)$, and we see that

$$\int_{-\infty}^{\infty} e^t \psi(t)\, dt = 1/2. \tag{3.6}$$

We thus have $\zeta(t) = 1 - 2\int_{-\infty}^t e^\tau \psi(\tau)\, d\tau$. For a radius $a > 0$, we let $\psi_a(t) = \psi(t - 2\log(a))$ and $\zeta_a(t) = \zeta(t - 2\log(a))$. We then let $\eta(t, \theta)$ be the solution of the wave equation $\eta_{tt} - \eta_{\theta\theta} - 2\eta_t = 0$ with initial data $\eta(t, 0) = 0$, $\frac{\partial\eta}{\partial\theta}(t, 0) = \zeta(t)$. (Note that $\theta$ is the *time* variable in the wave equation and $t$ is the *space* variable.) We then set $\eta_a(t, \theta) = \eta(t - 2\log(a), \theta)$. Since the range on $\theta$ is between $-\pi/2$ and $\pi/2$, we see from domain of dependence considerations that
$$\eta_a(t,\theta) = \theta - 2\lambda \frac{e^t}{a^2}\sin\theta$$
for $t < \log(a^2) - c - \pi/2$ and $\eta_a(t,\theta) = 0$ for $t > \log(a^2/2) + \pi/2$. We observe furthermore that the functions $F_a(t,\theta)$, $G_a(t,\theta)$ given by
$$G = \eta_\theta - \eta_{\theta t},\ G_a(t,\theta) = G(t - 2\log(a), \theta)$$
and
$$F = -1/2(\eta_{\theta t}\tilde{s}^{-1}\cos\theta + \eta_t\tilde{s}^{-1}\sin\theta) = -1/2e^{-t}(\eta_{\theta t}\cos\theta + \eta_t\sin\theta)$$
$$F_a(t,\theta) = F(t - 2\log(a), \theta)$$
satisfy $F = G = 0$ for $t > \log(a^2/2) + \pi/2$ and $F = \lambda$, $G = 1$ for $t < \log(a^2) - c - \pi/2$. We also note that the initial value of $F_a$, given by $F_a(t,0) = \psi_a(t)$, is specified with the constraints above. Thus we have
$$\mathrm{div}_0 X_{\eta_a} = -2G_a|\nabla^T\theta|^2 + \frac{4}{a^2}F_a.$$
We now choose $\eta = \eta_a - \eta_b$ where $0 < b < a$, and consider the first variation of area with vector field $X_\eta$. This gives by (3.4)
$$a^{-2}\int_\Sigma F_a da_0 - b^{-2}\int_\Sigma F_b da_0 = \tfrac{1}{2}\int_\Sigma (G_a - G_b)|\nabla^T\theta|^2 da_0 + \mathcal{R}(X_\eta). \tag{3.7}$$

We now describe the choice of $\zeta$. Let $c$ be a positive number such that $c + \log(1/2) > 2\pi e^{\pi/2}$ and then set $\tau = \tfrac{1}{4}(c + \log(1/2))$. Note that $\tau > \pi/2$. Let $\alpha$ be a nonincreasing, nonnegative function constructed as follows. Define a piecewise linear function $\beta$ by:
$$\beta(t) = \begin{cases} 1 & t < -c \\ -\frac{1}{4\tau}(t - \log(1/2)) & -c \leq t \leq \log(1/2) \\ 0 & t > \log(1/2) \end{cases}$$
Choose $t_0 = -c + 2\tau$ so that $\beta(t_0) = 1/2$. Then define $\alpha$ to be $\beta$ with the corners at $t = -c$ and $t = \log(1/2)$ smoothed and such that:

(i) $0 \leq \alpha \leq 1$.
(ii) $\alpha' \leq 0$.
(iii) $\alpha''(t) \leq 0$ for $t < t_0 + \tau$.
(iv) $\alpha''(t) \geq 0$ for $t > t_0 - \tau$.
(v) $\alpha$ has the symmetry $\alpha(t) = 1 - \alpha(2t_0 - t)$.

Finally, let $\zeta$ be the solution of the ODE $\zeta' - \zeta = -\alpha$ with $\zeta(t) = 0$ for $t > \log(1/2)$. Then it follows that $\zeta$ is a smooth nonincreasing function which is identically zero for $t > \log(1/2)$ and equal to $1 - \lambda e^t$ for $t < -c$ and some positive number $\lambda$.

The following result summarizes the conditions we will need concerning the functions $F$ and $G$.

**Proposition 3.1.** *With the choices above, the function $F$ is nonnegative, and $G$ satisfies $0 \leq G \leq 1$ for $-\pi/2 \leq \theta \leq \pi/2$. Furthermore, there is a fixed number $\theta_0 \in (0,1)$ such that $G_a - G_b \geq 0$ for $0 < b \leq \theta_0 a$*

*Proof.* The transformation $u = e^{-t}\eta$ changes the initial value problem for $\eta$ into the following initial value problem for $u$:

$$\begin{aligned} u_{tt} - u_{\theta\theta} - u &= 0 \\ u(t,0) &= 0 \\ u_\theta(t,0) &= e^{-t}\zeta(t). \end{aligned}$$

Utilizing the Riemann function it is possible to explicitly solve this initial value problem [CH]. Let $J_0(x)$ denote the Bessel function of the first kind of order zero. Then the solution is $u = e^{-t}\eta$ where

$$\eta(t,\theta) = \tfrac{1}{2} \int_{-\theta}^{\theta} J_0(\sqrt{(\theta-\mu)(\theta+\mu)})e^{-\mu}\zeta(\mu+t)d\mu \tag{3.8}$$

where $J_0(\sigma)$ is the solution of the ODE

$$J_0'' + \frac{1}{\sigma}J_0' + J_0 = 0, \ J_0(0) = 1. \tag{3.9}$$

The basic properties of $J_0$ which we shall need are that $J_0(\sigma) > 0$, $J_0'(\sigma) < 0$, and $J_0''(\sigma) < 0$ for $0 \leq \sigma \leq \pi/2$. The fact that $J_0''$ and $J_0'$ are negative up to the first zero of $J_0$ is shown in [CH, p. 495]. It follows that $J_0'' + J_0 > 0$ for $0 < \sigma < \sigma_0$ where $\sigma_0$ is the first zero of $J_0$. It remains to show that $\sigma_0 \geq \pi/2$. To see this, we set $v_\varepsilon(\sigma) = J_0(\sigma)/\cos(\varepsilon + \sigma)$ for $0 \leq \sigma < \pi/2 - \varepsilon$, and we compute $J_0''(\sigma) = v_\varepsilon''(\sigma)\cos(\varepsilon+\sigma) - 2v_\varepsilon'(\sigma)\sin(\varepsilon+\sigma) - J_0(\sigma)$. It follows that $v_\varepsilon''(\sigma)\cos(\varepsilon+\sigma) > 2v_\varepsilon'(\sigma)\sin(\varepsilon+\sigma)$ for $0 \leq \sigma < \min\{\sigma_0, \pi/2 - \varepsilon\}$. Thus the function $v_\varepsilon$ can not achieve a local maximum in $(0, \min\{\sigma_0, \pi/2 - \varepsilon\})$. On the other hand we have $v_\varepsilon'(0) = \sin(\varepsilon)/\cos^2(\varepsilon) > 0$ for $\varepsilon > 0$, so it follows that $v_\varepsilon'(\sigma) > 0$ for $\sigma$ small. Therefore we see that $v_\varepsilon$ is an increasing function for $0 \leq \sigma < \min\{\sigma_0, \pi/2 - \varepsilon\}$ for any $\varepsilon > 0$. In particular it follows that $\sigma_0 \geq \pi/2$ and setting $v = \lim_{\varepsilon \to 0} v_\varepsilon$ we have $v(\sigma) \geq \cos(\sigma)$ on this interval.

To verify the proposition, we first need only check that $F$ is nonnegative. We see from (3.8) that

$$\eta_t(t,\theta) = \frac{1}{2}\int_{-\theta}^{\theta} J_0(\sqrt{\theta^2 - \mu^2})e^{-\mu}\zeta'(\mu+t)d\mu.$$

Differentiating in $\theta$ we have

$$\begin{aligned} \eta_{t\theta} &= \tfrac{1}{2}\{e^{-\theta}\zeta'(t+\theta) + e^{\theta}\zeta'(t-\theta) \\ &+ \int_{-\theta}^{\theta} \frac{\partial}{\partial\theta}J_0(\sqrt{\theta^2-\mu^2})e^{-\mu}\zeta'(\mu+t)d\mu\}. \end{aligned}$$

Substituting $\zeta' = -2e^t\psi$, we then have the explicit formula for $F$

$$\begin{aligned} F(t,\theta) &= \tfrac{1}{2}(\psi(t+\theta) + \psi(t-\theta))\cos(\theta) + \\ &\tfrac{1}{2}\int_{-\theta}^{\theta}\{\sin(\theta)J_0(\sqrt{\theta^2-\mu^2})) + \cos(\theta)\frac{\partial}{\partial\theta}J_0(\sqrt{\theta^2-\mu^2})\}\psi(\mu+t)d\mu. \end{aligned}$$

Now we observe that the integrand above is nonnegative for $0 < \theta < \pi/2$. To see this, note that it suffices to show that the function $Q(\theta) = J_0(\sqrt{\theta^2 - \mu^2})/\cos(\theta)$ is increasing in this range. Recall that we showed above that the function $v(\sigma) = J_0(\sigma)/\cos(\sigma)$ is increasing for $0 < \sigma < \pi/2$, and note that we have

$$Q(\theta) = v(\sqrt{\theta^2 - \mu^2}) \frac{\cos(\sqrt{\theta^2 - \mu^2})}{\cos(\theta)},$$

and it is elementary to check that the ratio of cosine functions is increasing in $\theta$ for $0 < \theta < \pi/2$ and $|\mu| \leq \theta$. We have therefore shown that $F \geq 0$.

To analyze $G$ we begin by showing that $G(t, \theta) \geq 0$. We first observe that $u = e^{-t}G$ is a solution of the wave equation $u_{tt} - u_{\theta\theta} - u = 0$ with initial data given by $u(t, 0) = u_0(t) = e^{-t}(\zeta - \zeta') = e^{-t}\alpha$ and $u_\theta(t, 0) = 0$. We then have the explicit formula

$$u(t, \theta) = \frac{1}{2}\{u_0(t+\theta) + u_0(t-\theta) + \int_{-\theta}^{\theta} \frac{\partial}{\partial \theta} J_0(\sqrt{\theta^2 - \mu^2}) u_0(t+\mu) \, d\mu\}. \tag{3.10}$$

We note that we have

$$\cos(\theta) = 1 + \frac{1}{2}\int_{-\theta}^{\theta} \frac{\partial}{\partial \theta} J_0(\sqrt{\theta^2 - \mu^2}) \, d\mu$$

since both sides represent the solution of the Cauchy problem $v_{tt} - v_{\theta\theta} - v = 0$, $v(t, 0) = 1$, $v_\theta(t, 0) = 0$. Since $u_0 = e^{-t}\alpha$ and $\alpha'' \geq 0$ for $t \geq t_0 - \tau$, we have $u_0'' = e^{-t}(\alpha'' - 2\alpha' + \alpha) \geq 0$ for $t \geq t_0 - \tau$. It follows that the function $u_0$ is a convex function for $t \geq t_0 - \tau$. Suppose $t \geq t_0$ and $|\theta| \leq \pi/2$. If $\ell(t+\mu)$ denotes the linear function with $\ell(t+\theta) = u_0(t+\theta)$ and $\ell(t-\theta) = u_0(t-\theta)$, then we have $u_0(t+\mu) \leq \ell(t+\mu)$ for $|\mu| \leq |\theta|$. (Recall that $\tau > \pi/2$.) Hence, for $t \geq t_0$ and $|\theta| \leq \pi/2$, we have

$$G(t, \theta) \geq \frac{1}{2}e^t\{u_0(t+\theta) + u_0(t-\theta) + \int_{-\theta}^{\theta} \frac{\partial}{\partial \theta} J_0(\sqrt{\theta^2 - \mu^2}) \ell(t+\mu) \, d\mu\}.$$

Since $\frac{\partial}{\partial \theta} J_0(\sqrt{\theta^2 - \mu^2})$ is an even function of $\mu$ whose integral is greater than or equal to $-2$, it follows that the right hand side is nonnegative, and so $G(t, \theta) \geq 0$ for $t \geq t_0$. We next consider the case $t \leq t_0$. From (3.10), using that $\alpha$ is nonincreasing and $\frac{\partial}{\partial \theta} J_0(\sqrt{\theta^2 - \mu^2}) \leq 0$, we have:

$$G(t, \theta) \geq \frac{1}{2}\{e^{-\theta}\alpha(t+\theta) + e^{\theta}\alpha(t-\theta) + \alpha(t-\theta)\int_{-\theta}^{\theta} \frac{\partial}{\partial \theta} J_0(\sqrt{\theta^2 - \mu^2})e^{-\mu} \, d\mu\}, \tag{3.11}$$

Note that $u(t, \theta) = e^{-t}$ is the unique solution of the Cauchy problem: $u_{tt} - u_{\theta\theta} - u = 0$, with initial data given by $u(t, 0) = e^{-t}$ and $u_\theta(t, 0) = 0$. From the explicit formula (3.10) for the solution with this initial data we conclude:

$$1 = \frac{1}{2}\{e^{-\theta} + e^{\theta} + \int_{-\theta}^{\theta} \frac{\partial}{\partial \theta} J_0(\sqrt{\theta^2 - \mu^2})e^{-\mu} \, d\mu\}.$$

Using this in (3.11) we have:
$$G(t,\theta) \geq \alpha(t-\theta) + \frac{1}{2}e^{-\theta}(\alpha(t+\theta) - \alpha(t-\theta)).$$

If $t \leq t_0$ and $|\theta| \leq \pi/2$ this implies:
$$\begin{aligned} G(t,\theta) &\geq \alpha(t+\pi/2) - \frac{1}{2}e^{\pi/2}\frac{\pi}{c+\log(1/2)} \\ &\geq \alpha(t+\pi/2) - \frac{1}{4} > 0, \end{aligned}$$

since $\alpha(t+\pi/2) > \alpha(t+\tau) \geq \frac{1}{4}$. Hence we conclude that $G(t,\theta) \geq 0$ for all $t$ and $|\theta| \leq \pi/2$.

To see that $G \leq 1$, we observe that the function $v = e^{-t}(1-G)$ satisfies the same wave equation as $u$, and has initial conditions $v(t,0) = v_0(t) = e^{-t}(1-\alpha(t)) = e^{-t}\alpha(2t_0 - t)$ and $v_\theta(t,0) = 0$. We can then apply the previous argument to conclude that $(1 - G(t,\theta)) \geq 0$. □

We have shown the (approximate) monotonicity of the quantity $a^{-2}\int_\Sigma F_a$ as a function of the radius $a$. Observe that because of the normalization (3.6), we have $\lim_{a\to 0}(\pi a^2)^{-1}\int_\Sigma F_a = 1$ at a smooth point of $\Sigma$.

**Proposition 3.2.** *Let $\ell : D \to N$ be a weakly lagrangian, weakly conformal, contact stationary map having image in a Darboux chart (centered at $\ell(0)$). There exist positive constants $c_1, c_2$ depending only on $N$ and the area of $\ell$ such that for any $\sigma \in (0,1)$ we have*
$$c_1 \leq \sigma^{-2}A(\ell(\Sigma) \cap B_\sigma(0)) \leq c_2.$$

*Proof.* We choose $a = 1$ in (3.7), and observe that by (3.5) we have
$$|\mathcal{R}(X_\eta)| \leq |\mathcal{R}(X_{\eta_1})| + |\mathcal{R}(X_{\eta_b})| \leq c + cb^{-1}A(\ell(\Sigma) \cap B_{\lambda b})$$
where $\lambda > 1$ is fixed. From (3.7), we then have for any $b \in (0, \theta_0)$
$$b^{-2}A(\lambda^{-1}b) \leq c + cb^{-1}A(\lambda b).$$

This implies that
$$\sup_{t\in(0,s)} t^{-2}A(t) \leq c + cs \sup_{t\in(0,\lambda^2 s)} t^{-2}A(t)$$

for any $s \in (0, \lambda^{-1}\theta_0)$. Fixing $s$ such that $cs < 1/2$, we then get the bound
$$\sup_{t\in(0,s)} t^{-2}A(t) \leq c + \sup_{t\in(s,1)} t^{-2}A(t) \leq c$$

which is the desired upper bound. From the upper bound and (3.7) it follows that for each $P \in \Sigma$ and $\sigma > 0$ we have
$$\sup_{t\in(0,\theta_0\sigma)} (\pi t^2)^{-1}\int_\Sigma F_t da \leq (\pi\sigma^2)^{-1}\int_\Sigma F_\sigma da + c\sigma.$$

It follows that the limit
$$\Theta(P) = \lim_{\sigma \to 0} (\pi\sigma^2)^{-1} \int_\Sigma F_\sigma \, da$$
exists and is upper semicontinuous. At any Lebesgue point $P$ of the derivative of $\ell$ at which the derivative is nonzero, we have $\Theta(\ell(P)) \geq 1$. Since such points have full measure with respect to the induced measure $da$, we have $\Theta \geq 1$ at all points of the support of $\Sigma$. The lower bound now follows. $\square$

The above result implies the following important continuity estimate which is an improvement on Theorem 2.8. In the following results, the notation $\tilde{\ell}$ is used to denote the legendrian lift of the exact lagrangian map $\ell$. Note that $A(\ell) = A(\tilde{\ell})$, but that the balls we consider lie in $\mathbb{R}^5$, and hence the monotonicity formula exists there also.

**Proposition 3.3.** *Suppose that $\ell$ is a weakly conformal, exact, lagrangian stationary map from the unit disk $D_1$ into $N$. Assume further that there exists $r_0 > 0$ and a constant $c$ such that for $r \leq r_0$ we have $A(\tilde{\ell}(D_1)) \cap B_r(P) \leq cr^2$ for any $P \in N$. Then $\tilde{\ell}$ is Hölder continuous in $D_{1/2}$. Moreover, there is a constant $\varepsilon_0$ depending only on $N$ such that if $A(\ell(D_1)) \leq \varepsilon_0$, then there is a uniform bound on the Hölder modulus of continuity of $\tilde{\ell}$ in $D_{1/2}$.*

*Proof.* It suffices to prove the last statement, since we may cover $D_{1/2}$ with small balls of small area and apply the last statement to a rescaled version of $\ell$ to prove the Hölder continuity. To prove the last statement, we may rescale the map in the image so that $A(\tilde{\ell}(D_1)) = 1$, and this rescaling only flattens the metric on the image. We now claim that there is a constant $\theta < 1$ such that $A(\tilde{\ell}(D_{1/2})) \leq \theta$. To prove this, we need only show that there is a positive constant $\delta$ such that $A(\ell(D_1 \setminus D_{1/2})) \geq \delta$ for then our conclusion follows with $\theta = 1 - \delta$. Suppose to the contrary that $A(\ell(D_1 \setminus D_{1/2}))$ is small. We can then find a radius $\sigma \in [1/2, 1]$ such that the length $L(\ell(\partial D_\sigma))$ is small, say less than $\varepsilon$. Let $P$ be a point of the curve $\ell(\partial D_\sigma)$, and choose Darboux coordinates centered at $P$ in which the metric is Euclidean at the origin. Construct a smooth function $\zeta(r)$ such that $\zeta(r) = 0$ for $r \leq \varepsilon$, $\zeta(r) = 1$ for $r \geq 2\varepsilon$, and $\varepsilon|\zeta'(r)| + \varepsilon^2|\zeta''(r)| \leq c$. We then consider the hamiltonian $h = -\varphi \cdot \zeta \circ r_0$. The corresponding contact vector field satisfies $X \approx \zeta \circ r_0 (x \frac{\partial}{\partial x} + y \frac{\partial}{\partial y}) - \varphi X_{\zeta \circ r_0}$. Using this vector field in the first variation formula and using easy bounds we get:

$$\begin{aligned}
A(\tilde{\ell}(D_\sigma) \setminus B_{2\varepsilon}) &\leq c \int_{D_\sigma \cap \tilde{\ell}^{-1}(B_{2\varepsilon} \setminus B_\varepsilon)} (|D\varphi||D\zeta \circ r_0| + |\varphi||\operatorname{div}_\Sigma(X_{\zeta \circ r_0})|) \, d\mu \\
&\leq c \int_{D_\sigma \cap \tilde{\ell}^{-1}(B_{2\varepsilon} \setminus B_\varepsilon)} (\varepsilon \cdot \varepsilon^{-1} + \varepsilon^2 \cdot \varepsilon^{-2}) \, d\mu \\
&\leq cA(\tilde{\ell}(D_\sigma \cap B_{2\varepsilon})) \leq c\varepsilon^2.
\end{aligned}$$

This implies $A(\tilde{\ell}(D_\sigma)) \leq c\varepsilon^2$ which is a contradiction if $\varepsilon$ is small. This verifies our claim that $A(\tilde{\ell}(D_{1/2})) \leq \theta$.

We complete the proof by observing that we can iterate this inequality to obtain $A(\tilde{\ell}(D_{2^{-n}})) \leq \theta^n$ which implies that $A(\tilde{\ell}(D_r)) \leq cr^{2\alpha}$ for some $\alpha > 0$. This same inequality then holds for disks with any center point in $D_{1/2}$, and by Morrey's lemma (since $A(\tilde{\ell})$ is the Dirichlet integral of $\ell$) this implies that $\ell$ has a uniform Hölder estimate on $D_{1/2}$. Since $|D\varphi| \leq c$, and the Dirichlet integral is conformally invariant, we see that $\int_{D_r} |\nabla \varphi|^2 \, dt \leq cr^{2\alpha}$, and $\varphi$ is therefore also Hölder continuous. $\square$

We will need the following global Hölder estimate under the condition that the restriction of $\ell$ to $\partial D_1$ has finite energy.

**Proposition 3.4.** *Suppose that $\ell$ is a weakly conformal, exact, lagrangian stationary map from the unit disk $D_1$ into $N$. Assume that there exists $r_0 > 0$ and a constant $c$ such that for $r \leq r_0$ we have $A(\tilde{\ell}(D_1)) \cap B_r(P) \leq cr^2$ for any $P \in N$. Finally assume that the trace of $\ell$ as a map from the unit circle to $N$ has energy bounded by a constant $c$. Then $\tilde{\ell}$ is Hölder continuous on all of $D_1$. Moreover, there is a constant $\varepsilon_0$ depending only on $N$ such that if $A(\ell(D_1)) \leq \varepsilon_0$, then there is a uniform bound on the Hölder modulus of continuity of $\tilde{\ell}$ on $D_1$.*

*Proof.* The proof is similar to the previous one. We show that
$$\int_{D_r(t) \cap D_1} (|\nabla \ell|^2 + |\nabla \varphi|^2) \, dt \leq cr^{2\alpha}$$
for all $t \in D_1$. If a disk $D_r(t)$ intersects $\partial D_1$ in a point $t_0$, we may instead consider the larger disk $D_{2r}(t_0)$ and prove the bound there. Thus it suffices to consider disks which are either completely contained in $D_1$, or centered on $\partial D_1$. The first case is done above. To handle disks centered on $\partial D_1$, we need to use the fact that $\ell$ has a $C^{0,1/2}$ estimate on $\partial D_1$ which follows from the bound on the boundary energy. With the aid of this bound, the proof proceeds as above. $\square$

We close this section by analyzing the situation for which the monotonicity inequality is an equality.

**Proposition 3.5.** *Let $\ell$ be a weakly conformal, exact, lagrangian stationary map from $\Omega \subset R^2$ into $R^4$ with euclidean metrics. Suppose that $\ell(0) = 0$, that $\ell^{-1}(B_a)$ is compact and $a^{-2} \int_\Omega F_a \circ \ell \, da$ is constant for $a \leq a_0$. It follows that $C = \tilde{\ell}(\Omega) \cap B_{a_0}$ is contained in $\{\varphi = 0\}$, and is a cone over $0$ in the sense that for almost every $t \in \ell^{-1}(C)$, the position vector $(x, y) = \ell(t) \in C$, is tangent to $C$. If the map $\ell$ is a smooth immersion in a neighborhood $U$ of $t \in \ell^{-1}(C)$, then $\ell(U)$ is a portion of a cone.*

*Proof.* From (3.7), we see that $\theta$ is a constant function if equality holds in the monotonicity inequality. This means that $\varphi = as$ for some constant $a$. We want to show that $a = 0$, so we note that $r_0^2 = 2(1+a^2)s$, so we choose a hamiltonian function $\zeta \circ s$ where $\zeta$ is a monotone increasing, concave function with $\zeta(0) = 0$ and $\zeta = 1$ outside $B_{a_0}$. Using the first variation formula and the fact that $div(X_s) = 0$, we find $\int \zeta''(s) Ds \cdot X_s \, d\mu = 0$. If $a \neq 0$, this implies that $\int \zeta''(s) D\varphi \cdot X_s \, d\mu = 0$. On the other hand, $D\varphi \cdot X_s = |D\varphi|^2$, so since $\zeta'' > 0$ we see that $\varphi$ is a constant function. Since $\varphi(0) = 0$, we have shown that $\varphi$ is identically zero and hence $a = 0$ as claimed. The condition that $D\varphi = 0$ implies that the position vector is tangent to the image of $\ell$, and this condition easily implies that any regular piece of the image is a cone over 0. □

## 4. Regularity for The Two Dimensional Mapping Problem

In this section we discuss the regularity properties of maps $\ell$ from a surface to a symplectic 4-manifold $(N, \omega)$ with almost complex structure $J$ and compatible Riemannian metric $g$. Since the theorem is local, we assume that $\ell$ is defined on the unit disk $D$ in the plane with coordinates $(t^1, t^2)$. We use the notation $D_r(p)$ to denote the disk of radius $r$ centered at $p \in \mathbb{R}^2$ with $D_r$ denoting the ball with center at the origin. We make the following assumptions on the map where we use the notation $\ell_i$ to denote the weak partial derivative $\ell_{t^i}$:

i) The map $\ell$ is a $W^{1,2}$ map from $D$ to $N$ which is weakly conformal and weakly lagrangian in the sense that for almost every $(t^1, t^2) \in D$ we have $\|\ell_1\|^2 = \|\ell_2\|^2$, $\langle \ell_1, \ell_2 \rangle = 0$, and $\langle J\ell_1, \ell_2 \rangle = 0$.

ii) The map $\ell$ is exact on $D$.

iii) The restriction of $\ell$ to any subdisk is *stationary* as a contact map in the sense that for any $P \in D$ and $r \in (0, 1 - |P|]$, and any smooth contact vector field $X$ which vanishes in a neighborhood of $\ell(\partial D_r(P))$ we have

$$\frac{d}{dt} A(\ell_t(D_r(P)))|_{t=0} = 0 \tag{4.1}$$

where $\ell_t = F_t \circ \ell$ where $F_t$ is the flow generated by $X$.

Note that the restriction of $\ell$ to $\partial D_r(P)$ is a $W^{1,2}$ map for almost every $r \in (0, 1-|P|]$, and hence this restriction is continuous, and the length of the curve $\ell(\partial D_r)$ is finite. It is thus clear that condition (iii) allows for many variations. Using (ii), (iii) and monotonicity it follows (Theorem 2.8 and Proposition 3.3) that $\ell$ is continuous. We can thus suppose that for small $r$ the image of $\ell$ lies in a Darboux neighborhood $U$. Let $(x, y)$ be Darboux coordinates in $U$ centered at a point with coordinates $(0, 0)$ such that $g, J$ are standard at $(0, 0)$. Let $\tilde{U}$ denote the contact lifting of $U$. Then $\ell$ lifts to a legendrian map (again denoted $\ell$) from $D_r \to \tilde{U}$. With respect to the contact coordinates on $\tilde{U}$ we can write $\ell = (x, y, \varphi)$. On $D_r$ we have:

$$d\varphi = \sum_i (x_i dy_i - y_i dx_i).$$

This expresses the exactness of $\ell$ on $D_r$. We state the following regularity result.

**Theorem 4.1.** *There exists $\varepsilon_0 > 0$ such that if there is a conformal linear lagrangian map $\ell_0$ with $\ell_0(0) = 0$ and an $r < 1/2$ for which:*

$$r^{-2} \int_{D_r} |\nabla(\ell - \ell_0)|^2 \, dt + r^{-4} \int_{D_r} |\ell - \ell_0|^2 dt < \varepsilon_0,$$

$$r^{-2} \sup_{t \in D_r} |\ell_0(t) - (t, 0)|^2 < \varepsilon_0,$$

$$\|g - \delta\|_{C^2(B_1)} < \varepsilon_0$$

*then $u$ is a smooth embedding on $D_{r/2}$.*

The following result is necessary for controlling the parametrization of maps constructed as rescaled limits of $\ell$.

**Proposition 4.2.** *Assume that there is a constant $c_1$ such that $A(\ell(D_1)) \leq c_1$, and $\inf_{t \in D_{1/4}} A(\ell(D_{1/8}(t))) \geq c_1^{-1} A(\ell(D_1))$. There is a constant $c$ depending only on $c_1$ such that for all $t \in D_{1/4}$ and all $r \leq 1/4$ we have*

$$A(\ell(D_r(t))) \geq c^{-1} A(\ell(D_{2r}(t))).$$

*Proof.* Fix a point $t_0 \in D_{1/4}$, and let $A(r) = A(\ell(D_r(t_0)))$ which is a continuous, monotone increasing function. The derivative $A'(r)$ then exists for almost every $r$ and is equal to $\int_{\partial D_r(t_0)} |\nabla \ell|^2 \, ds$. By differentiation theory of monotone functions we have

$$\int_\sigma^{2\sigma} \frac{A'(r)}{A(r)} \, dr \leq \log\left(\frac{A(2\sigma)}{A(\sigma)}\right). \tag{4.2}$$

Now if we choose $r \in (0, 1/4]$, we may find a unique nonnegative integer $k$ such that $r = 2^{-k} r_0$ with $r_0 \in (1/8, 1/4]$. We let $r_j = 2^{-j} r_0$ for $j = 0, 1, \ldots, k$ so that $r = r_k$. We will prove that $A(r_j) \geq c^{-1} A(2r_j)$ for a fixed constant $c$ by induction on $j$. We begin with $j = 0$, and observe that the inequality follows with $c = c_1$ by hypothesis since $1/8 < r_0$ and $2r_0 \leq 1/2$ so that $D_{2r_0}(t) \subset D_1$. We now prove the inductive step. Assume that $A(r_{j-1}) \geq c^{-1} A(2r_{j-1})$ for some $j \geq 1$. Applying the mean value theorem with $\sigma = r_{j-1}$ in (4.2) and using the inductive assumption, we find $\sigma_0 \in [1, 2]$ with $r_{j-1} A'(\sigma_0 r_{j-1}) \leq \log(c) A(\sigma_0 r_{j-1})$. We now rescale the map $\ell$ by setting $\bar{\ell}(\tau) = (A(\sigma_0 r_{j-1}))^{-1/2} \ell(\sigma_0 r_{j-1} \tau + t_0)$. We then define $\bar{A}(r)$ by

$$\bar{A}(r) = A(\bar{\ell}(D_r)) = \frac{A(\sigma_0 r_{j-1} r)}{A(\sigma_0 r_{j-1})}.$$

Our conditions then imply $\bar{A}(1) = 1$ and $\bar{A}'(1) \leq \log(c)$, and we need to show that $\bar{A}(2^{-1} \sigma_0^{-1}) \geq c^{-1} \bar{A}(\sigma_0^{-1})$. Because $\sigma_0 \in [1, 2]$, it suffices to show that $\bar{A}(1/4) \geq c^{-1}$. By the mean value inequality we may choose $\sigma_1 \in [1/8, 1/4]$ such that $\bar{A}'(\sigma_1) \leq 8 \bar{A}(1/4)$. If we set $\delta = (\bar{A}(1/4))^{1/2}$, we may apply the Schwarz inequality to show that

$$L(\bar{\ell}(\partial D_{\sigma_1})) \leq (2\pi \sigma_1)^{1/2} (\bar{A}'(\sigma_1))^{1/2} \leq 2\pi^{1/2} \delta.$$

Let $\rho(P)$ for $P \in N$ be the (modified) distance function from a point $P_0 \in \bar{\ell}(\partial D_{\sigma_1})$ defining the balls $B_r(P_0)$ (as in Section 3). We may also estimate $\varphi$ by integration along $\bar{\ell}(\partial D_{\sigma_1})$. (Note that $\varphi(P_0) = 0$.) We find that $\bar{\ell}(\partial D_{\sigma_1}) \subset B_{\delta_1}(P_0)$ where $\delta_1 = c_1 \delta$ for a fixed constant $c_1$. Since $\bar{A}(1) = 1$, we may use the monotonicity formula to show that $\bar{A}(B_{\rho_0}(P_0)) \leq c_2 \rho_0^2$, and thus there is a fixed constant $\rho_0 > 0$ depending only on $c_1$ such that $\sup_{D_1} \rho \circ \bar{\ell} \geq 2\rho_0$. By the global Hölder estimate of Proposition 3.4 the map $\log(c)^{-1/2} \bar{\ell}$ satisfies a uniform Hölder estimate, and thus $\bar{\ell}$ satisfies the bound $|\bar{\ell}(t_1) - \bar{\ell}(t_2)| \leq \log(c)^{1/2} c_3 |t_1 - t_2|^\alpha$ for fixed positive constants $c_3, \alpha$. Using this estimate, we may find a point $t_1 \in D_1$ such that $\rho(\bar{\ell}(t)) \geq \rho_0$ for all $t \in D_1 \cap D_{\delta_2}(t_1)$ where $\delta_2 = c_4 \log(c)^{-1/(2\alpha)}$. We then consider the following function $\zeta$

$$\zeta(P) = \frac{\log(\rho(P)/\delta_1)}{\log(\rho_0/\delta_1)}$$

for $\delta_1 \leq \rho(P) \leq \rho_0$, and $\zeta(P) = 0$ for $\rho(P) < \delta_1$, $\zeta(P) = 1$ for $\rho(P) > \rho_0$. Using monotonicity and the coarea formula we can derive the bound $\int_{D_1} |\nabla \zeta \circ \bar{\ell}|^2 \, dt \leq c_5/\log(1/\delta_1)$. On the other hand, if we let $r, \theta$ be polar coordinates centered at $t_1$, we may use the fact that $\zeta \circ \bar{\ell}(t) = 0$ for $t \in D_{\sigma_1}$ and $\zeta \circ \bar{\ell}(t) = 1$ for $t \in D_1 \cap D_{\delta_2}(t_1)$ to conclude that $\int_{D_1} |\nabla \zeta \circ \bar{\ell}|^2 \, dt \geq c_6/\log(1/\delta_2)$. (This may be seen by using the Schwarz inequality and the fact that $\zeta \circ \bar{\ell}$ changes from 0 to 1 along a fixed fraction of the rays from $t_1$ as follows:

$$1 \leq \left(\int_{\delta_2}^1 |\nabla \zeta \circ \bar{\ell}| \, dr\right)^2 \leq \left(\int_{\delta_2}^1 r^{-1} \, dr\right)\left(\int_{\delta_2}^1 |\nabla \zeta \circ \bar{\ell}|^2 r \, dr\right).$$

The bound then follows by integrating along rays.) Combining these inequalities we conclude that $\bar{A}(1/4) \geq c_7 (\log(c))^{-\beta}$ for constants $c_7, \beta$ depending only on $c_1$. We then observe that if $c$ is chosen larger than a fixed constant, we have $c_7 (\log(c))^{-\beta} \geq c^{-1}$. This completes the proof. □

The next result deals with surfaces which lie near a plane, and shows that such surfaces can be well approximated by smooth conformal maps to the plane. We use the notation $\hat{A}(\ell(\Sigma) \cap B_r(P))$ to denote the monotone quantity $\int_\Sigma F_r \, da$ from Section 3.

**Proposition 4.3.** *Assume that $\ell = (x(t), y(t), \varphi(t)) : D_2 \to \mathbb{R}^5$ with $\ell(0) = 0$ is a stationary legendrian map with $A(\ell(D_2)) = 1$, $A(\ell(D_1)) \geq c_1^{-1}$, and $\int_{D_2} |\nabla y|^2 \, dt \leq \varepsilon$. Assume further that $\hat{A}(\ell(D_2) \cap B_r(P)) \leq (1 + \varepsilon)\pi r^2$ for any $P \in \mathbb{R}^5$. Given any $\delta > 0$, there is a constant $\varepsilon_2 > 0$ depending only on $\delta, c_1$ such that if $\varepsilon \leq \varepsilon_2$, then there is a 1-1 holomorphic or anti-holomorphic map $x_0 : D_1 \to \mathbb{R}^2$ such that $\int_{D_1} (|\nabla(\ell - x_0)|^2 + |\ell - x_0|^2) \, dt \leq \delta$.*

*Proof.* To prove this result by contradiction, it suffices to show that if $\ell^{(j)}$ is a sequence of maps satisfying the hypotheses with $\varepsilon_j \to 0$, then a subsequence of $\ell^{(j)}$ converges strongly in $W^{1,2}(D_1)$ to an injective conformal map to the $x$-plane. By the uniform continuity estimate, we may choose a subsequence, again denoted $\ell^{(j)}$ which converges uniformly to a continuous $W^{1,2}$ map $x_0$ on $D_{3/2}$. We first show that $x_0$ is 1-1 in $D_1$. If this were not true, then we can find a point $\xi$ and points $t_1, t_2 \in D_1$ with $t_1 \neq t_2$ and $x_0(t_1) = x_0(t_2) = \xi$. We can show that there is a $\delta > 0$ such that for $j$ sufficiently large $(\ell^{(j)})^{-1}(B_\delta(\xi))$ has distinct connected components $U_1, U_2$ containing $t_1, t_2$ respectively. This follows by the same type of capacity argument using the function $\zeta$ as in the previous proof, the point being that each connected component of $(\ell^{(j)})^{-1}(B_\delta(\xi))$ must have arbitrarily small diameter if $\delta$ is chosen sufficiently small. Since $\ell^{(j)}(t_i) \to \xi$, it follows from monotonicity that each component $U_i$, $i = 1, 2$ has area at least $\pi\delta^2 - o(j)$. For $j$ sufficiently large, this contradicts the assumption that $\hat{A}(\ell(D_2) \cap B_\delta(\xi)) \leq (1 + \varepsilon_j)\pi\delta^2$. Therefore we have shown that $x_0$ is 1-1 on $D_1$. It follows that $\Omega = x_0(D_1)$ is an open set bounded by the Jordan curve $x_0(\partial D_1)$.

We now show that $\ell^{(j)}$ converges strongly to $x_0$ in $W^{1,2}(D_1)$. We first show that $A(\ell^{(j)}(D_1)) \to A(x_0(D_1)) = A(\Omega)$. To see this, first note that if $B_r(\xi)$ is any open ball centered at a point $\xi \in \Omega$ with $\partial\Omega \cap B_r(\xi) = \phi$, then we have $A(\ell^{(j)}(D_1) \cap B_r(\xi)) \to A(\Omega \cap B_r(\xi)) = \pi r^2$. By the Vitali covering lemma, we can cover almost all of $\Omega$ with a disjoint collection of such balls. If $\mathcal{O}$ denotes the union of this collection of balls, we then have $A(\ell^{(j)}(D_1) \cap \mathcal{O}) \to A(\Omega \cap \mathcal{O}) = A(\Omega)$. Since $\bar{\Omega} \setminus \mathcal{O}$ may be covered by a finite number of balls $B_{r_k}(\xi_k)$ with $\sum r_k^2$ arbitrarily small, it follows that $\ell^{(j)}(D_1) \setminus \mathcal{O}$ will be contained in this finite union for $j$ sufficiently large, and therefore its area will be arbitrarily small. Thus we have shown that $A(\ell^{(j)}(D_1)) \to A(x_0(D_1))$ as required. By lowersemicontinuity of the energy we have the energy of $x_0$ on $D_1$ less than or equal to the liminf of the energies of the $\ell^{(j)}$. Since the $\ell^{(j)}$ are weakly conformal their energies agree with thier areas. Thus we have the energy of $x_0$ on $D_1$ is less than or equal to its area. This implies that the two are equal and that $x_0$ is weakly conformal. It also implies that there is no energy loss in the convergence of $\ell^{(j)}$ to $x_0$, and therefore this convergence is strong in $W^{1,2}(D_1)$. $\square$

We have the following important consequence of these propositions.

**Proposition 4.4.** *Assume that $A(\ell(D_2)) \leq c_1$, and that $\ell_0$ is a linear holomorphic map to a lagrangian plane. Assume further that*

$$\int_{D_2} \{|\nabla(\ell - \ell_0)|^2 + |\ell - \ell_0|^2\} \, dt \leq \delta.$$

*Given any $\varepsilon > 0$, there is a $\delta > 0$ depending only on $\varepsilon$ and $c_1$, and a fixed constant $c$ such that if the above inequality holds, then*

$$\int_{D_1} |\nabla(\ell - \ell_0)|^2 \, dt \leq c \int_{D_2} |\ell - \ell_0|^2 \, dt + \varepsilon \int_{D_2} |\nabla(\ell - \ell_0)|^2 \, dt.$$

*Proof.* We may assume by change of coordinates and dilation that in Darboux coordinates $(x, y)$ the map $\ell_0(t) = (t, 0)$, and we write $\ell(t) = (x(t), y(t))$. We begin by estimating the Dirichlet integral of $y(t)$. We observe that the hamiltonian $x \cdot y - \varphi$ yields the vector field $y\frac{\partial}{\partial y}$. We let $\zeta(x)$ be a cutoff function with $\zeta(x) = 1$ for $|x| \leq 3/2$, and $\zeta(x) = 0$ for $|x| \geq 7/4$. We then take the hamiltonian $h = \zeta(x)(x \cdot y - \varphi)$. In proving the desired inequality, we may assume that $\int_{D_2} |\ell - \ell_0|^2 \, dt$ is small, for otherwise the inequality is trivial. Thus it follows that $\ell$ is uniformly close to $\ell_0$ interior to $D_2$. Thus the hamiltonian $h \circ \ell$ has compact support in $D_2$. The vector field determined by $h$ is $\zeta(x)y\frac{\partial}{\partial y} + (x \cdot y - \varphi)\zeta_x\frac{\partial}{\partial y}$. Applying the first variation formula we get

$$\int_{D_2} \langle \nabla y, \nabla(\zeta(x(t))y(t) + (x(t) \cdot y(t) - \varphi)\zeta_x(x(t))) \rangle \, dt = 0$$

where $\nabla$ refers to derivatives with respect to $t$. This equation easily implies

$$\int_{D_2} \zeta(x(t))|\nabla y|^2 \, dt \leq c \int_{D_2} (|\zeta_x| + |\zeta_{xx}|)(|y| + |\varphi|)|\nabla y||\nabla x| \, dt$$

where we have used the fact that $\nabla(x \cdot y - \varphi) = 2y\nabla x$. Replacing $\zeta$ with $\zeta^4$, this implies

$$\int_{D_2} \zeta(x(t))^4|\nabla y|^2 \, dt \leq c \int_{D_2} \zeta^2(\zeta|\zeta_x| + \zeta|\zeta_{xx}| + |\zeta_x|^2)(|y| + |\varphi|)|\nabla y||\nabla x| \, dt.$$

We may then use the arithmetic-geometric mean inequality and a simple choice of $\zeta$ (with $0 \leq \zeta \leq 1$) to obtain

$$\int_{D_{5/4}} |\nabla y|^2 \, dt \leq c \int_{D_{15/8}} (|\ell - \ell_0|^2 + \varphi^2)|\nabla \ell|^2 \, dt.$$

Since $|\nabla \ell|^2 \leq 2|\nabla \ell_0|^2 + 2|\nabla(\ell - \ell_0)|^2$, and $|\nabla \ell_0|^2 = 1$, we have

$$\int_{D_{5/4}} |\nabla y|^2 \, dt \leq c \int_{D_{15/8}} (|\ell - \ell_0|^2 + \varphi^2) \, dt + \varepsilon \int_{D_2} |\nabla(\ell - \ell_0)|^2 \, dt \tag{4.3}$$

where we have used the fact that $|\ell - \ell_0|^2 + \varphi^2$ is pointwise small (which follows from Proposition 3.3).

To complete our proof, we now consider the conformality relation for the map $\ell$. If we introduce complex coordinates $\tau = t_1 + \sqrt{-1}t_2$, $\xi = x_1 + \sqrt{-1}x_2, \eta = y_1 + \sqrt{-1}y_2$ then we may write this relation $\xi_\tau \bar{\xi}_\tau + \eta_\tau \bar{\eta}_\tau = 0$. In particular we have

$$\int_{D_{9/8}} |\xi_\tau||\bar{\xi}_\tau| \, dt \leq \int_{D_{9/8}} |\nabla y|^2 \, dt. \tag{4.4}$$

We now define a measurable subset $B \subset D_{9/8}$, the 'bad set' to be the set of points for which $|\xi_\tau| \leq |\bar{\xi}_\tau|$. Observe that

$$\int_{D_{9/8} \setminus B} |\bar{\xi}_\tau|^2 \, dt \leq \int_{D_{9/8}} |\nabla y|^2 \, dt. \tag{4.5}$$

We now estimate the area of $\ell(B)$. We claim that if $t \in B$, there is a closed disk $D^t \subset D_{5/4}$ centered at $t$ such that $A(D^t) \leq c \int_{D^t} |\nabla y|^2 \, dt$. This follows from Proposition 4.3, since for any disk, the reverse inequality implies that $\ell$ is close to a holomorphic or anti-holomorphic map. Since for small disks centered at $t$, the map $\ell$ cannot be close to a holomorphic map, while for the disk of radius $1/8$, the map is close to the holomorphic map $\ell_0$, there must be a disk in which our inequality holds. We consider the collection $\mathcal{B}$ of such closed disks, and apply the Besicovitch covering lemma to obtain disjoint subcollections $\mathcal{B}_1, \cdots, \mathcal{B}_q$ whose union covers $B$. Summing the inequalities over each $\mathcal{B}_k$, and then summing over $k$ we obtain

$$\int_B |\nabla \ell|^2 \, dt \leq cq \int_{D_{5/4}} |\nabla y|^2 \, dt. \tag{4.6}$$

Combining the inequalities (4.5) and (4.6) we get

$$\int_{D_{9/8}} |\bar{\xi}_\tau|^2 \, dt \leq c \int_{D_{5/4}} |\nabla y|^2 \, dt. \tag{4.7}$$

The standard $L^2$ estimate for the Cauchy-Riemann operator then implies

$$\int_{D_1} |\nabla(x-t)|^2 \, dt \leq c \int_{D_{9/8}} |\bar{\xi}_\tau|^2 \, dt + c \int_{D_2} |x-t|^2 \, dt.$$

Combining this inequality with (4.3) and (4.7), we obtain

$$\int_{D_1} |\nabla(\ell - \ell_0)|^2 \, dt \leq c \int_{D_2} (|\ell - \ell_0|^2 + \varphi^2) \, dt + \varepsilon \int_{D_2} |\nabla(\ell - \ell_0)|^2 \, dt. \tag{4.8}$$

To complete the proof we estimate the $L^2$ norm of $\varphi$ as follows. Since we are free to subtract a constant from $\varphi$, we may assume that

$$\int_{D_2} (\varphi - y \cdot t) \, dt = 0,$$

and apply the Poincaré inequality to obtain

$$\int_{D_2} (\varphi - y \cdot t)^2 \, dt \leq c \int_{D_2} |\nabla(\varphi - y \cdot t)|^2 \, dt.$$

Since $d\varphi = x\,dy - y\,dx$, we see that $d(\varphi - y \cdot t) = y\,d(x-t) - (x-t)\,dy - 2y\,dt$, and this implies

$$\int_{D_2} (\varphi - y \cdot t)^2 \, dt \leq c \left[ \int_{D_2} y^2 \, dt + \int_{D_2} |\ell - \ell_0|^2 |\nabla(\ell - \ell_0)|^2 \, dt \right],$$

and therefore

$$\int_{D_2} \varphi^2 \, dt \leq c \left[ \int_{D_2} y^2 \, dt + \int_{D_2} |\ell - \ell_0|^2 |\nabla(\ell - \ell_0)|^2 \, dt \right].$$

Combining this inequality with (4.8) then gives the result (using that the pointwise norm of $\ell - \ell_0$ is as small as we wish). □

We now begin the proof of Theorem 4.1. For a map $\ell$, a linear conformal map $\ell_0$, and a radius $r$, define the following excess-type quantity

$$E(\ell, \ell_0, r) = \max\{r^{-2} \int_{D_r} |\nabla(\ell - \ell_0)|^2 \, dt, r \sup_{D_r} |\partial g| \circ \ell\}$$

where $\partial g$ refers to ambient derivatives of $g$. We will prove the desired result by showing that $E$ decays like a power of $r$ provided that we change $\ell_0$ appropriately as a function of $r$. We will exploit our freedom to rescale the problem and change coordinates. We may do a unitary change of coordinates so that the lagrangian plane which is the image of $\ell_0$ becomes the $x$ plane. We then write our map as $\ell(t) = (x(t), y(t))$. We may also rescale the map so that we may assume that $r = 1$. We will prove the theorem by establishing the following decay result for $E$. Observe that the hypotheses of the theorem imply

$$E(\ell, \ell_0, r) < \varepsilon_1 \tag{4.9}$$

for the original map.

**Proposition 4.5.** *There exists $\varepsilon_1 > 0$ and a number $\theta \in (0, 1/2)$ such that if (4.9) holds for $\varepsilon_1$, then there is a linear conformal map $\ell_1$ with image a lagrangian plane such that*

$$E(\ell, \ell_1, \theta) \le \frac{1}{2} E(\ell, \ell_0, 1). \tag{4.10}$$

*Proof.* The proof is by a blow-up argument. We will take $\theta$ to be a small constant to be determined. Assume for the sake of contradiction that the conclusion of the theorem fails. We consider sequences $\ell^{(j)}, g^{(j)}$ with $g^{(j)}$ converging in $C^1$ norm to the euclidean metric $\delta$, and

$$\lim_{j \to \infty} E(\ell^{(j)}, (t, 0), 1) = 0.$$

We may write our maps $\ell^{(j)} = (x^{(j)}, y^{(j)})$. We then consider the sequence of vector valued maps $u^{(j)} = \varepsilon_j^{-1/2}(x^{(j)} - t) - \alpha_j$ and $v^{(j)} = \varepsilon_j^{-1/2} y^{(j)} - \beta_j$ where $\varepsilon_j = E(\ell^{(j)}, (t, 0), 1)$. Therefore we have

$$\int_{D_1} (|\nabla u^{(j)}|^2 + |\nabla v^{(j)}|^2) \, dt \le 1,$$

and assume that $\alpha_j, \beta_j$ are chosen so that $\int u^{(j)} = \int v^{(j)} = 0$. The sequences $u^{(j)}, v^{(j)}$ are bounded in $W^{1,2}$, so we may choose a subsequence of $u^{(j)}, v^{(j)}$ which converges in $L^2$ and weakly in $W^{1,2}$ on the unit disk to limit maps $u, v$. We now show that $u$ is holomorphic and $v = \nabla w$ for a biharmonic function $w$.

We first consider $u$. We may use the conformality relation for the map $\ell^{(j)}$ in the following way. Observe that for the Euclidean metric $\delta$ we may write conformality in the complex form $\frac{\partial x}{\partial t} \frac{\partial \bar{x}}{\partial t} + \frac{\partial y}{\partial t} \frac{\partial \bar{y}}{\partial t} = 0$ where we are using complex notation $x = x^1 + ix^2$, $y = y^1 + iy^2$, $t = t^1 + it^2$. The conformality relation for $\ell^{(j)}$ using the metric $g^{(j)}$ is then:

$$\frac{\partial x^{(j)}}{\partial t} \frac{\partial \bar{x}^{(j)}}{\partial t} + \frac{\partial y^{(j)}}{\partial t} \frac{\partial \bar{y}^{(j)}}{\partial t} = 0.$$

On the other hand we have:

$$|\int_{D_1} \frac{\partial (x^{(j)} - t)}{\partial t} \frac{\partial (\bar{x}^{(j)} - \bar{t})}{\partial t} + \frac{\partial y^{(j)}}{\partial t} \frac{\partial \bar{y}^{(j)}}{\partial t}| \le c\varepsilon_j.$$

Since $x = t$ is holomorphic, we have
$$\frac{\partial x^{(j)}}{\partial t}\frac{\partial \bar{x}^{(j)}}{\partial t} = \frac{\partial(\bar{x}^{(j)} - \bar{t})}{\partial t} + \frac{\partial(x^{(j)} - t)}{\partial t}\frac{\partial(\bar{x}^{(j)} - \bar{t})}{\partial t}.$$

We may combine these to obtain
$$\int_{D_1} |\frac{\partial(\bar{x}^{(j)} - \bar{t})}{\partial t}| dt \leq c\varepsilon_j.$$

Thus we may divide by $\sqrt{\varepsilon_j}$ to obtain
$$\int_{D_1} |\frac{\partial \bar{u}^{(j)}}{\partial t}| dt \leq c\sqrt{\varepsilon_j}.$$

It follows that the weak limit $u$ is a (weakly) holomorphic function, as claimed.

To analyze the limit of $v^{(j)}$, we consider hamiltonians on $N$ of the form $h(x)$, so that,
$$X_h = \sum h_{x_i} \frac{\partial}{\partial y_i}$$

We use the first variation formula together with the fact that the metric is close to the euclidean metric to obtain,
$$|\int_{D_1} \nabla y^{(j)} \cdot \nabla X^{(j)} \, dt| \leq c\varepsilon_j \tag{4.11}$$

where $X^{(j)} = X_h \circ x^{(j)}$ (and $\nabla$ represents the $t$ derivative). By the Schwarz inequality, and the fact that $h$ is a smooth function of $x$ we have
$$|\int_{D_1} \nabla y^{(j)} \cdot DDh \nabla(x^{(j)} - t) dt| \leq c\varepsilon_j, \tag{4.12}$$

where $D$ represents derivatives with respect to $x$. By the chain rule $\nabla X^{(j)} = DDh \cdot \nabla x^{(j)}$, so we may combine with the previous inequality to obtain
$$|\int_{D_1} \nabla y^{(j)} \cdot DDh \, dt| \leq c\varepsilon_j. \tag{4.13}$$

We define the function $h(t)$ on $D_1$ by substituting $(t_1, t_2)$ for $(x_1, x_2)$, and observe that the functions $DDh$ then converge to $\nabla\nabla h$ in $L^2$ norm. We may then divide by $\sqrt{\varepsilon_j}$, and conclude that the weak limit $v$ satisfies
$$\int_{D_1} \nabla v \cdot \nabla\nabla h \, dt = 0 \tag{4.14}$$

for any $C^2$ function $h$ with compact support on $D_1$. On the other hand, the map $(x^{(j)}(t), v^{(j)}(t))$ is weakly lagrangian for each $j$. Since this sequence is bounded and converges weakly in $W^{1,2}$, it follows that the limit map $(t, v(t))$ is also weakly lagrangian. This condition implies that there is a function $w \in W^{2,2}(D_1)$ such that $v = \nabla w$. The equation above then shows that $w$ is a biharmonic function.

We may now complete the proof of the proposition by using standard estimates for holomorphic and biharmonic functions. In particular we have for $r < 1$,

$$\sup_{D_r} |\nabla^s u|^2 \leq c \int_{D_1} |\nabla u|^2 \, dt \tag{4.15}$$

for $s \geq 0$. Let $u_0$ denote the first order part of the Taylor expansion of $u$ at $t = 0$. Then using the second order Taylor series expansion of $u$ and the fact that $u$ is bounded in $L^2$ we have,

$$\sup_{D_r} |u - u_0|^2 \leq cr^2 \int_{D_1} |\nabla u|^2 \, dt. \tag{4.16}$$

Similarly, let $v_0$ denote the first order part of the Taylor expansion of $v$ at $0$. Then,

$$\sup_{D_r} |v - v_0|^2 \leq cr^2 \int_{D_1} |\nabla v|^2 \, dt. \tag{4.17}$$

It follows that

$$\int_{D_{\bar{\theta}}} [(u - u_0)^2 + (v - v_0)^2] \, dt \leq c\bar{\theta}^4 \int_{D_1} (|\nabla u|^2 + |\nabla v|^2) \, dt \leq c\bar{\theta}^4,$$

and hence for $j$ sufficiently large

$$\int_{D_{\bar{\theta}}} [(u^{(j)} - u_0)^2 + (v^{(j)} - v_0)^2] \, dt \leq c\bar{\theta}^4.$$

Multiplying through by $\varepsilon_j^2 = E(\ell^{(j)}, \ell_0^{(j)}, 1)$ we get

$$\int_{D_{\bar{\theta}}} (\ell^{(j)} - \hat{\ell}_1^{(j)})^2 \, dt \leq c\bar{\theta}^4 E(\ell^{(j)}, \ell_0^{(j)}, 1)$$

where $\hat{\ell}_1^{(j)} = \ell_0^{(j)} + \varepsilon_j(u_0 + \alpha_j, v_0 + \beta_j)$. We observe that we can modify $\hat{\ell}_1^{(j)}$ by a term of order $0(\varepsilon_j^2)$ to make it a conformal map into a lagrangian plane. First note that $v_0 + \beta_j$ is the gradient of a quadratic function of $t$, so we can replace $v_0 + \beta_j$ by the gradient of the same quadratic function of the variable $x_0^{(j)} + \varepsilon_j(u_0 + \alpha_j)$. The resulting perturbation of $\hat{\ell}_1^{(j)}$ is of order $\varepsilon_j^2$, and makes the image a lagrangian plane. Since $x_0^{(j)} + \varepsilon_j(u_0 + \alpha_j)$ is a holomorphic map, it follows that the Hopf differential of $\hat{\ell}_1^{(j)}$ is of order $\varepsilon_j^2$. A perturbation of order $\varepsilon_j^2$ will then produce a linear conformal map, which is then the desired map $\ell_1^{(j)}$. We then have

$$\int_{D_{\bar{\theta}}} (\ell^{(j)} - \ell_1^{(j)})^2 \, dt \leq c\bar{\theta}^2(\bar{\theta}^2 + \varepsilon_j^2) E(\ell^{(j)}, \ell_0^{(j)}, 1) \tag{4.18}$$

and

$$\int_{D_{\bar{\theta}}} |\nabla(\ell^{(j)} - \ell_1^{(j)})|^2 \leq \int_{D_{\bar{\theta}}} |\nabla(\ell^{(j)} - \ell_0^{(j)})|^2 + \int_{D_{\bar{\theta}}} |\nabla(\ell_0^{(j)} - \ell_1^{(j)})|^2 \leq c\varepsilon_j^2. \tag{4.19}$$

We may now apply (4.18), (4.19) and Proposition 4.4 to assert that for $j$ sufficiently large

$$\int_{D_\theta} |\nabla(\ell^{(j)} - \ell_1^{(j)})|^2 \leq c\theta^3 E(\ell^{(j)}, \ell_0^{(j)}, 1)$$

where $\theta = \frac{1}{2}\bar{\theta}$. On the other hand we clearly have

$$\theta \sup_{D_\theta} |\partial g| \circ \ell^{(j)} \leq \theta E(\ell^{(j)}, \ell_0^{(j)}, 1).$$

Thus, for $j$ sufficiently large, the inequality (4.10) holds in any sequence chosen as above. It follows that for $\varepsilon_1$ sufficiently small the inequality (4.10) holds for a fixed $\theta$. This completes the proof of the proposition. □

Now to complete the proof of the theorem, one iterates the inequality (4.10) with varying center point to show that for any $r$ small and $t \in D_{1/2}$, there is a conformal map $\ell_{t,r}$ with image a lagrangian plane such that

$$r^{-2} \int_{D_r(t)} |\nabla(\ell - \ell_{t,r})|^2 \, dt \leq cr^{2\alpha}$$

for a positive exponent $\alpha$. This implies that $\ell$ is $C^{1,\alpha}$ in $D_{1/2}$. The smoothness now follows from the following result.

**Proposition 4.6.** *Assume that $\ell$ is a $C^{1,\alpha}$ lagrangian stationary map defined on $D = D_1$ with $|\nabla \ell|(t) \neq 0$ for all $t \in D$. It follows that $\ell$ is a smooth immersion on $D$.*

*Proof.* It clearly suffices to show that $\ell$ is smooth in a neighborhood of $t = 0$. To do this we assume by choice of Darboux coordinates $(x, y)$ that $\ell(t) = (x(t), y(t))$ with $x(0) = y(0) = 0$ and $dy(0) = 0$, $dx(0) = I$. We may then describe a neighborhood of the origin as a $C^{1,\alpha}$ graph $y = y(x)$ with $y^i = \frac{\partial u}{\partial x^i}$ for $i = 1, 2$ for a $C^{2,\alpha}$ function $u(x)$. The stationary condition applied to the hamiltonian vector field $X = \eta_x \frac{\partial}{\partial y}$ implies

$$\int \langle \nabla \sum_i y^i \frac{\partial}{\partial y^i}, \nabla \sum_j \eta_j \frac{\partial}{\partial y^j} \rangle \, dt = 0$$

for any smooth $\eta(x)$ with compact support. We may do a change of variable to rewrite the equation in the form

$$\int \sum_{i,j,k,l} a^{ijkl}(x, u_x, u_{xx}) u_{x^i x^k} \eta_{x^j x^l} \, dx = 0$$

where $a^{ijkl}(x, p, q)$ is a smooth function of its arguments, and $a^{ijkl}(0, 0, 0) = \delta_{ik}\delta jl$. Since $u$ is $C^{2,\alpha}$ this is a linear elliptic equation (in divergence form) with $C^{0,\alpha}$ coefficients. In order to get a gain in regularity, we can form the difference quotient $u^{(h)}(x) = h^{-1}(u(x+he) - u(x))$ where $e$ is any (euclidean) unit vector and $h \neq 0$. We then see that $u^{(h)}$ satisfies the equation

$$\int \sum (a^{ijkl}(x, u_x, u_{xx}) u_{x^i x^k})^{(h)} \eta_{x^k x^l} \, dx = 0.$$

By a standard manipulation using the Fundamental Theorem of Calculus we have

$$\begin{aligned}(a^{ijkl}(x,u_x,u_{xx})u_{x^ix^k})^{(h)} &= a^{ijkl}(x+he, u_x(x+he), u_{xx}(x+he))u^{(h)}_{x^ix^k} \\ &+ (a^{ijkl}(x,u_x,u_{xx}))^{(h)} u_{x^ix^k} \\ &= \tilde{a}^{ijkl}(x)u^{(h)}_{x^ix^k} + \tilde{b}^{jlk}u^{(h)}_{x^k} + \tilde{c}^{jl}\end{aligned}$$

where $\tilde{a}, \tilde{b}, \tilde{c}$ are $C^{0,\alpha}$ and $\tilde{a}^{ijkl}(0) = \delta_{ik}\delta_{jl}$. We may then apply the divergence form Schauder theory (see [Si2]) for a recent treatment) to get a uniform $C^{2,\alpha}$ estimate on $u^{(h)}$. This implies that $u$ is a $C^{3,\alpha}$ function in a neighborhood of 0. This argument can be repeated to show that $u$ is smooth. □

We now consider the global regularity of minimizing lagrangian maps. We expect the results to hold for stationary maps, but the proof is somewhat more complicated, so for the remainder of this section we assume that our maps are minimizing. The following strong compactness theorem will be important for our methods.

**Proposition 4.7.** *Assume that $\ell_j : D \to N$ is a weakly conformal, minimizing lagrangian map for each $j = 1, 2, \ldots$ with $A(\ell_j(D)) \leq c$ for a constant $c$. There is a subsequence of $\{\ell_j\}$ which converges strongly in $W^{1,2}_{loc}(D)$ to a minimizing lagrangian map $\ell$.*

The proof of this result will require the following comparison result.

**Lemma 4.8.** *Suppose that $\ell_0, \ell_1 : S^1 \to \mathbb{R}^5$ are continuous maps with lengths $L_i = L(\ell_i(S^1))$ for $i = 1, 2$, and with $\sup_{S^1} d(\ell_0, \ell_1) \leq \varepsilon \ll L_0 + L_1$. There exists a legendrian map $\ell : S^1 \times [0,1] :\to \mathbb{R}^5$ with $\ell(\theta, 0) = \ell_0(\theta)$, $\ell(\theta, 1) = \ell_1(\theta)$, and $A(\ell(S^1 \times [0,1])) \leq c(L_0 + L_1)\varepsilon$ for a constant $c$.*

*Proof.* Choose an integer $N$ so that $(L_0 + L_1)/N \in [\varepsilon, 2\varepsilon)$, and divide $S^1$ into $N$ intervals $I_p = [\theta_{p-1}, \theta_p)$, $p = 1, \ldots, N$ such that $L(\ell_0(I_p)) + L(\ell_1(I_p)) = (L_0 + L_1)/N$. For each $p$, join the points $\ell_0(\theta_p), \ell_1(\theta_p)$ by a horizontal curve $\ell_p(t)$, $t \in [0,1]$ of length bounded by $\varepsilon$. By the isoperimetric inequality Proposition 2.7 we may span the closed curve $\Gamma_p = \ell_0(I_p) \cup \ell_p \cup \ell_1(I_p) \cup \ell_{p-1}$ with a legendrian disk with area bounded by $c\varepsilon^2$ since the length of this curve is bounded by $4\varepsilon$. We parametrize this disk on $I_p \times [0,1]$ so that it agrees with the specified boundary parametrization. This defines the map $\ell$ on $S^1 \times [0,1]$ with area at most $N \cdot c\varepsilon^2 \leq c(L_0 + L_1)\varepsilon$. □

We now give the proof of Proposition 4.7.

*Proof.* We may choose a subsequence also denoted $\{\ell_j\}$ which converges to a limit $\ell$ uniformly on compact subsets of $D$, and such that for almost all $r \in (0,1)$, the sequence $L(\ell_j(\partial D_r))$ is bounded. We let $r_0 < 1$, and we show that $\ell_j$ converges strongly to $\ell$ in $W^{1,2}(D_{r_0})$, and that $\ell$ is minimizing in $D_{r_0}$. We may choose $r \in (r_0, 1)$ so that $L(\ell_j(\partial D_r)) \leq c$. Given any $\varepsilon > 0$, we may choose $j$ sufficiently large so that $\sup_{D_r} d(\ell_j, \ell) < \varepsilon$. We may then use Lemma 4.8 to construct a map $\hat{\ell}$ which agrees with $\ell(2t)$ in $D_{r/2}$, and agrees with $\ell_j$ on $\partial D_r$ such that $A(\hat{\ell}(D_r)) \leq A(\ell(D_r)) + c\varepsilon$. Since $\ell_j$ is conformal and the area is bounded by the energy, it follows that $\limsup E_{D_r}(\ell_j) \leq E_{D_r}(\ell)$, and hence $\ell_j$ converges strongly to $\ell$ in $W^{1,2}(D_r)$.

To show that $\ell$ is minimizing in $D_{r_0}$, let $\tilde{\ell}$ be any $W^{1,2}$ map which agrees with $\ell$ outside $D_{r_0}$. Construct a map $\tilde{\ell}_j$ which agrees with $\tilde{\ell}$ inside $D_{r_0}$, and agrees with $\ell_j$ on $\partial D_r$ so that $A(\tilde{\ell}_j(D_r)) \leq A(\tilde{\ell}(D_r)) + o(j)$. Since $\ell_j$ is area minimizing, it follows that $A(\ell_j(D_r)) \leq A(\tilde{\ell}_j(D_r))$. Thus, letting $j \to \infty$, we have $A(\ell(D_r)) \leq A(\tilde{\ell}(D_r))$ as required. □

We now use this result to construct parametrized tangent cones which are again area minimizing. Given any point $P \in \Sigma$, we may choose conformal coordinates $t$ centered at $P$, and for any sequence $\varepsilon_j \to 0$ we consider the map $\ell_j(t) = \delta_j^{-1}\ell(\varepsilon_j t)$ where $\delta_j \to 0$ is chosen so that $A(\ell_j(D_2)) = 1$. There is a constant $c_1$ (depending on $\ell$) such that $A(\ell(D_1)) \leq c_1$ and $\inf_{t \in D_{1/4}} A(\ell(D_{1/8}(t))) \geq c_1^{-1} A(\ell(D_1))$, so by Proposition 4.2, we have $A(\ell_j(D_1)) \geq c^{-1}$. Furthermore, we have $A(\ell_0(D_{2a})) \leq cA(\ell_0(D_a))$ for $a \geq 1$, and hence it follows that $A(\ell_0(D_r)) \leq cr^p$ for some $p > 0$ and $r \geq 2$. Applying Proposition 4.7 we see that a subsequence of $\ell_j$ converges strongly in $W^{1,2}_{loc}(\mathbb{R}^2)$ to a nonconstant minimizing, weakly conformal, lagrangian map $\ell_0 : \mathbb{R}^2 \to \mathbb{R}^4$. We have the following properties of $\ell_0$.

**Lemma 4.9.** *The map $\ell_0$ is a proper map from $\mathbb{R}^2$ into $\mathbb{R}^4$, $\ell_0^{-1}(0) = 0$, and $\ell_0(\mathbb{R}^2)$ is a cone in $\mathbb{R}^4$.*

*Proof.* We assume that $\ell_0(0) = 0$, and we show that $\ell_0^{-1}(0)$ is a discrete set of points. To see this we recall the construction of the logarithmic function $\zeta_\delta(P)$ in the proof of Proposition 4.2 which vanishes for $\rho(P) \leq \delta$, and is equal to 1 for $\rho(P) \geq 1$. The function $\zeta_\delta \circ \ell_0$ then vanishes in a neighborhood of $\ell_0^{-1}(0)$, and the Dirichlet integral (in any disk centered at 0) of $\zeta_\delta \circ \ell_0$ tends to zero as $\delta \to 0$. Since $\ell_0$ is not constant, it follows that the closed set $\ell_0^{-1}(0)$ has zero logarithmic capacity, and hence is a totally disconnected zero dimensional set. To show that it is a discrete set, choose any radius $r$, and assume by perturbing $r$ slightly if necessary that $d(\ell_0(\partial D_r), 0) = \delta > 0$. Now if $\sigma < \delta$, we see from monotonicity that there are a bounded number of components of $\ell_0^{-1}(B_\sigma)$ whose image under $\ell_0$ contains 0. Since this number can only increase as $\sigma$ decreases, it must remain constant for $\sigma \leq \sigma_0$ for some $\sigma_0 > 0$. Each of these components must then contain exactly one point of $\ell_0^{-1}(0)$, and it is a discrete set.

In order to show that $\ell_0^{-1}(0) = \{0\}$, we go back to the original map $\ell$, and we note that again by monotonicity (upper and lower area bounds), for any $\sigma > 0$ there is a fixed bound on the number of components of $\ell^{-1}(B_\sigma)$ whose image under $\ell$ intersects $B_{\sigma/2}$. For a given $\sigma$, denote these components by $C_1^{(\sigma)}, \ldots, C_k^{(\sigma)}$ where $C_1^{(\sigma)}$ is the connected component containing 0. Observe that for $\rho < \sigma$, any component $C_j^{(\rho)}$ is contained in one of the components $C_l^{(\sigma)}$. Thus the number of components increases as $\sigma$ decreases, so there is a $\sigma_0 > 0$ such that the number of components is constant for $\sigma \leq \sigma_0$. Thus there is a neighborhood $\Omega$ of $0 \in \mathbb{R}^2$ such that for any $\sigma \leq \sigma_0$, there is precisely one component of $\ell^{-1}(B_\sigma) \cap \Omega$ whose image under $\ell$ intersects $B_{\sigma/2}$. We may now show that $\ell_0^{-1}(0) = \{0\}$, for if there were another $t_0 \in \mathbb{R}^2$ with $\ell_0(t_0) = 0$, there would be $\sigma > 0$ such that $\ell_0^{-1}(B_\sigma)$ has at least two connected components whose image contains 0. Thus for $j$ sufficiently large there would be at least two connected components of $\ell_j^{-1}(B_\sigma)$ whose image intersects $B_{\sigma/2}$. This contradicts the above condition on $\ell$ at radius $\delta_j \sigma$, since $\delta_j \sigma$ is arbitrarily small when $j$ is large. Thus we have shown that $\ell_0^{-1}(0) = \{0\}$.

We now observe that the area $A(\ell_0(\mathbb{R}^2))$ is infinite. To see this, observe that if it were finite, we could choose an arbitrarily large radius $R$ with the length of the curve $\ell_0(\partial D_R)$ as small as we choose. The isoperimetric inequality together with the minimizing property of $\ell_0$ would then show that $\ell_0$ is a constant map.

We can now show that the map $\ell_0$ is proper. In fact we claim that there is a constant $\varepsilon_0 > 0$ such that $d(0, \ell(\partial D_R)) \geq \varepsilon_0 \sqrt{A(\ell_0(D_R))}$. To prove this, we argue by contradiction, and suppose there is a sequence $R_i$ with $d(0, \ell(\partial D_{R_i})) \leq \frac{1}{i}\sqrt{A(\ell_0(D_{R_i}))}$. We renormalize the map setting $\ell^{(i)}(t) = (\sqrt{A(\ell(D_{R_i}))})^{-1} \ell_0(R_i t)$ so that we now have $A(\ell^{(i)}(D_1)) = 1$ and $d(0, \ell^{(i)}(\partial D_1)) \leq \frac{1}{i}$. We then choose a subsequence, again denoted $\ell^{(i)}$, converging strongly in $W^{1,2}_{loc}(\mathbb{R}^2)$ to a minimizing map $\tilde{\ell}$. As above the map $\tilde{\ell}$ is nonconstant and satisfies $\tilde{\ell}^{-1}(0) = 0$. This contradicts the fact that $d(0, \tilde{\ell}(\partial D_1)) = \lim d(0, \ell^{(i)}(\partial D_1)) = 0$, and completes the proof of properness.

The fact that $\ell_0(\mathbb{R}^2)$ is a cone follows from Proposition 3.5 □

We now state the main global regularity theorem for minimizing lagrangian maps.

**Theorem 4.10.** *Let $\ell : D_2 \to N$ be an area minimizing, weakly conformal lagrangian map. There is a finite subset $S$ of $D_1$ such that $\ell$ is a smooth immersion on $D_1 \setminus S$. A point $t \in S$ is either a branch point of $\ell$, or a singularity at which $\ell$ has a nonflat tangent cone. The map $\ell$ is smooth across the branch points, and is Lipschitz at the nonflat singularities.*

*Proof.* We first consider the set of points $\Omega \subset \Sigma$ at which every (parametrized) tangent cone is flat; i.e. has image a lagrangian plane. We show that $\Omega$ is an open set and that there is a discrete set of points $B \subset \Omega$ such that $\ell$ is a smooth immersion on $\Omega \setminus B$. Finally we show that $\ell$ is a smooth map on all of $\Omega$ whose differential vanishes at the points of $B$. First consider any point $P \in \Omega$, and choose coordinates $t$ centered at $P$. Let $\ell_0$ be any tangent cone at $P$, and observe that since $\ell_0(\mathbb{R}^2)$ is a lagrangian plane and $\ell_0$ is an energy minimizing lagrangian map, $\ell_0$ is a smooth harmonic map to its image plane. Since $\ell_0$ is also weakly conformal, it follows that $\ell_0$ is either holomorphic or anti-holomorphic. By changing orientation on the image if necessary, we assume that $\ell_0$ is holomorphic. Since $\ell_0$ is proper and $\ell_0^{-1}(0) = \{0\}$, we see that we must have $\ell_0(\tau) = a\tau^n$, where $\tau = t_1 + it_2$, $a$ is a nonzero complex number and $n$ is a positive integer. If, for any tangent cone $n = 1$, it follows from Theorem 4.1 that $\ell$ is a smooth immersion in a neighborhood of $P$. In any case, we claim that there is a disk $D_r$ about 0 such that $\ell$ is a smooth immersion in $D_r \setminus \{0\}$. To verify this, we show that for any sequence $t^{(j)} \to 0$ ($t^{(j)} \neq 0$), the map $\ell$ is a smooth immersion in a neighborhood of $t^{(j)}$ for $j$ sufficiently large. If we had a sequence that violated this condition, we could extract a subsequence again denoted $t^{(j)}$ such that $\ell$ is *not* a smooth immersion in any neighborhood of $t^{(j)}$ for all $j$. We set $\varepsilon_j = |t^{(j)}|$, and extract a subsequence so that the corresponding rescaled sequence $\ell_j(t) = (\delta_j)^{-1}\ell(\varepsilon_j t)$ converges to a tangent cone. Now for the $\ell_j$, there is a point $t^{(j)}/\varepsilon_j$ on the unit circle, near which, $\ell_j$ is not a smooth immersion. Since the tangent map is a smooth immersion on the unit circle, this contradicts Theorem 4.1, and shows that $\ell$ is a smooth immersion in a deleted neighborhood of any $P \in \Omega$. This shows that $\Omega$ is open and $B$ is discrete.

Now consider a point $P \in B$, and choose coordinates $t$ centered at $P$ as above. We wish to show that $\ell$ is a smooth map in a neighborhood of $P$. We first show that $\nabla \ell$ decays near 0. Precisely, we claim that $\sup_{D_r} |\nabla \ell| \leq c r^\alpha$ for any $\alpha \in (0, 1)$. To see this, we show that for $r$ small enough we have $\sup_{\partial D_{r/2}} |\nabla \ell| \leq \theta \sup_{\partial D_r} |\nabla \ell|$ for any $\theta > 1/2$. The proof of this is by contradiction, supposing that there is a sequence $r_i \to 0$ for which the opposite inequality holds, and forming the corresponding rescaling and tangent map construction. The corresponding tangent map is $\ell_0 = a\tau^n$ where $n \geq 2$ (since $P \in B$) and $|a|$ is determined by the condition that $A(\ell_0(D_1)) = 1$. By direct computation, for this map it is true that

$$\sup_{\partial D_{1/2}} |\nabla \ell_0| = (1/2)^{n-1} \sup_{\partial D_1} |\nabla \ell_0|.$$

On the other hand, since $\ell_0$ is the limit of the rescaled sequence,

$$\sup_{\partial D_{1/2}} |\nabla \ell_0| \geq \theta \sup_{\partial D_1} |\nabla \ell_0|.$$

and this gives a contradiction for $\theta > 1/2$. Thus it follows that $\sup_{\partial D_{r/2}} |\nabla \ell| \leq \theta \sup_{\partial D_r} |\nabla \ell|$, and by a simple iteration we get $\sup_{D_r} |\nabla \ell| \leq c r^\alpha$ for any $\alpha \in (0, 1)$ as required. It follows that the map $\ell$ is $C^{1,\alpha}$ in a neighborhood of $P$ for any $\alpha \in (0, 1)$.

To get the higher regularity of $\ell$ near $P \in B$, we need to control the mean curvature near $P$. We claim that

$$\lim_{r \to 0} \int_{D_r \setminus D_{r/2}} |H|^2 \, d\mu = 0 \tag{4.20}$$

where $d\mu = \lambda dt$ is the area form with $\lambda = \frac{1}{2}|\nabla\ell|^2$. In fact, this statement follows from a simple blow-up argument using the facts that: i) $\int |H|^2 \, d\mu$ is scale invariant, ii) that all tangent maps satisfy $H = 0$, and iii) that the convergence to tangent maps is in the $C^2$ topology away from the origin. We now consider the first order elliptic system satisfied by the mean curvature (see the appendix). We let $\sigma$ be the 1-form associated with the mean curvature vector, so that $\sigma(v) = g(\ell_*(v), JH)$ for a tangent vector $v$ on $\mathbb{R}^2$ where $H = \lambda^{-1}\Delta\ell$ and $\Delta\ell$ denoted the Laplacian (tension field) of $\ell$. The form $\sigma$ is then smooth away from $t = 0$, and satisfies the first order elliptic system

$$d\sigma = d*\ell^*\tau + \ell^*(\text{Ric}), \quad \delta\sigma = 0. \tag{4.21}$$

Set,

$$\gamma = \sigma - *(\ell^*\tau).$$

Then (4.21) becomes,

$$d\gamma = \ell^*(\text{Ric}), \quad \delta\gamma = *(\ell^*(d\tau)). \tag{4.22}$$

In order to gain regularity, we must show that this system is satisfied distributionally across $t = 0$. Note that $\int |\sigma|^2 \, dt = \int |H|^2 \, d\mu$ since the $L^2$ norm of 1-forms is conformally invariant. Now we choose a cutoff function $\zeta_r(t)$ satisfying $\zeta_r(t) = 0$ for $|t| \leq r/2$, $\zeta_r(t) = 1$ for $t \geq r$, $|\nabla\zeta_r| \leq 2/r$, and let $\phi(t)$ be any smooth function with compact support in $D_1$. We then have, from (4.22), the equations,

$$\int \delta(\zeta_r\phi dt) \cdot \gamma \, dt = \int \zeta_r\phi\ell^*(\text{Ric}), \quad \int d(\zeta_r\phi) \cdot \gamma \, dt = \int \zeta_r\phi\ell^*(d\tau).$$

The first equation implies

$$\left|\int \zeta_r\{\delta(\phi dt) \cdot \gamma \, dt - \phi\ell^*(\text{Ric})\}\right| \leq$$

$$2(\sup|\phi|)r^{-1}\int_{D_r\setminus D_{r/2}} |\sigma| \, dt + 2(\sup|\phi|)r^{-1}\int_{D_r\setminus D_{r/2}} |\ell^*\tau|\delta t.$$

By the Schwarz inequality, the first term on the right hand side is bounded by a constant times $\int_{D_r\setminus D_{r/2}} |\sigma|^2 \, dt$ which converges to zero by (4.20). The second term clearly converges to zero. It follows from this and the dominated convergence theorem that the first equation of (4.22) is satisfied distributionally. A similar argument applies to establish the weak form of the second equation. Since the right hand side of (4.22) is $C^{0,\alpha}$, it follows from elliptic regularity theory that $\gamma$ is $C^{1,\alpha}$ and that therefore $\sigma$ is $C^{0,\alpha}$. We now observe that the mean curvature may be written in terms of the components of $\sigma = \sum_{j=1}^2 \sigma_j(t)dt^j$

$$H = \lambda^{-1}\Delta\ell = \lambda^{-1}\sum_{j=1}^2 \sigma_j(t)J(\ell_*(\frac{\partial}{\partial t^j})),$$

which implies

$$\Delta\ell = \sum_{j=1}^{2}\sigma_j(t)J(\ell_*(\frac{\partial}{\partial t^j})).$$

As an equation for $\ell$, we see that since $\sigma_j$ is $C^{0,\alpha}$, we may conclude that $\ell$ is $C^{2,\alpha}$, and then we may feed this information back into (4.22) to get $\gamma$ in $C^{2,\alpha}$, and then we may continue inductively to get $\ell$ smooth in a neighborhood of $P$.

To complete the proof we must show that $\Sigma \setminus \Omega$ is a discrete set of points. For this purpose, we consider a point $P$ in this set, and observe that there must be tangent map $\ell_0$ at $P$ whose image is not a plane. Observe that if we take $t_0 \in \mathbb{R}^2 \setminus \{0\}$ in the domain of $\ell_0$, and we take a tangent map $\ell_1$ at $t_0$, then the image of $\ell_1$ is the product of a one dimensional cone in $\mathbb{R}^2$ with a line. The one dimensional cone is a geodesic in $\mathbb{R}^2$, and hence is a line. Therefore the image of $\ell_1$ is a plane. By our previous discussion it follows that $\ell_0$ is a smooth immersion except for a discrete set of points in $\mathbb{R}^2 \setminus \{0\}$. From this it follows that $\ell_0(\mathbb{R}^2) \setminus \{0\}$ is a smooth surface since it is a cone, and the intersection with $\partial B_\sigma$ is an immersed curve for typical values of $\sigma$. It follows that $\ell_0(\mathbb{R}^2)$ is the cone over any such immersed curve. Thus $C = \ell_0(\mathbb{R}^2)$ is a geometric cone, and is conformally equivalent to $\mathbb{R}^2$ via some map $\xi : C \to \mathbb{R}^2$. It follows as above that $\xi \circ \ell_0(\tau) = a\tau^n$. In particular, $\ell_0$ is a smooth immersion on all of $\mathbb{R}^2 \setminus \{0\}$. We may then prove as above that $\ell$ is a smooth immersion on $D_r \setminus \{0\}$ for some $r > 0$. It follows that $\Sigma \setminus \Omega$ is a discrete set as claimed. To prove that $\ell$ is Lipschitz in a neighborhood of $P \in \Sigma \setminus \Omega$, we need to use the fact that for any tangent map $\ell_0$ at $P$ we have $|\nabla \ell_0|(t)| \leq c|t|^\alpha$ for some $\alpha > 0$ (and then we may argue as above). This is a consequence of the structure of the tangent cones discussed in Section 7. In particular, it follows from the fact that the length of the intersection of any such cone with the unit sphere has length strictly greater than $2\pi$. □

## 5. Existence

Let $\alpha$ be a homotopy class of maps $[\Sigma, N]$ (or a homology class in $H_2(N, \mathbb{Z})$). We say a $W^{1,2}$ map $f : \Sigma \to N$ represents $\alpha$ if given any sequence of smooth maps $\{f_i\}$ that strongly approximates $f$ in $W^{1,2}$, for $i$ sufficiently large each smooth map $f_i$ represents $\alpha$. The fact that the $f_i$ lie in the same homotopy class for $i$ sufficiently large follows from the smoothing results of [ScU]. We begin with three lemmas on the homology classes and free homotopy classes of $N$.

**Lemma 5.1.** *There is a constant $C_0 > 0$, depending only on $N$, such that if $f : S^2 \to N$ is a $W^{1,2}$ map with $E(f) < C_0$ then $[f]$ is null in $\pi_2(N)$.*

*Proof.* Suppose not. Then for any $c > 0$ there is a $W^{1,2}$ map $f : S^2 \to N$ with $E(f) < c$ representing a non-zero element of $\pi_2(N)$. By harmonic map theory [SU 1] it follows that there is a harmonic map $h : S^2 \to N$ with $E(h) < c$ representing a non-zero element of $\pi_2(N)$. For $c$ sufficiently small the image of such a harmonic map lies in a coordinate neighborhood and the map is therefore trivial. □

**Lemma 5.2.** *Let $C$ be a positive constant. There are at most finitely many free homotopy classes in $\pi_2(N)$ that can be represented by $W^{1,2}$ maps $f : S^2 \to N$ with $E(f) < C$.*

*Proof.* Suppose not. Then there are infinitely many classes $\alpha_\lambda \in \pi_2(N), \lambda \in \Lambda$ each represented by a $W^{1,2}$ map $f_\lambda : S^2 \to N$ with $E(f_\lambda) < C$. Minimizing energy in each homotopy class we can represent each $\alpha_\lambda$ by a finite family $\{f_{\lambda_i} : i = 1, \ldots, n_\lambda\}$ of harmonic maps $S^2 \to N$ with $[f_{\lambda_i}] \neq 0$ and $\sum_i^{n_\lambda} E(f_{\lambda_i}) < C$. Since $E(f_{\lambda_i}) \geq C_0$, there must be infinitely many different free homotopy classes in the set $\{[f_{\lambda_i}] : \lambda \in \Lambda, i = 1, \ldots, n_\lambda\}$. Hence we can assume that the classes $\alpha_\lambda$ are each represented by a harmonic map $f_\lambda : S^2 \to N$ with $E(f_\lambda) < C$. Choosing a subsequence we can suppose that the $\{f_\lambda\}$ converge up to bubbling. In particular, the images of the $\{f_\lambda\}$ converge in Hausdorff distance to the image of a set of harmonic maps of $S^2$ and connecting curves. Thus there is a $\lambda_0$ such that for $\lambda > \lambda_0$ the $\{f_\lambda\}$ represent the same free homotopy class. The contradiction proves the lemma. □

**Lemma 5.3.** *Let $C$ be a positive constant. There are at most finitely many homology classes in $H_2(N; \mathbb{Z})$ that can be represented by cycles $\Gamma$ with $A(\Gamma) < C$.*

*Proof.* Recall that a cycle $\Gamma$ determines a closed integral current $T_\Gamma$ by:
$$T_\Gamma(\sigma) = \int_\Gamma \sigma,$$
for $\sigma$ a smooth 2-form of compact support in $N$. If $A(\Gamma) < C$ then $T_\Gamma$ satisfies the mass bound $M(T_\Gamma) < C$. The integral homology class determined by $T_\Gamma$ is $[\Gamma]$. The lemma now follows from the compactness theorem of Federer-Fleming and is explicitly stated in [FF, 9.6]. □

In this paper, it will be important to use the area functional rather than the energy for the purpose of proving strong convergence and regularity results. This is because we don't know if such results hold for the energy functional on the space of lagrangian maps. The following proposition will allow us to pass back and forth between these functionals. Let $\Sigma$ be a smooth closed oriented surface, and $\alpha$ a homotopy class of $W^{1,2}$ lagrangian maps from $\Sigma$ to $N$. Let $\mathcal{A}(\alpha)$ denote the infimum of area taken over the maps in this class; i.e.
$$\mathcal{A}(\alpha) = \inf\{A(\ell) : \ell \in W^{1,2}_L(\Sigma, N) \cap \alpha\}.$$

On the other hand, if we fix a metric $h$ on $\Sigma$ with curvature $1, 0$, or $-1$, then we may consider the infimum $\mathcal{E}(\alpha, h)$ given by

$$\mathcal{E}(\alpha, h) = \inf\{E(\ell, h) : \ \ell \in W^{1,2}_L(\Sigma, N) \cap \alpha\}.$$

Of course, we have $\mathcal{A}(\alpha) \leq \mathcal{E}(\alpha, h)$ for any metric $h$. We will refer to the metrics $h$ above as admissible metrics, and we note that each conformal structure on $\Sigma$ has an admissible metric which is unique for genus$(\Sigma) > 1$, unique up to scale for genus$(\Sigma) = 1$, and unique up to conformal transformation if genus$(\Sigma) = 0$. The following proposition is fairly standard, but we need a slightly more general version than is common in classical minimal surface theory, so we record it here.

**Proposition 5.4.** *Given any $\ell \in W^{1,2}_L(\Sigma, N)$ and any $\varepsilon > 0$, there exists an admissible metric $h$ and a homeomorphism $f$ of $\Sigma$ homotopic to the identity so that $E(\ell \circ f, h) \leq A(\ell) + \varepsilon$. In particular, $\mathcal{A}(\alpha) = \inf_h \{\mathcal{E}(\alpha, h)\}$ where the infimum is taken over all admissible metrics on $\Sigma$.*

*Proof.* We consider the pullback metric $\tau = \sum \tau_{ij} dx^i dx^j$ where

$$\tau_{ij} = \langle \ell_*(\partial/\partial x^i), \ell_*(\partial/\partial x^j) \rangle.$$

The tensor $\tau$ defines an $L^1$ degenerate Riemannian metric on $\Sigma$. In order to uniformize this metric, we need to make it bounded and positive definite. We let $h_1$ be any smooth Riemannian metric on $\Sigma$, $\delta > 0$, and set $\tau_\delta = \tau + \delta h_1$. Now, if we choose $\Lambda$ sufficiently large, and we consider the set

$$E_\Lambda = \{P \in \Sigma : \ Tr_{h_1}(\tau)(P) > \Lambda\},$$

then we will have $E(\ell, E_\Lambda, h_1) + A(E_\Lambda, h_1) < \varepsilon/4$. We then define $\hat{\tau}$ to be $\tau_\delta$ on $\Sigma \setminus E_\Lambda$, and equal to $h_1$ on $E_\Lambda$. The metric $\hat{\tau}$ may then be uniformized by the measurable Riemann Mapping Theorem of Morrey ([M2]). Thus we may find a homeomorphism $f$ homotopic to the identity, and an admissible metric $h$ such that $f^*(\hat{\tau})$ is conformal to $h$. Thus it follows that

$$A(\Sigma, f^*(\hat{\tau})) = A(\Sigma, \hat{\tau}) = 1/2 \int_\Sigma Tr_h(f^*(\hat{\tau})) \ d\mu_h.$$

Now we have $A(\Sigma, \hat{\tau}) \leq A(\Sigma, \tau) + c\delta + A(E_\Lambda, h_1) < A(\ell) + c\delta + \varepsilon/4$. The energy of $\ell \circ f$ with respect to $h$ taken over $f^{-1}(E_\Lambda)$ is equal to the energy of $\ell$ with respect to $h_1$ taken over $E_\Lambda$ (which is less than $\varepsilon/4$) since $f$ is a conformal map from $(f^{-1}(E_\Lambda), h)$ to $(E_\Lambda, h_1)$. On the other hand

$$E(\ell \circ f, \Sigma \setminus E_\Lambda, h) = 1/2 \int_{\Sigma \setminus E_\Lambda} Tr_h(f^*(\tau)) \ d\mu_h \leq 1/2 \int_\Sigma Tr_h(f^*(\hat{\tau})) \ d\mu_h.$$

Combining these results we see that $E(\ell \circ f, h) \leq A(\ell) + c\delta + \varepsilon/2$ which gives the desired conclusion for $\delta$ small enough. □

**Minimizing in a homotopy class**

Let $\Sigma$ be a compact surface and let $\alpha \in [\Sigma, N]$ be a lagrangian homotopy class. Let $h$ be an admissible metric on $\Sigma$ for which $\mathcal{E}(\alpha, h) = \mathcal{A}(\alpha)$. Note that if $\Sigma = S^2$, then any constant curvature one metric has this property. Denote the induced volume form by $d\zeta$. Let $\{\ell_i\}$ be an energy minimizing sequence of $W^{1,2}$ lagrangian maps $\Sigma \to N$ representing $\alpha$. Choosing a subsequence we can suppose that $\{\ell_i\}$ converges weakly in $W^{1,2}(\Sigma, N)$, strongly in $L^2(\Sigma, N)$ and pointwise almost everywhere to a weakly lagrangian map $\ell \in W^{1,2}(\Sigma, N)$. Of course the convergence may not be strongly in $W^{1,2}$. To keep track of the energy loss we describe the convergence in measure theoretic language. Consider the energy measures, $\zeta_i = e(\ell_i)d\zeta = \frac{1}{2}|\nabla \ell_i|^2 d\zeta$, on $\Sigma$. Using the weak convergence theorem for measures and choosing a subsequence we can suppose that these measures converge weakly to a limit measure $\eta$. Denote the energy measure associated to the limit map $\ell$ by, $\eta_\ell = e(\ell)d\zeta$. Set,
$$\eta = \eta_\ell + \eta_Z.$$
$\eta_Z$ is a measure, with support $Z$, that records the failure of strong convergence. We will call $\eta_Z$ the *defect measure*. We will show in the next propositions that if $\ell_i$ is an area minimizing sequence, then the defect measure is a finite sum of point masses.

**Proposition 5.5.** *Assume that $h$ is an admissible metric on $\Sigma$ for which $\mathcal{E}(\alpha, h) = \mathcal{A}(\alpha)$, and suppose that $\{\ell_i\}$ is a sequence of weakly lagrangian maps in $W^{1,2}(\Sigma, N)$ that minimize energy in a lagrangian homotopy class $\alpha$. Let $\ell$ be the weak limit. There is an $\varepsilon_0 > 0$, depending only on $N$, such that if $D_r \subset \Sigma$ is a disc of radius $r$ and $E_{D_r}(\ell_i) < \varepsilon_0$ then a subsequence (still denoted $\{\ell_i\}$) converges strongly to $\ell$ in $W^{1,2}(D_{\frac{r}{2}}, N)$.*

*Proof.* Let $\eta$ be the limit energy measure. For any $\varepsilon > 0$ there is at most a finite set $S$ of points, such that if $Q \in S$ then $\eta(\{Q\}) \geq \varepsilon$. For $P \notin S$ there is an $r_0 > 0$ such that if $r < r_0$ and $i$ is sufficiently large then $E_{D_r(P)}(\ell_i) \leq 2\varepsilon$. Fix such a radius $r$ and suppose that for any subsequence the convergence on $D_{\frac{r}{2}}(P)$ results in energy loss. In particular there is a $\delta > 0$ such that,
$$E_{D_{\frac{r}{2}}}(\ell) + \delta < \liminf_i E_{D_{\frac{r}{2}}}(\ell_i).$$
We can choose $\rho \in [r/2, r]$ and a subsequence of $\{\ell_i\}$ ( that we continue to denote $\{\ell_i\}$) such that $\ell_{i|\partial D_\rho}$ has energy bounded by $\frac{2\varepsilon}{r}$. It follows that $\ell_{i|\partial D_\rho}$ is continuous, rectifible and has length bounded by a uniform constant times $\varepsilon$. Thus a subsequence converges uniformly to $\ell_{|\partial D_\rho}$. Using the collar construction of Lemma 4.8, for each sufficiently large $i$, we can construct a map $\tilde{\ell}_i \in W^{1,2}(\Sigma, N)$ that agrees with $\ell$ on $D_{\frac{r}{2}}$, agrees with $\ell_i$ outside $D_r$ and with area,
$$A(\tilde{\ell}_i) + \frac{\delta}{2} < E(\ell_i).$$
The maps $\tilde{\ell}_i$ and $\ell_i$ differ by a $W^{1,2}$ map $s_i : S^2 \to N$ with $A(s_i) < c\varepsilon$, where $c$ is a uniform constant. By Proposition 5.4, this map can be reparametrized to have energy less than $2c\varepsilon$. Choose $\varepsilon < C_0/(2c)$, where $C_0$ is the constant in Lemma 5.1 . Then $s_i$ has energy less than $C_0$ and hence is trivial in homotopy. Thus $\{\tilde{\ell}_i\}$ is a sequence in the same lagrangian homotopy class as $\{\ell_i\}$ but with strictly smaller area, contradicting the area minimizing property of $\{\ell_i\}$. □

The following corollary is an immediate consequence.

**Corollary 5.6.** *Assume that $h$ is an admissible metric on $\Sigma$ for which $\mathcal{E}(\alpha, h) = \mathcal{A}(\alpha)$, and suppose that $\{\ell_i\}$ is a sequence of weakly lagrangian maps in $W^{1,2}(\Sigma, N)$ that minimize energy in a lagrangian homotopy class $\alpha$. Suppose $\ell$ be the weak limit and that the energy measures $\zeta_i = e(\ell_i)d\zeta$ converge weakly as measures to $\eta$. Then,*

$$\eta = \eta_\ell + \eta_Z.$$

*$Z$, the support of $\eta_Z$, consists of a finite number of points $\{x_1, \ldots, x_k\} \subset \Sigma$ and*

$$\eta_Z = \sum_{j=1}^k m_{x_j} \delta_{x_j},$$

*where the masses $m_{x_j} \geq \varepsilon_0$. The constant $\varepsilon_0$ is given in Proposition 5.5.*

**Proposition 5.7.** *Each point $x \in Z$ can be used to construct a finite set of lagrangian stationary, weakly lagrangian maps $\ell_{x_\lambda} \in W^{1,2}(S^2, N)$. Each map, $\ell_{x_\lambda}$, is minimizing in some free lagrangian homotopy class. The construction allows no energy loss:*

$$\liminf_i E(\ell_i) = E(\ell) + \sum_{x_\lambda} E(\ell_{x_\lambda}).$$

*Proof.* Let $x \in Z$. There is a sequence of radii $r_i \to 0$ such that,

$$\liminf_i E_{D_{r_i}(x)}(\ell_i) \geq \varepsilon_0,$$

since otherwise, the convergence in a neighborhood of $x$ is strong. By the Courant-Lebesgue Lemma we can suppose that $L(\ell_i(\partial D_{r_i})) < \varepsilon_i$ where the $\varepsilon_i \to 0$. Set $C_{r_i} = \partial D_{r_i}(x)$. We can suppose that the curves $\{\ell_i(C_{r_i})\}$ all lie in Darboux coordinate neighborhoods $U_i$. The lagrangian isoperimetric inequality, Proposition 2.7 then implies that there are lagrangian maps $m_i : D \to U_i$ of the unit disc $D$ with $\partial m_i(D) = \ell_i(C_{r_i})$ and with,

$$A(m_i(D)) \leq C \left(L(\ell_i(C_{r_i}))\right)^2.$$

The constant $C$ is universal. Rescale the maps $\ell_{i|D_{r_i}}$ so that they are defined on the unit disc. By stereographic projection we consider the rescaled $\ell_i$ defined on the southern hemisphere of $S^2$ and the $m_i$ defined on the northern hemisphere. Denote these maps $(\ell_x)_i : S^2 \to N$. Note that $E((\ell_x)_i) \geq \varepsilon_0$ and so each $(\ell_x)_i$ is nontrivial. Choosing a subsequence we can suppose that each $(\ell_x)_i$ represents a fixed lagrangian homotopy class $\alpha_x$. Then, since the original sequence is minimizing and $A(m_i(D)) \to 0$, the sequence $\{(\ell_x)_i\}$ must be minimizing. To each map $(\ell_x)_i$ we can associate an energy measure $(\zeta_x)_i = e((\ell_x)_i)d\zeta$. We next *balance* each map $(\ell_x)_i$ as follows: Consider the domain $S^2 = \{\vec{x} \in \mathbb{R}^3 : |\vec{x}| = 1\}$. Reparametrize each $(\ell_x)_i$ so that the center of mass of the measure $(\zeta_x)_i$ on $S^2$ is the origin of $\mathbb{R}^3$. We call such a map and its associated energy measure *balanced*. We will denote the balanced sequence $\{(\ell_x)_i\}$.

Choosing a subsequence we can suppose that the balanced sequence $\{(\ell_x)_i\}$ converges weakly in $W^{1,2}(S^2, N)$, strongly in $L^2(S^2, N)$ and pointwise almost everywhere to a limit weakly lagrangian map $\ell_x \in W^{1,2}(S^2, N)$. The measures $(\zeta_x)_i$ converge to a balanced measure $\zeta_x$. If the convergence to $\ell_x$ is strong the lemma is proved. If not, we consider the points where strong convergence fails and iterate the argument. If at each step the map $\ell_x$ is nontrivial then, since the total energy is bounded, the process must terminate. If $\ell_x$ is trivial then the balanced condition implies that there are at least two points in the support of the defect measure of $\zeta_x$ each with mass $\geq \varepsilon_0$ and thus again the process must terminate after finitely many steps. Note that since $A(m_i(D)) \to 0$ there is no energy loss. Clearly $\ell$ and each map $\ell_{x_\lambda}$ minimizes energy in some free lagrangian homotopy class. $\square$

Using the standard nomenclature we will call the maps $\ell_{x_\lambda} : S^2 \to N$, *bubbles*, and the limit process *bubbling*.

Let $\Sigma$ denote a Riemann surface and $f \in W^{1,2}(\Sigma, N)$. Then $f$ determines a current $T_f$ as follows: Let $\tau$ be a smooth 2-form on $N$ with compact support. Define,
$$T_f(\tau) = \int_\Sigma f^*(\tau).$$
We claim that $T_f$ is an integral current and as such it represents a class $[f] \in H_2(N; \mathbb{Z})$. To verify this, first note that if a sequence $\{f_i\}$ in $W^{1,2}(\Sigma, N)$ converges strongly to $f$ then the currents $T_{f_i}$ converge to $T_f$ in the weak (current) topology. Since $W^{1,2}(\Sigma, N)$ is dense in $C^\infty(\Sigma, N)$ in the strong $W^{1,2}$-topology there is a sequence $\{f_i\}$ of smooth maps $\Sigma \to N$ that converge strongly in $W^{1,2}$ to $f$. The currents $\{T_{f_i}\}$ are integral currents that satisfy a mass bound. By the compactness theorem [F] they converge in the flat norm topology to an integral current $T$. It then follows from [Si1,31.2] that $T = T_f$. Thus, $T_f$ is an integral current. Suppose that $\{f_i\}$ is a sequence of maps with $\|f_i\|_{W^{1,2}} < E$ that converges strongly to $f \in W^{1,2}(\Sigma, N)$ and that $[f_i] = \alpha$ for all $i$. Then, as above, the integral currents $T_{f_i}$ converge in the flat norm topology to an integral current $T_f$. By a result of Federer-Fleming [FF] the current $T_f$ also represents $\alpha$ in $H_2(N; \mathbb{Z})$. In other words, for maps from surfaces, strong $W^{1,2}$ convergence preserves homology.

Suppose that $\Sigma = S^2$. Let $\{\ell_i\}$ be an area minimizing sequence of lagrangian maps $S^2 \to N$ representing a free homotopy class and let $\ell$ be the weak limit. Note that $\ell$ is stationary under compactly supported variations of the domain, since such variations preserve the lagrangian condition. It follows by a standard argument (see, for example, [S]) that the Hopf differential is holomorphic and hence, since the domain is $S^2$, that the Hopf differential vanishes. In particular, $\ell$ is weakly conformal.

**Theorem 5.8.** *Let $\alpha \in \pi_2(N)$ be a lagrangian homotopy class. Let $\{\ell_i\}$ be an area minimizing sequence of $W^{1,2}$ lagrangian maps $S^2 \to N$ each representing $\alpha$. Then a subsequence of $\{\ell_i\}$ converges without energy loss to a finite collection $\{\ell_\lambda\}_{\lambda \in \Lambda}$ of lagrangian stationary, weakly conformal, weakly lagrangian maps in $W^{1,2}(S^2, N)$. Each map $\ell_\lambda$ minimizes area among $W^{1,2}$ lagrangian maps in some free homotopy class and*

$$\sum_\lambda [\ell_\lambda] = \alpha,$$

*where $[\ell_\lambda]$ denotes the homology class determined by the map $\ell_\lambda$.*

*Proof.* The first two statements follows from the proof of Proposition 5.7. To prove the third statement we first consider the case where the limit consists of a single map $\ell \in W^{1,2}(S^2, N)$. In particular, a subsequence of $\{\ell_i\}$ converges strongly in $W^{1,2}(S^2, N)$ to $\ell$. Thus, as discussed above, the homology class is preserved in the limit. Next consider the case where the limit consists of more than one map $\ell_\lambda \in W^{1,2}(S^2, N)$. By construction, each map $\ell_\lambda$ is the weak $W^{1,2}$ limit of a minimizing sequence $(\ell_\lambda)_i$ of lagrangian maps $S^2 \to N$. Since there is no energy loss, each sequence must converge strongly. Note that, for each $i$, the homology class $\sum_\lambda [(\ell_\lambda)_i]$ is equal to the class $\alpha$. Thus,

$$\sum_\lambda [\ell_\lambda] = \sum_\lambda [(\ell_\lambda)_i] = \alpha.$$

□

**Non-collapsible surfaces**

Let $\Sigma$ be a Riemann surface with conformal structure $\mu$ and let $\phi : \Sigma \to N$ be a smooth lagrangian map. We will say that a map $\ell : \Sigma \to N$ has the same action on $\pi_1$ as $\phi$ if for $p \in \Sigma$ there is a path $\rho$ from $\ell(p)$ to $\phi(p)$ such that,

$$\ell_* = \rho_* \circ \phi_* \circ \rho_*^{-1},$$

where $\phi_* : \pi_1(\Sigma, p) \to \pi_1(N, \phi(p))$ and $\ell_* : \pi_1(\Sigma, p) \to \pi_1(N, \ell(p))$ are the induced maps on $\pi_1$. We recall that in [SY] it is shown that the induced map on $\pi_1$ is well-defined for any map in $W^{1,2}(\Sigma, N)$ and that the induced map on $\pi_1$ is preserved in the weak limit. We define

$$\mathcal{L}_\phi = \ \{\ell \in W^{1,2}(\Sigma, N) : \ell \text{ is weakly lagrangian and}$$
$$\ell \text{ has the same action on } \pi_1 \text{ as } \phi\}.$$

Set,
$$E_\mu = \inf\{E(\ell, \mu) : \ell \in \mathcal{L}_\phi\}.$$

**Proposition 5.9.** *There exists an $\ell \in \mathcal{L}_\phi$ such that $E(\ell) = E_\mu$.*

*Proof.* Let $\{\ell_i\} \in \mathcal{L}_\phi$ be an energy minimizing sequence. Choosing a subsequence we can suppose that $\{\ell_i\}$ converges weakly in $W^{1,2}(\Sigma, N)$, strongly in $L^2(\Sigma, N)$ and pointwise almost everywhere to a weakly lagrangian map $\ell \in W^{1,2}(\Sigma, N)$ with the same induced action on $\pi_1$ as $\phi$. Since the class $\mathcal{L}_\phi$ is weakly closed, we have $\ell \in \mathcal{L}_\phi$, and hence $\lim E(\ell_i) \leq E(\ell) = E_\mu$. It follows that there is no energy loss, and hence $\ell_i$ converges strongly to $\ell$ in $W^{1,2}$. □

Suppose that $\phi : \Sigma \to N$ is a piecewise $C^1$ lagrangian map. Let $\gamma : S^1 \to \Sigma$ be an oriented simple closed curve on $\Sigma$. Suppose that $\phi_*([\gamma]) = 1$, where $\phi_*$ is the induced map on $\pi_1$. The following discussion applies to any such curve $\gamma$ but is only interesting when $[\gamma]$ is non-trivial in $\pi_1(\Sigma)$. Also, we could carry through the following discussion for any closed curve on $\Sigma$ with similar results. We confine ourselves to simple closed curves because they are suited to our applications. We define the *period of $\gamma$* as follows. Let $D$ be an oriented immersed disc in $N$ with $\partial D = \phi(\gamma)$. Set:
$$\text{period}_D(\phi, \gamma) = \int_D \omega.$$
This is well-defined independent of the choice of oriented representative $\gamma \in [\gamma]$. If $D'$ is another oriented immersed disc in $N$ with $\partial D' = \phi(\gamma)$ the pair $D, D'$ determine an element $\alpha$ of $\pi_2(N)$. Then $\text{period}_D(\gamma) - \text{period}_{D'}(\gamma) = \omega(\alpha)$. The period of $\gamma$ is thus well-defined only modulo the values of $\omega$ on $\pi_2(N)$. If the values of $\omega$ on $\pi_2(N)$ are rationally related then we can define:
$$\text{period}(\phi, \gamma) = \inf_{\{D:\, \partial D = \phi(\gamma)\}} |\int_D \omega|.$$
However, if the values of $\omega$ on $\pi_2(N)$ are not rationally related then this definition yields $\text{period}(\phi, \gamma) = 0$. Accordingly, in the general case, we choose a constant $C > 0$ and consider:
$$\text{per}_C(\phi, \gamma) = \inf_{\{D:\partial D=\phi(\gamma), A(D)<C\}} |\int_D \omega|.$$
If there are no discs $D \subset N$ with $\partial D = \phi(\gamma)$ and $A(D) < C$ then we set $\text{per}_C(\phi, \gamma) = \infty$. Because we restrict the area of the spanning disc to be less than $C$ the integral of $\omega$ over a spanning disc is well-defined modulo the values of $\omega$ on classes in $\pi_2(N)$ that can be represented by two-spheres with energy less than $2C$. By Lemma 5.2 there are at most finitely many such classes. Hence the infimum is taken over finitely many classes. In particular, if $\text{per}_C(\phi, \gamma) = 0$ then there is an immersed disc $D$ with $\partial D = \phi(\gamma)$ and $\int_D \omega = 0$. Unfortunately $\text{per}_C(\phi, \gamma)$ depends on the representative $\gamma \in [\gamma]$. Define:

$$\text{period}_C(\phi, [\gamma]) = \inf_{\gamma \in [\gamma]} \text{per}_C(\phi, \gamma).$$

Consider the energy functional $E(S^2, N)$ on piecewise $C^1$ maps $S^2 \to N$. We will require that $2C$ be a regular value of $E(S^2, N)$. The next lemma describes the behavior of $\text{period}_C$ under hamiltonian variations.

**Lemma 5.10.** *Let $\Sigma$ be a compact surface. Let $[\gamma]$ be a non-trivial homotopy class that can be represented by a simple closed curve. Suppose that $\phi : \Sigma \to N$ is a piecewise $C^1$ lagrangian map such that $\phi_*([\gamma]) = 0$. Let $h_t$, $0 \leq t \leq 1$ be a 1-parameter family of hamiltonian diffeomorphisms with $h_0 = Id$. Then there is an $\varepsilon > 0$ such that if $t < \varepsilon$ then:*

$$\text{period}_C(h_t \circ \phi, [\gamma]) = \text{period}_C(\phi, [\gamma]).$$

*Proof.* Choose a simple closed curve $\gamma_0$ on $\Sigma$ that represents $[\gamma]$. Set,

$$S(\phi, \gamma_0, C) = \{|\int_D \omega| : A(D) < C, \partial D = \phi(\gamma_0)\}.$$

For any $C$ this set of non-negative real numbers is finite. For any simple closed curve $\gamma$ on $\Sigma$ that represents $[\gamma]$ we have:

$$S(\phi, \gamma, C) \subset S(\phi, \gamma_0, C + A(\phi(\Sigma))).$$

It follows that $\cup_{\gamma \in [\gamma]} S(\phi, \gamma, C)$ is a finite set of values and each value is achieved by a simple closed curve in $[\gamma]$. Denote a set of such curves by $\{\gamma_\lambda : \lambda \in \Lambda\}$ and a set of spanning discs by $\{D_\lambda : \lambda \in \Lambda\}$, where $\Lambda$ labels the elements of $\cup_{\gamma \in [\gamma]} S(\phi, \gamma, C)$. Note that under hamiltonian isotopy the period of a curve does not change. Thus the set $\cup_{\gamma \in [\gamma]} S(\phi, \gamma, C)$ can change under $h_t$ only because of a change in the area of the spanning discs. Clearly for sufficiently small $t$ and each $\lambda \in \Lambda$, $A(h_t(D_\lambda)) < C$. It follows that for sufficiently small $t$,

$$\cup_{\gamma \in [\gamma]} S(\phi, \gamma, C) \subset \cup_{\gamma \in [\gamma]} S(h_t \circ \phi, \gamma, C).$$

By the choice of $C$ there is a $\delta > 0$ such that the classes in $\pi_2(N)$ that can be represented by piecewise $C^1$ maps with energy less than $2(C + \delta)$ are the same as the classes that can be represented piecewise $C^1$ maps with energy less than $2(C - \delta)$. It follows that for each closed curve $\gamma$ on $\Sigma$ that represents $[\gamma]$ we have:

$$S(\phi, \gamma, C) = S(\phi, \gamma, C + \delta).$$

Thus for sufficiently small $t$,

$$\cup_{\gamma \in [\gamma]} S(h_t \circ \phi, \gamma, C) \subset \cup_{\gamma \in [\gamma]} S(\phi, \gamma, C).$$

The result follows. □

**Remark:** The lemma remains true for any 1-parameter family of diffeomorphisms that preserves the lagrangian condition and preserves the periods of curves on the lagrangian surface. In particular, the lemma is true for 1-parameter families with variational vector field $X$ of the type used in the proof of the monotonicity formula. Such vector fields have support in a disc on $\Sigma$ and have lifts that are contact vector fields. Thus they satisfy: $d(X \lrcorner \, \omega)|_\Sigma = 0$ and therefore $X \lrcorner \, \omega$ is a closed 1-form on $\Sigma$ with support in a disc. It follows that $\int_\gamma X \lrcorner \, \omega = 0$ for any closed curve $\gamma$. This implies that the period of $\gamma$ is preserved. On the other hand the lemma is not true for arbitrary 1-parameter families of symplectic diffeomorphisms since these do not preserve the periods of curves.

**Corollary 5.11.** *If* $\mathrm{period}_C(\phi, [\gamma]) = 0$ *then there is a simple closed curve* $\gamma \in [\gamma]$ *and an immersed disc* $D$ *with* $\partial D = \phi(\gamma)$ *and* $\int_D \omega = 0$.

*Proof.* From the proof of the lemma we have that $\cup_{\gamma \in [\gamma]} S(\phi, \gamma, C)$ is a finite set of values and each value is achieved by a simple closed curve in $[\gamma]$. The corollary follows. □

The $\mathrm{period}_C$ invariant is defined for weakly lagrangian maps in $W^{1,2}(\Sigma, N)$. Let $\ell : \Sigma \to N$ be such a map and suppose $\ell_*([\gamma]) = 1$. Let $C_r^0(S^1, N)$ denote the curves $S^1 \to N$ that are continuous and rectifiable. Define:

$$\mathrm{period}_C(\ell, [\gamma]) = \inf_{\{\gamma \in [\gamma]: \, \ell \circ \gamma \in C_r^0(S^1, N)\}} \mathrm{per}_C(\ell, \gamma).$$

This invariant is preserved under weak convergence. In particular, let $\{\ell_i\}$ be a sequence of weakly lagrangian maps in $W^{1,2}(\Sigma, N)$ with $\ell_{i*}([\gamma]) = 1$ for all $i$ and $E(\ell_i) < E$ for some constant $E$. Suppose that $\mathrm{period}_C(\ell_i, [\gamma]) = p$ for each $i$. Let $\ell$ be the weak limit in $W^{1,2}(\Sigma, N)$ of a subsequence of $\{\ell_i\}$ (that we will continue to denote $\{\ell_i\}$). Let $\gamma_0 \in [\gamma]$ be a simple closed curve such that $\ell(\gamma_0)$ is continuous and rectifiable. Let $T$ be a tubular neighborhood of $\gamma_0$ in $\Sigma$ with $\Gamma : S^1 \times [-1, 1] \to T$ a smooth immersion. For $s \in [-1, 1]$, let $\gamma_s : S^1 \to \Sigma$ be the curve $\gamma_s(t) = \Gamma(t, s)$. For almost every $s$, the curve $\ell(\gamma_s)$ is continuous and rectifiable. Given any $\varepsilon > 0$, if $s$ is sufficiently small, the curves $\ell(\gamma_0)$ and $\ell(\gamma_s)$ bound an annulus of area less than $\varepsilon$. Moreover, for almost every $s$, $\{\ell_i(\gamma_s)\}$ is continuous, rectifiable and has uniformly bounded length. Thus, there is a subsequence of $\{\ell_i\}$ (that we will continue to denote $\{\ell_i\}$) that converges uniformly to $\ell(\gamma_s)$. In particular, for $i$ sufficiently large, $\ell_i(\gamma_s)$ is uniformly close to $\ell(\gamma_s)$. Therefore, for $i$ sufficiently large, the curves $\ell_i(\gamma_s)$ and $\ell(\gamma_0)$ bound an annulus of area less than $2\varepsilon$. It follows that:

$$\mathrm{period}_C(\ell, [\gamma]) \leq p.$$

The proof of the opposite inequality is similar.

**Remarks:**
(i) Lemma 5.10 applies to weakly lagrangian maps in $W^{1,2}(\Sigma, N)$. The proof is essentially the same.

(ii) Suppose that $\{\ell_i\}$ is a minimizing sequence of maps in $W^{1,2}(\Sigma, N)$ all with the same period$_C$ invariants and with weak limit $\ell$. Then it follows that $\ell$ also has these period$_C$ invariants. The replacement argument of Proposition 5.5 can be made to preserve the period$_C$ invariants and therefore constructs an admissible comparison map among maps with fixed periods. It follows that a subsequence of $\{\ell_i\}$ converges strongly to $\ell$ on $\Sigma \setminus S$, where $S$ is a finite set of points. The details are left to the reader.

Let $\Sigma$ be a closed surface of genus greater than one. Fix $\delta > 0$. Suppose that $\phi : \Sigma \to N$ is a $C^1$ lagrangian map such that either $\phi_* : \pi_1(\Sigma, p) \to \pi_1(N, \phi(p))$ contains no nontrivial element represented by a simple closed curve in its kernal or if such an element $[\gamma]$ does exist then period$_C(\phi, [\gamma]) \geq \delta > 0$. We require that $C$ satisfy $2A(\phi) < C$ and that $2C$ is a regular value of the energy functional $E(S^2, N)$ on piecewise $C^1$ maps $S^2 \to N$. We call such a lagrangian map *non-collapsible*. The notion is a generalization of incompressibility as formulated in [SU 2].

In analogy to the above definitions we define:

$$\tilde{\mathcal{L}}_\phi = \ \{\ell \in W^{1,2}(\Sigma, N) : \ell \text{ is weakly lagrangian,}$$
$$\ell \text{ has the same action on } \pi_1 \text{ as } \phi \text{ and}$$
$$\text{the period}_C \text{ invariants of } \ell \text{ and } \phi \text{ are equal.}\}$$

Set,
$$\tilde{E}_\mu = \inf\{E(\ell, \mu) : \ell \in \tilde{\mathcal{L}}_\phi\}.$$

**Proposition 5.12.** *There exists an $\ell \in \tilde{\mathcal{L}}_\phi$ such that $E(\ell, \mu) = \tilde{E}_\mu$.*

*Proof.* Same as above. $\square$

Suppose that genus$(\Sigma) = g \geq 2$. Denote the Teichmüller space of $\Sigma$ by $\mathcal{T}(\Sigma)$. $\mathcal{T}(\Sigma)$ can be described as the equivalence classes of pairs $(\Sigma_1, f_1)$ where $f_1 : \Sigma \to \Sigma_1$ is a homeomorphism and two pairs $(\Sigma_1, f_1)$ and $(\Sigma_2, f_2)$ are equivalent if $f_2 \circ f_1^{-1}$ is homotopic to a conformal map. We can mark $\Sigma$ by choosing a set of simple closed curves on $\Sigma$ whose homotopy classes generate $\pi_1(\Sigma)$. A marked surface represents an element of Teichmüller space under the equivalence relation that two marked surfaces are equivalent if there is a homeomorphism respecting the markings and homotopic to a conformal map. Let $\mathcal{R}(\Sigma)$ denote the Riemann moduli space of conformal structures on $\Sigma$. Then $\mathcal{R}(\Sigma) = \mathcal{T}(\Sigma)/\text{Mod}(\Sigma)$ where $\text{Mod}(\Sigma)$, the Teichmüller modular group of $\Sigma$, is the mapping class group modulo its subgroup of ineffective elements. The following is a result of Abikoff [A], (quoted in [SU 2]).

**Theorem 5.13.** *There exists a compactification $\bar{\mathcal{R}}(\Sigma)$ of $\mathcal{R}(\Sigma)$ such that the boundary points $\bar{\mathcal{R}}(\Sigma) \setminus \mathcal{R}(\Sigma)$ correspond precisely to the Riemann surfaces with nodes that can be obtained from $\Sigma$ by collapsing a set of admissible closed curves on $\Sigma$.*

Let $\Sigma_0$ denote a fixed Riemann surface of genus $g \geq 2$ and let $\Sigma_0$ be the base surface of $\mathcal{T}(\Sigma_0)$. We define an energy function on Teichmüller space as follows: Given a lagrangian map $\phi : \Sigma_0 \to N$ consider a point $\mu = (\Sigma, f) \in \mathcal{T}(\Sigma_0)$. Define:

$$\tilde{E}_\mu = \inf\{E(\ell, \mu) : \ell \circ f \in \tilde{\mathcal{L}}_\phi\}.$$

The conformal invariance of $E$ insures that $\tilde{E}_\mu$ is well-defined independent of the representative of $\mu$. This defines a map:

$$\tilde{E}: \mathcal{T}(\Sigma_0) \to \mathbb{R}$$
$$\mu \mapsto \tilde{E}_\mu.$$

This function is lower semi-continuous.

**Proposition 5.14.** *$\tilde{E}_\mu$ is attained at each $\mu \in \mathcal{T}(\Sigma_0)$ by some weakly lagrangian $W^{1,2}(\Sigma, N)$ map $\ell_\mu$ which has the same induced map on $\pi_1(\Sigma)$ and the same* period$_C$ *invariants as $\phi$.*

**Theorem 5.15.** *Let $\ell_i : \Sigma_{\mu_i} \to N$ be a sequence of maps in $W^{1,2}$ for the conformal structure $\mu_i$ on $\Sigma$ such that $E(\ell_i, \mu_i) < B$. Suppose that $\pi(\Sigma_{\mu_i}) \to \Sigma_\infty$ where $\pi : \mathcal{T}(\Sigma_0) \to \mathcal{R}(\Sigma_0)$ is the quotient map and $\Sigma_\infty \in \bar{\mathcal{R}}(\Sigma_0) \setminus \mathcal{R}(\Sigma_0)$. Then for $i$ sufficiently large, $(\ell_i)_*([\gamma]) = 1$ for at least one $[\gamma] \in \pi_1(\Sigma)$ represented by a simple closed curve. Moreover,*

$$\liminf_i \text{period}_B(\ell_i, [\gamma]) = 0.$$

*Proof.* The argument is based on the proof of a similar result in [SU 2]. The sequence $\Sigma_{\mu_i}$ gives a deformation of $\Sigma_0$ onto a Riemann surface $\Sigma_\infty$ with nodes $w_m$, $m = 1, \ldots, r$, where each node corresponds to the pinching of a homotopically nontrivial simple closed curve $\gamma_m$. There is a sequence $\{D_j^m\}$ of closed annular neighborhoods of $\gamma_m$ such that $D_j^m$ converges to the node $w_m$ of $\Sigma_\infty$ and, for each fixed $j$, the change in conformal structure on $\Sigma_0$ as $\Sigma_{\mu_i} \to \Sigma_\infty$ is restricted to the interior of $\cup_{m=1}^r D_j^m$. Let $\Sigma_j = \Sigma_0 \setminus \cup_m D_j^m$. Then, for fixed $j$, $\{\ell_{i|\Sigma_j}\}$ is a sequence of $W^{1,2}$ maps with bounded $\mu_j$ energy. Choosing a subsequence of the sequence $\{\ell_i\}$, that we continue to denote $\{\ell_i\}$, we can suppose that $\{\ell_{i|\Sigma_1}\}$ converges weakly in $W^{1,2}$, strongly in $L^2$ and pointwise almost everywhere to a $W^{1,2}$ map $\ell^{(1)} : \Sigma_1 \to N$. Choosing successive subsequences for, $j = 2, 3, \ldots$, we construct a sequence of $W^{1,2}$ maps $\ell^{(j)} : \Sigma_j \to N$, with the property that $\ell^{(k)}$ extends $\ell^{(j)}$ from $\Sigma_j$ to $\Sigma_k$ for each $k > j$. Let $j \to \infty$. We get a $W^{1,2}$ map $\ell : \Sigma'_\infty \to N$, where $\Sigma'_\infty$ is the punctured Riemann surface $\Sigma_0 \setminus \{w_1, \ldots, w_r\}$ and $\ell^{(j)} \to \ell$ strongly in $W^{1,2}_{loc}$. Since $E(\Sigma_j, \mu_j) < B$ it follows that $E(\ell) < B$. Trivially we can extend $\ell$ to a $W^{1,2}$ map $\tilde{\Sigma}_\infty \to N$ where $\tilde{\Sigma}_\infty = \Sigma'_\infty \cup \{(x_1, y_1), \ldots, (x_r, y_r)\}$ is the closed Riemann surface obtained by adding the pair of points $(x_m, y_m)$ at the two punctures of $\Sigma'_\infty$ corresponding to each node $w_m \in \Sigma_\infty$.

Fix $m$. Let $\gamma$ be a curve homotopic to $\gamma_m$ and contained in $D_j^m$ for some large $j$. Then $\gamma \subset \Sigma_\infty$ is homotopically trivial. Let $T$ be a tubular neighborhood of $\gamma$ in $\Sigma_\infty \setminus \{w_m\}$ with $\Gamma : S^1 \times [-1, 1] \to T$ a smooth immersion. For $s \in [-1, 1]$, let $\gamma_s : S^1 \to \Sigma$ be the curve $\gamma_s(t) = \Gamma(t, s)$. For almost every $s$, $\ell(\gamma_s)$ is continuous, rectifiable and homotopically trivial. Also, for almost every $s$, there is a subsequence of $\{\ell_i\}$, that we continue to denote $\{\ell_i\}$, such that the curves $\{\ell_i(\gamma_s)\}$ converge uniformly to $\ell(\gamma_s)$. Therefore, for sufficiently large $i$ and almost every $s$, $\ell_i(\gamma_s)$ is homotopically trivial.

We next show that $\liminf_i \operatorname{period}_B(\ell_i, [\gamma]) = 0$. Choose any $\varepsilon > 0$. Then $\gamma \subset \Sigma_\infty$ can be chosen so that $\ell(\gamma)$ is continuous and rectifiable and so that a subsequence of the curves $\{\ell_i(\gamma)\}$ converges uniformly to $\ell(\gamma)$. $\ell(\gamma)$ bounds the image of a disc $D$ under the weakly lagrangian map $\ell$. Clearly, $E(\ell_{|D}) < B$. Using smooth approximation, we can find a smoothly immersed disc $D'$ spanning $\ell(\gamma)$ with $\int_{D'} \omega < \varepsilon$ and $A(D') < B$. It follows that for $i$ sufficiently large, there is a smoothly immersed disc $D'_i$ spanning $\ell_i(\gamma)$ with $\int_{D'_i} \omega < 2\varepsilon$ and $A(D'_i) < B$. The result follows. □

**Lemma 5.16.** *If $\mu \in \mathcal{T}(\Sigma)$ is a critical point for the function $\tilde{E}$, then the weakly lagrangian map $\ell_\mu$ is weakly conformal.*

*Proof.* A variation of the domain shows that the Hopf differential is holomorphic and hence smooth. The result then follows from an argument in [SU 1]. □

**Theorem 5.17.** *Suppose $\Sigma$ is a surface of genus $g \geq 2$ and $\phi : \Sigma \to N$ is a lagrangian map that is non-collapsible. Then there is a conformal structure on $\Sigma$ and a weakly conformal, weakly lagrangian map $\ell \in \tilde{\mathcal{L}}_\phi$ that minimizes area among all maps in $\tilde{\mathcal{L}}_\phi$.*

*Proof.* Let $\nu$ be the conformal structure on $\Sigma$ for the metric induced by $\phi$. Then $E(\phi, \nu) = A(\phi) < C$. In particular, $\tilde{E}_\nu < C$. To prove the theorem it suffices to show, by Proposition 5.4, that a minimum of $\tilde{E}$ is attained in the interior of $\mathcal{T}(\Sigma)$. Let $\mu_i \in \mathcal{T}(\Sigma)$ be a minimizing sequence for $\tilde{E}$. We can suppose that $\tilde{E}_{\mu_i} < C$. Then there is a subsequence, still denoted $\mu_i$, and a sequence of elements $\tau_i \in \operatorname{Mod}(\Sigma)$ such that $\xi_i = \tau_i(\mu_i)$ either converges to an element in the interior of $\mathcal{T}(\Sigma)$ or the marked surface $\Sigma_{\xi_i}$ representing $\xi_i$ converges to a surface with nodes [A]. Let $\ell_i$ be a minimizing map in $\tilde{\mathcal{L}}_{\phi \cdot \tau_i^{-1}}$ for the conformal structure $\xi_i$. Clearly $E(\ell_i, \xi_i) = \tilde{E}_{\mu_i} < C$. Therefore since a nontrivial simple closed curve in $\Sigma$ either maps to a nontrivial curve in $N$ or, if not, to a curve with $\operatorname{period}_C$ bounded below, by Theorem 5.15 we have $\lim_{i \to \infty} \xi_i = \xi \in \mathcal{T}(\Sigma)$. Choosing a subsequence, we can suppose that $\{\ell_i\}$ converges weakly in $W^{1,2}$ to $\ell$ and $\liminf_{i\to\infty} \tilde{E}_{\mu_i} = E(\ell, \xi)$. It follows from Lemma 5.16 that $\ell$ is weakly conformal. □

We next consider the case where genus $(\Sigma) = 1$. Let $T^2$ be a closed surface of genus one. We include the period condition in the following definition however it will not be necessary in our applications. Fix $\delta > 0$. Suppose that $\phi : T^2 \to N$ is a piecewise $C^1$ lagrangian map such that if $[\gamma]$ is a generator of $\pi_1(T^2, p)$ that lies in the kernal of $\phi_* : \pi_1(T^2, p) \to \pi_1(N, \phi(p))$ then $\text{period}_C(\phi, [\gamma]) \geq \delta > 0$. We require, in addition, that $C$ satisfy $2A(\phi) < C$ and that $2C$ is a regular value of the energy functional $E(S^2, N)$ on piecewise $C^1$ maps $S^2 \to N$. We call such a lagrangian map *non-collapsible*. Note that every generator of $\pi_1(T^2, p)$ can be represented by a simple closed curve.

In analogy to the above definitions we define:

$$\tilde{\mathcal{L}}_\phi = \ \{\ell \in W^{1,2}(T^2, N) : \ell \text{ is weakly lagrangian,}$$
$$\ell \text{ has the same action on } \pi_1 \text{ as } \phi \text{ and}$$
$$\text{the period}_C \text{ invariants of } \ell \text{ and } \phi \text{ are equal.}\}$$

The proof of the following theorem is also a variation of the work of [SU2] and is left to the reader.

**Theorem 5.18.** *Suppose $T^2$ is a torus and $\phi : T^2 \to N$ is a non-collapsible lagrangian map. Then there is a conformal structure on $T^2$ and a weakly conformal, weakly lagrangian map $\ell \in \tilde{\mathcal{L}}_\phi$ that minimizes area among all maps in $\tilde{\mathcal{L}}_\phi$.*

**Theorem 5.19.** *Let $\Sigma$ be a surface of genus $g \geq 1$ and $\phi : \Sigma \to N$ be a lagrangian map that is non-collapsible. Let $\alpha$ be a lagrangian homology class. Suppose that there are maps $f \in \tilde{\mathcal{L}}_\phi$ such that $[f] = \alpha$. Then there is a conformal structure on $\Sigma$ and weakly conformal, weakly lagrangian maps $\ell : \Sigma \to N$ and $s_i : S^2 \to N$, $i = 1, \ldots, \nu$, such that:*

1. *Each map $s_i$ minimizes area among lagrangian maps in some free homotopy class.*

2. *The map $\ell$ is non-collapsible and minimizes area among all maps in $\tilde{\mathcal{L}}_\phi$.*

3. $[\ell] + \sum_i^\nu [s_i] = \alpha$.

*Proof.* The result follows by applying the argument in the proof of Proposition 5.7 to an area minimizing sequence of maps in $\tilde{\mathcal{L}}_\phi$ that represent $\alpha$. □

**Minimizing in a homology class**

**Theorem 5.20.** *Suppose that $\alpha$ is a lagrangian class in $H_2(N; \mathbb{Z})$. Then $\alpha$ can be represented by a finite number of piecewise $C^1$ lagrangian maps of surfaces into $N$ each of which is non-collapsible.*

*Proof.* Without loss of generality, we assume that $\alpha$ cannot be written as a sum of two nontrivial classes each of which is lagrangian. We represent $\alpha$ by a piecewise $C^1$ map $\ell : \Sigma \to N$. We can suppose that $\ell$ is an immersion except in the neighborhood of finitely many points and, by perturbation, we can suppose that, $\ell$ has, at worst, double points. We take the surface $\Sigma$ to have minimal genus among such surfaces representing $\alpha$. Since $\alpha$ is a lagrangian class,

$$\int_\Sigma f^*\omega = 0, \tag{5.1}$$

If $\gamma$ is a simple closed curve on $\Sigma$ that does not separate and for which $\ell(\gamma)$ is trivial in $\pi_1(N)$ then $\ell(\gamma)$ can be contracted without changing the homology class of $\ell(\Sigma)$ and preserving (5.1). Thus we can assume that all such curves are nontrivial. If for every simple closed curve $\gamma$ on $\Sigma$ which separates we have that $\ell_*([\gamma]) \neq 0$ then the theorem is proved. So suppose that $\gamma$ be a simple closed curve on $\Sigma$ which separates and for which $\ell_*([\gamma])$ is trivial in $\pi_1(N, x)$. We can also suppose $\ell(\gamma)$ is imbedded. $\ell(\gamma)$ bounds an immersed disc $D$ in $N$. If this disc is lagrangian then $\ell(\Sigma)$ can be cut along $\ell(\gamma)$ into two lagrangian surfaces, contradicting the indecomposability of $\alpha$. If $\int_D \omega = 0$ then $D$ is homotopic rel $\partial D = \ell(\gamma)$ to an immersed lagrangian disc. Thus we can suppose that for any such $\gamma$ there are no immersed discs $D$ in $N$ spanning $\ell(\gamma)$ such that $\int_D \omega = 0$. In particular, $\text{period}_C(\ell, [\gamma]) \neq 0$, where we choose $C > 0$ to satisfy $2A(\ell(\Sigma)) < C$ and such that $2C$ is a regular value of the energy functional on piesewise $C^1$ maps $S^2 \to N$.

There can be many homotopy classes in $\pi_1(\Sigma, p)$ that are represented by simple closed separating curves. Let $[\gamma_1]$ and $[\gamma_2]$ be two such classes and let $\gamma_1$ and $\gamma_2$ be simple closed curves that are representatives. Recalling the notation introduced in Lemma 5.10 , for $i = 1, 2$, let

$$S(\ell, [\gamma_i], C) = \cup_{\gamma \in [\gamma_i]} S(\ell, \gamma, C).$$

Let $D_i$ be a disc in $N$ with $A(D_i) < C$ and $\partial D_i = \ell(\gamma_i)$. Since $\gamma_1$ and $\gamma_2$ are both separating curves they bound a region $U \subset \Sigma$. In $N$ construct the cycle

$$\Sigma_{D_1, D_2} = D_1 \cup_{\ell(\gamma_1)} \ell(U) \cup_{\ell(\gamma_2)} D_2.$$

Since $\ell$ is lagrangian,

$$\omega([\Sigma_{D_1, D_2}]) = \int_{D_1} \omega - \int_{D_2} \omega.$$

Making different choices of spanning discs $D_1$ and $D_2$ we construct different cycles and get different numbers after pairing with $[\omega]$. However note that $A(\Sigma_{D_1, D_2}) < 2C + A(\ell(\Sigma))$. Using Lemma 5.3 there are at most finitely many homology classes in $H_2(N, \mathbb{Z})$ that can be represented by cycles with area less than $2C + A(\ell(\Sigma))$. Pairing these homology classes with $[\omega]$ we get finitely many scalars $T = \{t_1, \ldots, t_\lambda\}$. It follows that each element of $S(\ell, \gamma_1, C)$ differs from one of $S(\ell, \gamma_2, C)$ by addition or subtraction by some element of $T$. This is true for each representative $\gamma_1$ of $[\gamma_1]$ and $\gamma_2$ of $[\gamma_2]$ and therefore it is true for the elements of the sets $S(\ell, [\gamma_1], C)$ and $S(\ell, [\gamma_2], C)$. Fix a homotopy class $[\gamma_0] \in \pi_1(\Sigma, p)$ that can be represented by a simple closed separating curve. Let $[\gamma] \in \pi_1(\Sigma, p)$ be any other such class. Then the elements of $S(\ell, [\gamma], C)$ can be obtained from the elements of $S(\ell, [\gamma_0], C)$ by addition or subtraction by an element of $T$. Hence the union of the sets $S(\ell, [\gamma], C)$, over all classes $[\gamma] \in \pi_1(\Sigma, p)$ that can be represented by simple closed separating curves, has finitely many elements. Denote the infimum by $\delta$. If $\delta = 0$ then $\ell(\Sigma)$ can be cut into two lagrangian surfaces, contradicting the indecomposibility of $\alpha$. This proves the theorem. □

Let $\alpha$ be a lagrangian homology class. By the previous result we can represent $\alpha$ by a sum of non-collapsible lagrangian surfaces and lagrangian 2-spheres. Applying the minimization arguments of this section to each map in this decomposition we have:

**Theorem 5.21.** *Let $\alpha \in H_2(N; \mathbb{Z})$ be a lagrangian homology class. Then there exist Riemann surfaces $\Sigma_1, \ldots, \Sigma_\mu$ and hamiltonian stationary, weakly conformal, weakly lagrangian maps $\ell_j : \Sigma_j \to N$, $j = 1, \ldots, \mu$, and $s_i : S^2 \to N$, $i = 1, \ldots, \nu$ such that:*

1. *Each map $s_i$ minimizes area among lagrangian maps in a free homotopy class.*

2. *Each map $\ell_j$ minimizes area among lagrangian maps $\Sigma_j \to N$ with fixed induced action on $\pi_1$ and fixed periods.*

3. *The maps $\ell_j : \Sigma_j \to N$ are non-collapsible.*

4. *$\sum_j^\mu [\ell_j] + \sum_i^\nu [s_i] = \alpha$, where $[\ell_j], [s_i] \in H_2(N; \mathbb{Z})$.*

## 6. Second Variation Formula

In this section we derive the second variation formula for immersed lagrangian submanifolds that are hamiltonian stationary. However, to begin we digress to state and prove the general second variation formula for smoothly immersed submanifolds possibly with boundary. We do not assume that the variation vector field is normal or that the immersed submanifold is stationary.

Let $M$ be a Riemannian manifold of dimension $m$ with metric $g$. Let $\ell : \Sigma \to M$ be a smooth immersion, where $\Sigma$ is a compact manifold of dimension $n$ with boundary (possibly empty). A smooth mapping
$$L : (-\varepsilon, \varepsilon) \times \Sigma \to M$$
satisfying

(i) each map $\ell_t = L(t, -) : \Sigma \to M$ is a smooth immersion
(ii) $\ell_0 = \ell$

will be called a *smooth variation* of $\ell$. Let $\frac{\partial}{\partial t}$ denote the vector field along the $(-\varepsilon, \varepsilon)$ factor; and denote $X = L_*(\frac{\partial}{\partial t})|_{t=0}$. $X$ is a vector field along $\ell(\Sigma)$ called the *variation vector field*.

Let $p \in \Sigma$. With respect to the metric induced on $\Sigma$ by $\ell$, choose an orthonormal frame $\{e_1 \cdots e_n\}$ to satisfy:
$$\nabla^\Sigma_{e_i} e_j(p) = 0 \qquad \text{for all } i, j,$$
where $\nabla^\Sigma$ denotes the Levi-Civita connection on $\Sigma$ for the metric induced by $\ell$. For each $t \in (-\varepsilon, \varepsilon)$, the vector fields $\{\ell_{t*}e_1, \ldots, \ell_{t*}e_n\}$ give a framing along $\ell_t(\Sigma)$. Note that the metric on $\Sigma$ induced by $\ell_t$ can be written
$$g_{ij}(t) = \langle \ell_{t*}e_i, \ell_{t*}e_j \rangle, \quad 1 \le i, j \le n.$$
For brevity of notation we will write
$$e_i(t) = \ell_{t*}e_i, \quad 1 \le i \le n.$$
Let $dv_t$ denote the volume form of the metric induced by $\ell_t$. Set
$$V(t) = \int_\Sigma dv_t.$$
Let $B$ denote the second fundamental form on $\ell(\Sigma)$. $B$ is a symmetric bilinear form on $T\Sigma$ with values in the normal bundle. $R$ will denote the curvature tensor of $g$ on $M$ and $H$ the mean curvature vector.

**Theorem 6.1. (Second Variation Formula)** *Let $\ell : \Sigma \to M$ be a smoothly immersed submanifold, possibly with boundary. Let $L : (-\varepsilon, \varepsilon) \times \Sigma \to M$ be a smooth variation of $\ell$ with variation vector field $X$. Then*

$$\left.\frac{d^2V}{dt^2}\right|_{t=0} = \int_\Sigma \left[ \sum_i |(\nabla_{e_i} X)^\perp|^2 + \sum_i \langle R_{e_i, X} e_i, X \rangle - \langle \nabla_X X, H \rangle \right.$$

$$\left. - \sum_{i,j} \langle \nabla_{e_j} X, e_i \rangle \langle \nabla_{e_i} X, e_j \rangle + \sum_{i,j} \langle \nabla_{e_i} X, e_i \rangle \langle \nabla_{e_j} X, e_j \rangle \right] dv_0$$

$$- \int_{\partial \Sigma} \langle \nabla_X X, \nu \rangle dv_{\partial \Sigma}$$

*where $\{e_1, \ldots, e_n\}$ is an orthonormal frame on $\ell_0(\Sigma)$ and $\nu$ is the inward pointing conormal vector field along $\ell_0(\partial \Sigma)$.*

*Proof.* The well-known computation of the first variation formula shows that:

$$\frac{d}{dt}dv_t = \sum_{i,j} g^{ij}(t)\langle \nabla_{e_i} X, e_j\rangle dv_t,$$

where $g^{ij}(t)$ is the inverse matrix of the metric $g_{ij}(t) = \langle e_i(t), e_j(t)\rangle$. It follows that,

$$\begin{aligned}
\left.\frac{d^2 V}{dt^2}\right|_{t=0} &= \left.\frac{d}{dt}\int_\Sigma \sum_{i,j} g^{ij}(t)\langle \nabla_{e_i} X, e_j\rangle dv_t\right|_{t=0} \\
&= \int_\Sigma \frac{d}{dt}\left[\sum_{i,j} g^{ij}(t)\langle \nabla_{e_i} X, e_j\rangle\right]\bigg|_{t=0} dv_0 + \sum_i \langle \nabla_{e_i} X, e_i\rangle \frac{d}{dt}dv_t\bigg|_{t=0}
\end{aligned}$$

Since $g^{ij}(0) = \delta_{ij}$,

$$\left.\frac{d}{dt}dv_t\right|_{t=0} = \sum_i \langle \nabla_{e_i} X, e_i\rangle dv_0,$$

Using $[X, e_i] = 0, \quad i = 1, \ldots, n$ we have,

$$\frac{dg^{ij}}{dt}(0) = -\frac{dg_{ij}}{dt}(0) = -X\langle e_i, e_j\rangle = -\langle \nabla_{e_i} X, e_j\rangle - \langle e_i, \nabla_{e_j} X\rangle.$$

Thus,

$$\begin{aligned}
\left.\frac{d^2 V}{dt^2}\right|_{t=0} = \int_\Sigma &\bigg[\sum_{ij}(-\langle \nabla_{e_i} X, e_j\rangle - \langle e_i, \nabla_{e_j} X\rangle)\langle \nabla_{e_i} X, e_j\rangle \\
&+ \sum_i \langle \nabla_X \nabla_{e_i} X, e_i\rangle + \sum_i \langle \nabla_{e_i} X, \nabla_X e_i\rangle + \sum_{i,j}\langle \nabla_{e_i} X, e_i\rangle\langle \nabla_{e_j} X, e_j\rangle\bigg] dv_0
\end{aligned}$$

Note that,

$$\langle \nabla_X \nabla_{e_i} X, e_i\rangle = \langle R_{X,e_i} X, e_i\rangle + \langle \nabla_{e_i} \nabla_X X, e_i\rangle.$$

Thus at a point $p \in M$ where $\nabla_{e_i} e_j(p) = 0$ we have,

$$\begin{aligned}
\sum_i \langle \nabla_{e_i} \nabla_X X, e_i\rangle &= \sum_i e_i \langle \nabla_X X, e_i\rangle - \sum_i \langle \nabla_X X, \nabla_{e_i} e_i\rangle \\
&= \sum_i e_i \langle \nabla_X X, e_i\rangle - \langle \nabla_X X, H\rangle,
\end{aligned}$$

where $H$ is the mean curvature. It follows that,

$$\left.\frac{d^2V}{dt^2}\right|_{t=0} = \int_\Sigma \left[ -\sum_i |(\nabla_{e_i} X)^\top|^2 + \sum_i |\nabla_{e_i} X|^2 + \sum_i \langle R_{X,e_i} X, e_i \rangle \right.$$

$$\left. -\langle \nabla_X X, H \rangle - \sum_{i,j} \langle e_i, \nabla_{e_j} X \rangle \langle e_j \nabla_{e_i} X \rangle + \sum_{i,j} \langle \nabla_{e_i} X, e_i \rangle \langle \nabla_{e_j} X, e_j \rangle + \sum_i e_i \langle \nabla_X X, e_i \rangle \right] dv_0$$

The result now follows by applying the divergence theorem to give:

$$\int_\Sigma \sum_i e_i \langle \nabla_X X, e_i \rangle dv_0 = -\int_{\partial \Sigma} \langle \nabla_X X, \nu \rangle dv_0,$$

where $\nu$ is the inward pointing conormal. $\square$

We now return to the lagrangian case. Let $M$ be a Kähler manifold of complex dimension $n$ with Kähler form $\omega$, Kähler metric $g$ and complex structure $J$. Let $\ell : \Sigma \to M$ be a lagrangian immersion, where $\Sigma$ is a compact manifold with boundary (possibly empty). A smooth mapping

$$L : (-\varepsilon, \varepsilon) \times \Sigma \to M$$

satisfying

(i) each map $\ell_t = L(t, -) : \Sigma \to M$ is a lagrangian immersion
(ii) $\ell_0 = \ell$

will be called a *smooth lagrangian variation* of $\ell$. The variation vector field $X = L_*(\frac{\partial}{\partial t})|_{t=0}$ will be called a *lagrangian variation*. Analogous terms will be used in the hamiltonian case.

**Lemma 6.2.** *Suppose $\Sigma$ is a hamiltonian stationary submanifold. If $X$ denotes the variational vector field along $\Sigma$ of a hamiltonian variation then*

$$\int_\Sigma (\langle H, \nabla_X X \rangle - \langle X, \nabla_{JH} JX \rangle) \, dvol = 0$$

*where $H$ denotes the mean curvature vector field along $\Sigma$ and $\nabla$ is the Levi-Civita connection of $g$ on $M$.*

*Proof.* Since $X$ is hamiltonian there is a smooth function $h$ such that $JX = \nabla h$. Define a 1-form $\sigma$ on $\Sigma$ by:

$$\sigma(V) = \langle J(\nabla_X X), V \rangle - \langle X, \nabla_V JX \rangle,$$

for $V \in T(\Sigma)$. Then,

$$\begin{aligned}
\sigma(V) &= \langle J(\nabla_X X), V\rangle - \langle X, \nabla_V JX\rangle \\
&= X\langle JX, V\rangle - \langle JX, \nabla_X V\rangle + \langle JX, \nabla_V X\rangle \\
&= X\langle JX, V\rangle + \langle JX, [V, X]\rangle. \\
&= X\langle \nabla h, V\rangle + \langle \nabla h, [V, X]\rangle. \\
&= X(V(h)) + [V, X](h). \\
&= V(X(h)).
\end{aligned}$$

It follows that $\sigma = d(X(h))_{|\Sigma}$, in particular, $\sigma$ is an exact 1-form. Using $\sigma$ in the first variational formula we get:

$$\int_\Sigma (\langle H, \nabla_X X\rangle - \langle X, \nabla_{JH} JX\rangle) \mathrm{dvol} = 0.$$

This proves the lemma. $\square$

Using a slightly different argument the following result can be proved. However since we won't require it we leave the proof to the reader.

**Lemma 6.3.** *Suppose $\Sigma$ is a lagrangian stationary submanifold. If $X$ denotes the variational vector field along $\Sigma$ of a lagrangian variation then*

$$\int_\Sigma (\langle H, \nabla_X X\rangle - \langle X, \nabla_{JH} JX\rangle) \ dvol = 0$$

*where $H$ denotes the mean curvature vector field along $\Sigma$ and $\nabla$ is the Levi-Civita connection of $g$ on $M$.*

The next theorem is due to Y. G. Oh [O] who derived it, and applied it to prove the stability of certain hamiltonian stationary submanifolds. We include the proof here for the convenience of the reader.

**Theorem 6.4.** *(Y. G. Oh) Let $\ell : \Sigma \to M$ be a hamiltonian stationary submanifold. Let $L : (-\varepsilon, \varepsilon) \times \Sigma \to M$ be a smooth hamiltonian variation of $\ell$ with variation vector field $X$, that leaves the boundary fixed. Denote by*

$X^\perp$ the normal part of $X$. Then

$$\left.\frac{d^2V}{dt^2}\right|_{t=0} = \int_\Sigma \left[\sum_i |\nabla^\Sigma_{e_i} JX^\perp|^2 + \sum_i \langle R_{e_i,X^\perp}e_i, X^\perp\rangle + \langle X^\perp, H\rangle^2\right.$$

$$\left. - \sum_{i,j}\langle X^\perp, B_{e_i,e_j}\rangle^2 - \langle X^\perp, B_{JH,JX^\perp}\rangle\right] dv_0$$

$$= \int_\Sigma \left[|\delta\sigma|^2 - \mathrm{Ric}(X^\perp, X^\perp) + \langle X^\perp, H\rangle^2 - 2\langle X^\perp, B_{JH,JX^\perp}\rangle\right] dv_0$$

where $\sigma = X^\perp \lrcorner\, \omega$ is a closed 1-form, $\delta$ is the adjoint of $d$ on $\Sigma$ and Ric is the Ricci curvature on $N$.

*Proof.* The first equality follows from the lemma. To see the second equality we note that since $\Sigma$ is lagrangian, for vector fields $X, Y, Z$ on $\Sigma$, the second fundamental form satisfies the symmetry:

$$\langle B_{X,Y}, JZ\rangle = \langle B_{X,Z}, JY\rangle.$$

Using this and the Gauss equation we have:

$$\sum_i \langle \mathcal{R}_{e_i,JX^\perp}e_i, JX^\perp\rangle \tag{6.1}$$

$$= \sum_i \langle R_{e_i,JX^\perp}e_i, JX^\perp\rangle - \langle X^\perp, B_{JH,JX^\perp}\rangle + \sum_{i,j}\langle X^\perp, B_{e_i,e_j}\rangle^2,$$

where $\mathcal{R}$ denotes the curvature tensor on $\Sigma$. The Weitzenböck formula relating the Hodge and covariant laplacians on 1-forms on $\Sigma$ gives:

$$\Delta_h \sigma = -\nabla^2 JX^\perp + \mathcal{R}ic(JX^\perp, JX^\perp), \tag{6.2}$$

where $\Delta_h$ is the Hodge laplacian and $\mathcal{R}ic$ denotes the Ricci curvature on $\Sigma$. Combining (6.1), (6.2) and integrating by parts gives the result. □

## 7. Two Dimensional Stationary Lagrangian Cones

In this section we give a complete description of the stationary lagrangian cones in $\mathbb{R}^4$. We apply the second variation formula to study their stability properties and conclude with a discussion of their minimization properties.

To begin consider the unit sphere in $\mathbb{R}^{2n}$, endowed with its induced $CR$ structure as well as its induced metric. To describe this $CR$ structure, let $\vec{p}$ denote the position vector at a point $p \in \mathbb{S}^{2n-1}$. Then $T = J\vec{p}$ is a tangent vector field, and the orthogonal complement $\Pi_p = T^\perp$ in the sphere is a $J$-invariant $2n-2$ dimensional subspace. This is the contact (CR) distribution in the sphere. If we consider a cone $C(\Gamma)$ over an $n-1$ dimensional submanifold $\Gamma$ of the unit sphere, then we see that $C(\Gamma)$ is lagrangian if and only if $\Gamma$ is legendrian in the sphere. This is true because the tangent plane to the cone at points of $\Gamma$ is spanned by the tangent space to $\Gamma$ together with the position vector, and the legendrian condition guarantees that the position vector is taken by $J$ to a vector orthogonal to the cone. In particular, we see that a two dimensional cone is lagrangian if and only if $\Gamma$ is a legendrian curve.

We now consider a curve $\gamma$ in $\mathbb{S}^3$. We introduce complex coordinates $z_1, z_2$ and we write $\gamma(s) = (\gamma_1(s), \gamma_2(s))$ where $\gamma_j$ are $l$ periodic complex valued functions and $s$ is an arclength parameter with $l = \text{length}(\gamma)$. The condition that the cone over $\gamma$ be lagrangian is that the $2 \times 2$ matrix with columns $\gamma$ and $\dot{\gamma}$ is a unitary matrix. We have

$$e^{i\beta} = \gamma_1 \dot{\gamma}_2 - \gamma_2 \dot{\gamma}_1$$

where $\beta$ is the lagrangian angle. In order that the cone be hamiltonian stationary, $\beta$ must be a harmonic function. On the other hand $\beta$ is homogeneous of degree 0 on the cone. It follows that $\beta$ is a linear function of $s$ which we can take to be $2as$ for a real constant $a$. These conditions imply that

$$\dot{\gamma}_1 = -e^{i2as}\bar{\gamma}_2, \ \dot{\gamma}_2 = e^{i2as}\bar{\gamma}_1. \tag{7.1}$$

If we differentiate the first equation and substitute in from the second we get

$$\ddot{\gamma}_1 = i2a\dot{\gamma}_1 - \gamma_1.$$

Thus we have

$$\gamma_1(s) = c_1 e^{i(a+\sqrt{a^2+1})s} + c_2 e^{i(a-\sqrt{a^2+1})s}$$

for complex constants $c_1, c_2$. It follows that

$$\gamma_2(s) = i(a + \sqrt{a^2+1})\bar{c}_1 e^{i(a-\sqrt{a^2+1})s} + i(a - \sqrt{a^2+1})\bar{c}_2 e^{i(a+\sqrt{a^2+1})s}.$$

Since $\gamma$ is to be a curve on the unit sphere, we must have the condition

$$|\gamma_1|^2 + |\gamma_2|^2 = |c_1|^2(1 + (a + \sqrt{a^2+1})^2) + |c_2|^2(1 + (a - \sqrt{a^2+1})^2) = 1.$$

In order that $\gamma$ be a closed curve, it must be true that $a + \sqrt{a^2+1}$ and $a - \sqrt{a^2+1}$ are rationally related. One can check that this condition is equivalent to saying that there are relatively prime positive integers $p, q$ such that $a + \sqrt{a^2+1} = \sqrt{p/q}$ and $a - \sqrt{a^2+1} = -\sqrt{q/p}$. We then have that the length of the curve satisfies $l = 2\pi\sqrt{pq}$, and if we set $\theta = s/\sqrt{pq}$, then we have $\gamma(\theta)$ is a $2\pi$ periodic parametrization. Furthermore we have $\beta = (\sqrt{p/q} - \sqrt{q/p})s = (p-q)\theta$. Thus we find that $\beta$ has a period around $\gamma$ which is equal to $2\pi(p-q)$. Applying a unitary transformation we can take $c_1$ to be real and $c_2 = 0$, so we see that $\gamma$ is unitarily equivalent to the curve,

$$\gamma(s) = \frac{1}{\sqrt{p+q}}\left(\sqrt{q}e^{i\sqrt{\frac{p}{q}}s}, i\sqrt{p}e^{-i\sqrt{\frac{q}{p}}s}\right), \tag{7.2}$$

where $0 \leq s \leq 2\pi\sqrt{pq}$. We note that these curves are not great circles except in the special case in which $p = q = 1$. In general they lie on Clifford tori in $\mathbb{S}^3$. When both $p > 1$ and $q > 1$ the curves are knotted. If either $p = 1$ or $q = 1$ they are unknotted. Thus we have shown:

**Theorem 7.1.** *The lagrangian cones $C_\gamma$ in $\mathbb{R}^4$ that are hamiltonian stationary are the cones over the curves $\gamma$ of (7.2). They are parameterized by a pair of relatively prime positive integers $(p, q)$.*

The second variation formula takes a particularly simple form when applied to variations of the cone, $C_\gamma$, over $\gamma$ whose variation vector field $X$ is hamiltonian, i.e. when $X = J\nabla f$ for $f \in C_c^\infty(C_\gamma)$.

**Proposition 7.2.** *Suppose that the cone $C_\gamma$ over $\gamma$ is lagrangian and hamiltonian stationary. Let $X$ be a compactly supported hamiltonian variation, $X = J\nabla f$, $f \in C_c^\infty(C_\gamma)$ that leaves the cone vertex fixed. Then the second variation formula with variation vector field $X$ is:*

$$\left.\frac{d^2 V}{dt^2}\right|_{t=0} = \int_{C_\gamma} \left[ (\Delta f)^2 - r^{-4}\frac{(q-p)^2}{pq}(f_s)^2 \right] r\, dr\, ds. \qquad (7.3)$$

*where $\Delta$ is the Laplacian on $C_\gamma$.*

*Proof.* Let $r$ denote the radial coordinate on $C_\gamma$. With respect to the coordinates $r$ and $s$ the induced metric on $C_\gamma$ has the form:

$$d\sigma^2 = dr^2 + \frac{1}{r^2} ds^2$$

Set,

$$e_1 = \gamma, \quad e_2 = \dot\gamma$$

Then $\{e_1, e_2\}$ is an orthonormal frame on $C_\gamma$. With respect to this frame an easy computation shows that the second fundamental form satisfies:

$$B_{e_1,e_1} = B_{e_1,e_2} = 0, \quad B_{e_2,e_2} = \frac{1}{r}(\ddot\gamma + \gamma). \qquad (7.4)$$

In particular, $H = \frac{1}{r}(\ddot\gamma + \gamma)$. From (7.4) it follows immediately that

$$\langle H, X \rangle^2 - \sum_{i,j} \langle X, B_{e_i,e_j}\rangle^2 = 0$$

To compute $\langle X, \nabla_{JH} JX \rangle$, we note that, from the explicit form of $\gamma$ (7.2), we have:

$$JH = \frac{1}{r}\frac{q-p}{\sqrt{pq}}\dot\gamma = \frac{1}{r}\frac{q-p}{\sqrt{pq}}e_2.$$

Clearly,

$$JX = -\nabla f = -f_r e_1 - \frac{1}{r}f_s e_2.$$

Thus,
$$\begin{aligned}\langle X, \nabla_{JH} JX\rangle &= \langle X, B_{JH,JX}\rangle \\ &= \langle X, -\frac{1}{r^2}\frac{(q-p)}{\sqrt{pq}} f_s B_{e_2,e_2}\rangle \\ &= -\frac{1}{r^2}\frac{(q-p)}{\sqrt{pq}} f_s \langle JX, JH\rangle \\ &= \frac{1}{r^4}\frac{(q-p)^2}{pq}(f_s)^2\end{aligned}$$

The result follows. $\square$

**Remark 7.1.** The coordinate functions are easily seen to be Jacobi fields.

We next apply Proposition 7.2 to study the stability properties of the lagrangian cones for various pairs $p, q$.

**Proposition 7.3.** *The stationary lagrangian cones with $|p - q| > 1$ are strictly unstable. This is true for lagrangian variations fixing a neighborhood of the cone vertex.*

*Proof.* Choose a positive integer $\ell$ satisfying:
$$\ell(-|p-q| + \ell) < pq < \ell(|p-q| + \ell) \tag{7.5}$$
Set
$$f(r,s) = \zeta(r)\cos\frac{\ell s}{\sqrt{pq}} \qquad 0 \leq s \leq 2\pi\sqrt{pq}$$
where $\zeta(r) \in C_c^\infty(\mathbb{R}_+)$ will be specified later. Then we have,
$$\begin{aligned}\Delta f &= \Delta\zeta\cos\frac{\ell s}{\sqrt{pq}} - r^{-2}\zeta\frac{\ell^2}{pq}\cos\frac{\ell s}{\sqrt{pq}} \\ &= \left(\zeta'' + \frac{\zeta'}{r} - \zeta\frac{\ell^2}{pq}r^{-2}\right)\cos\frac{\ell s}{\sqrt{pq}}\end{aligned}$$
Using $f$ in the second variation formula (7.3) we have:
$$\left.\frac{d^2 V}{dt^2}\right|_{t=0} = \tag{7.6}$$
$$\int\left[\left(\zeta'' + \frac{\zeta'}{r} - \zeta\frac{\ell^2}{pq}r^{-2}\right)^2\cos^2\frac{\ell s}{\sqrt{pq}} - r^{-4}\frac{(p-q)^2}{pq}\zeta^2\frac{\ell^2}{pq}\sin^2\frac{\ell s}{\sqrt{pq}}\right]r\,dr\,ds$$

Set
$$g(r,s) = \zeta(r)\sin\frac{\ell s}{\sqrt{pq}} \qquad 0 \leq s \leq 2\pi\sqrt{pq}$$

where $\zeta(r)$ is as above. Using $g$ in the second variation formula (7.3) we have:

$$\left.\frac{d^2V}{dt^2}\right|_{t=0} = \tag{7.7}$$

$$\int\left[\left(\zeta'' + \frac{\zeta'}{r} - \zeta\frac{\ell^2}{pq}r^{-2}\right)^2 \sin^2\frac{\ell s}{\sqrt{pq}} - r^{-4}\frac{(p-q)^2}{pq}\zeta^2\frac{\ell^2}{pq}\cos^2\frac{\ell s}{\sqrt{pq}}\right]r\,dr\,ds$$

We wish to show that for suitable choice of $\zeta$ one of (7.6), (7.7) is negative. To do this add the integrals to give:

$$2\pi\sqrt{pq}\int\left[\left(\zeta'' + \frac{\zeta'}{r} - \zeta\frac{\ell^2}{pq}r^{-2}\right)^2 - r^{-4}\frac{(p-q)^2}{(pq)^2}\ell^2\zeta^2\right]r\,dr \tag{7.8}$$

Next choose $\varepsilon$, $0 < \varepsilon < 1$ and define $\zeta$:

$$\zeta(r) = \begin{cases} \delta(r) & 0 \leq r \leq \varepsilon \\ r & \varepsilon \leq r \leq 1 \\ \eta(r) & r \geq 1 \end{cases}$$

where $\eta$ satisfies:

(i) $\eta$ has support on $[1,2]$
(ii) $\eta(1) = 1$, $\eta'(1) = 1$, $\eta''(1) = 0$
(iii) $0 \leq \eta \leq c$, $|\eta'| \leq c$, $|\eta''| \leq c$ for some constant $c > 0$.

and where $\delta$ satisfies:

(i) $\delta$ has support on $[\frac{\varepsilon}{2}, \varepsilon]$
(ii) $\delta(\varepsilon) = \varepsilon$, $\delta'(\varepsilon) = 1$, $\delta''(\varepsilon) = 0$
(iii) $0 \leq \delta \leq \varepsilon$, $|\delta'| \leq 4$, $|\delta''| \leq \frac{4}{\varepsilon}$

Thus (7.8) becomes three integrals:

$$2\pi\sqrt{pq}\int_{\frac{\varepsilon}{2}}^{\varepsilon}\left[\left(\delta'' + \frac{\delta'}{r} - \delta\frac{\ell^2}{pq}r^{-2}\right)^2 - r^{-4}\frac{(q-p)^2}{(pq)^2}\ell^2\delta^2\right]r\,dr$$

$$+\; 2\pi\sqrt{pq}\int_{\varepsilon}^{1}\left[(pq-\ell^2)^2 - \ell^2(p-q)^2\right]\frac{1}{(pq)^2}\frac{dr}{r}$$

$$+\; 2\pi\sqrt{pq}\int_{1}^{2}\left[\left(\eta'' + \frac{\eta'}{r} - \eta\frac{\ell^2}{pq}r^{-2}\right)^2 - r^{-4}\frac{(q-p)^2}{(pq)^2}\ell^2\eta^2\right]r\,dr$$

The absolute value of the first integral is bounded by
$$C\int_{\frac{\varepsilon}{2}}^{\varepsilon}\frac{dr}{\varepsilon} = \frac{C}{2}$$

where $C$ depends on $\ell, p, q$ but is independent of $\varepsilon$. The third integral is clearly bounded. The second integral equals

$$-2\pi(pq)^{-\frac{3}{2}}\left[(pq-\ell^2)^2 - \ell^2(p-q)^2\right]\ln\varepsilon$$

By our choice of $\ell$ this expression is negative and so, for $\varepsilon$ sufficiently small, the sum of all three integrals is negative. We conclude that one (or both) of (7.6), (7.7) is negative. □

Suppose now that $|p - q| = 1$. Without loss of generality we can suppose that $q = p + 1$.

**Proposition 7.4.** *The stationary lagrangian cones with $|p - q| = 1$ are strictly stable for lagrangian variations fixing the cone vertex.*

*Proof.* Without loss of generality we can assume that the lagrangian variation is a normal hamiltonian variation, i.e., has the form $X = J\nabla f$ for some $f \in C_c^\infty(C_\gamma \setminus \{0\})$ satisfying $f(0) = f_r(0) = 0$. We first argue that the second variation is positive for variations of the form $X = J\nabla f_\ell$ or $X = J\nabla g_\ell$ where,

$$f_\ell(r, s) = \zeta(r)\cos\frac{\ell s}{\sqrt{pq}}, \quad 0 \leq s \leq 2\pi\sqrt{pq},$$

or

$$g_\ell(r, s) = \zeta(r)\sin\frac{\ell s}{\sqrt{pq}}, \quad 0 \leq s \leq 2\pi\sqrt{pq}.$$

for any positive integer $\ell$. $\zeta(r)$ is assumed to have compact support and to satisfy:

$$\zeta(0) = \zeta'(0) = 0$$

(so that $X(0) = 0$). The computation in the proof of Proposition 7.3 shows that the second variation for both $f_\ell$ and $g_\ell$ reduces to a positive multiple of:

$$\int_0^\infty \left[\left(\zeta'' + \frac{\zeta'}{r} - \zeta\frac{\ell^2}{pq}r^{-2}\right)^2 - r^{-4}\frac{1}{(pq)^2}\ell^2\zeta^2\right]r\,dr. \tag{7.9}$$

Make the substitutions,

$$r = e^t, \quad \zeta(r) = e^t\rho(t).$$

The boundary conditions on $\rho$ are:

$$\rho(t) \to 0, \quad \rho'(t) \to 0, \quad \text{as} \quad t \to \pm\infty.$$

Using the same argument as in the proof of Proposition 7.4 the integral (7.9) becomes:

$$\int_{-\infty}^\infty \left[(\rho'')^2 + (2 + \frac{2p}{p+1})(\rho')^2 + \frac{1}{(p(p+1))^2}[(p(p+1) - \ell^2)^2 - \ell^2]\rho^2\right]dt.$$

Since,

$$[(p(p+1) - \ell^2)^2 - \ell^2] \geq 0,$$

for any $\ell \in \mathbb{Z}_+$, this integral is positive for any $\rho$.

If $f \in C_c^\infty(C_\gamma \setminus \{0\})$ satisfies $f(0) = f_r(0) = 0$ then $f$ has a Fourier series expansion:

$$f(r,s) = \sum_{\ell=0}^\infty \left( \zeta_\ell(r) \cos \frac{\ell s}{\sqrt{pq}} + \eta_\ell(r) \sin \frac{\ell s}{\sqrt{pq}} \right),$$

with

$$\zeta_\ell(0) = \zeta_\ell'(0) = \eta_\ell(0) = \eta_\ell'(0) = 0, \quad \zeta, \eta \in C_c^\infty((0,\infty)).$$

It follows using the $L^2$ perpendicularity of sine and cosine and the above computation that,

$$dV(J\nabla f, J\nabla f) > 0.$$

The result follows. $\square$

The situation for multiply covered cones is markedly different.

**Proposition 7.5.** *The multiply covered stationary lagrangian cones for any $p,q$ are strictly unstable. This is true for lagrangian variations fixing a neighborhood of the cone vertex.*

*Proof.* The k times covered $(p,q)$ cones are parameterized by:

$$\frac{1}{\sqrt{p+q}} \left( r\sqrt{q} e^{i\sqrt{\frac{p}{q}}s}, ir\sqrt{p} e^{-i\sqrt{\frac{q}{p}}s} \right), \tag{7.10}$$

where $0 \leq s \leq 2k\pi\sqrt{pq}$, $0 \leq r < \infty$ and $k > 1$ is an integer. Without loss of generality we can suppose $q > p$. Put:

$$\ell = p + \frac{1}{k}. \tag{7.11}$$

Set:

$$f(r,s) = \zeta(r) \cos \frac{\ell s}{\sqrt{pq}} \qquad 0 \leq s \leq 2k\pi\sqrt{pq}$$

and

$$g(r,s) = \zeta(r) \sin \frac{\ell s}{\sqrt{pq}} \qquad 0 \leq s \leq 2k\pi\sqrt{pq}$$

where $\zeta(r) \in C_c^\infty(\mathbb{R}_+)$ is as given in the proof of Proposition 7.3. Then both $f$ and $g$ are compactly supported hamiltonians on the k covered $(p,q)$ cone. Using the argument in the proof of Proposition 7.3 it suffices to show

$$(pq - \ell^2)^2 - \ell^2(p-q)^2 < 0 \tag{7.12}$$

to conclude that the second variation with respect to $f$ or $g$ (or both) is negative. Using $\ell = p + \frac{1}{k}$ and $q > p$ this follows easily. $\square$

We now show that none of the cones are minimizing if we allow comparisons which are lagrangian and nonorientable. A convenient way to do this is to use the following generalization of Allcock's theorem [Al] due to Weiyang Qiu [Q]. The idea here is that Allcock's method can be extended to give a lagrangian homotopy (with controlled area) of any curve $\Gamma$ to a twice covered circle in a complex plane with the same enclosed symplectic area as that of $\Gamma$. This then represents a lagrangian Möbius band.

**Lemma 7.6.** *(Qiu) Let $\Gamma$ be any rectifiable curve in $\mathbb{R}^4$. There exists a lagrangian Möbius band $\Sigma$ with boundary $\Gamma$ and such that $A(\Sigma) \leq cL(\Gamma)^2$ for a fixed constant $c$.*

The same argument as in the proof of Lemma 4.8 now shows the following.

**Lemma 7.7.** *Suppose that $\ell_0, \ell_1 : S^1 \to \mathbb{R}^4$ are continuous maps with lengths $L_i = L(\ell_i(S^1))$ for $i = 1, 2$, and with $\sup_{S^1} d(\ell_0, \ell_1) \leq \varepsilon < L_0 + L_1$. There exists a (nonorientable) surface $\Sigma$ with two boundary components $C_0, C_1$ and a lagrangian map $\ell : \Sigma \to \mathbb{R}^4$ with $\ell = \ell_i$ on $C_i$ for $i = 1, 2$, and $A(\ell(\Sigma)) \leq c(L_0 + L_1)\varepsilon$.*

With this preparation we may now prove the following.

**Proposition 7.8.** *For any $(p, q)$ with $(p, q) \neq (1, 1)$, the cone $C$ does not minimize area among nonorientable lagrangian comparison surfaces. In fact, there exists a nonorientable lagrangian surface $\Sigma$ with $\partial\Sigma = C \cap \partial B_1$, and $A(\Sigma) < A(C \cap B_1)$.*

*Proof.* Because of the homogeneity of the cone, it clearly suffices to work in a ball of any radius. We let $\Lambda > 1$ be a large number to be determined, and observe that the mean curvature vector $H$ of $C$ is homogeneous of degree $-1$ and is therefore bounded on $C \cap (B_\Lambda \setminus B_1)$. We may then do a normal variation of $C$ by an amount $\varepsilon_0 H$ for a fixed $\varepsilon_0 > 0$ to produce a lagrangian perturbation $C_1$ of $C \cap (B_\Lambda \setminus B_1)$ with $A(C_1) \leq A(C \cap (B_\Lambda \setminus B_1)) - c\log(\Lambda)$ (since the integral of $H^2$ on $C \cap (B_\Lambda \setminus B_1)$ is of order $\log(\Lambda)$). Now we write $\partial C_1 = \Gamma_1 \cup \Gamma_2$ where $\Gamma_1$ is a curve of bounded distance to $C \cap \partial B_1$, and $\Gamma_2$ is of distance at most $c\Lambda^{-1}$ to $C \cap \partial B_\Lambda$. Applying Lemma 7.7, we may find nonorientable lagrangian "strips" $\Sigma_1, \Sigma_2$ of bounded area and with $\partial\Sigma_1 = \Gamma_1 \cup (C \cap \partial B_1)$, $\partial\Sigma_2 = \Gamma_2 \cup (C \cap \partial B_\Lambda)$. Thus we may form the nonorientable surface $\Sigma$ which is equal to a union of $\Sigma_1$, $C_1$, $\Sigma_2$, and $C \cap B_1$. We then have $\partial\Sigma = C \cap \partial B_\Lambda$, and $A(\Sigma) \leq c + A(C_1) \leq A(C \cap B_\Lambda) + c - c\log(\Lambda)$. If $\Lambda$ is chosen sufficiently large this gives us a comparison surface which is smaller than the cone. $\square$

We next consider the minimizing properties (among orientable comparisons) of the non-multiply covered $(p, p+1)$ cones. Set,

$$\gamma(s) = \frac{1}{\sqrt{2p+1}}\left(\sqrt{p+1}e^{i\sqrt{\frac{p}{p+1}}s}, i\sqrt{p}e^{-i\sqrt{\frac{p+1}{p}}s}\right),$$

where $0 \leq s \leq 2\pi\sqrt{p(p+1)}$.

**Theorem 7.9.** *For at least one integer $p \geq 1$, the hamiltonian stationary lagrangian $(p, p+1)$ cone minimizes area among disk type lagrangian comparison surfaces with boundary $\gamma(s)$.*

*Proof.* Let $\alpha \in H_2(N, \mathbb{Z})$ be a lagrangian homology class with $c_1(N)(\alpha) \neq 0$. The existence theorem allows us to represent $\alpha$ by a finite set of weakly lagrangian maps that are area minimizers. For at least one such map $\ell : \Sigma \to N$ we have $c_1(N)(\ell(\Sigma)) \neq 0$. The regularity theorem implies that $\ell$ is a branched immersion except at finitely many cone-type isolated singularities. By Propostion 2.2 twice the sum of the local Maslov indices at these singularities equals $c_1(N)(\ell(\Sigma))$ and is therefore not zero. Hence there must be a non-trivial area minimizing cone. □

**Example:** Perhaps the simplest example of a Kähler surface $(N, \omega)$ and a homology class $\alpha \in H_2(N, \mathbb{Z})$ satisfying $[\omega](\alpha) = 0$ and $c_1(N)(\alpha) \neq 0$ is given as follows: Let $N$ be $\mathbb{P}^2$ with one point blown up. Let $L$ be the class of a line in $\mathbb{P}^2$ and $E$ be the class of the exceptional curve. There is a Kähler form $\omega$ on $N$ satisfying $[\omega](L) = [\omega](E)$ so that $L - E$ is a lagrangian class. The canonical class is $K = 3L - E$ and therefore $c_1(N)(L - E) = 2$. There are, of course, many such examples.

In fact it is possible to show that the $(p, p+1)$ cones occur in situations where $c_1(N) = 0$. To discuss these examples we digress to describe some results due to Hitchin [H] on Einstein four-manifolds. We are indebted to C. Lebrun for referring us to the relevent literature. Hitchen proves the following result: Let $N$ be a compact four-dimensional Einstein manifold with signature $\sigma$ and Euler characteristic $\chi$. Then

$$|\sigma| \leq \frac{2}{3}\chi.$$

If equality occurs then $\pm N$ is either flat or its universal cover is a $K3$ surface. If the universal cover of $N$ is a $K3$ surface then $N$ is a $K3$ surface, an Enriques surface or the quotient of an Enriques surface by a free antiholomorphic involution. In a note added in proof, Hitchin shows that there is a $K3$ surface $X$, given as a complete intersection in $\mathbb{CP}^5$, that admits a free holomorphic involution $\tau_+$ and a free antiholomorphic involution $\tau_-$ that commute and such that $\tau_+ \tau_-$ is free. Also both involutions are isometries of the metric $h$ induced on $X$ by the Fubini-Study metric.

Denote the Kähler form of $h$ by $\omega_h$. Since $\tau_+$ is a holomorphic isometry, $\tau_+^* \omega_h = \omega_h$ and therefore $\omega_h$ descends to a Kähler form on the quotient four manifold $Y_+ = X/\tau_+$. It follows easily that $Y_+$ is an Enriques surface. By Yau's theorem [Y] there exists a unique Ricci-flat Kähler metric $g_+$ on $Y_+$ with Kähler form $\omega_{g_+}$ in the same cohomology class as $\omega_h$. This metric lifts to a Ricci-flat Kähler metric on $X$ that we denote by $g$. Note that $\tau_+$ is a holomorphic isometry of $g$. Next consider the antiholomorphic involution $\tau_-$. Let $\sigma$ denote a non-zero holomorphic $(2,0)$-form on $X$. With respect to $g$, $\sigma$ is a parallel section of the canonical bundle. Then $\tau_-^*\sigma$ is a holomorphic $(2,0)$-form for the conjugate complex structure and so is a parallel section of the anticanonical bundle. Thus $\tau_-^*\sigma = \pm\bar\sigma$. It follows that either Re $\sigma$ or Im $\sigma$ is invariant under $\tau_-$. Suppose that Re $\sigma$ is invariant under $\tau_-$. (The case that Im $\sigma$ is invariant under $\tau_-$ is entirely similar.) By scaling we can assume that Re $\sigma$ has unit length. Then Re $\sigma$ determines a complex structure $J_\sigma$ that lies on the same hyperKähler line as the original complex structure and Re $\sigma$ is the Kähler form for the metric $g$ and $J_\sigma$. Since both Re $\sigma$ and $J_\sigma$ are invariant under $\tau_-$ so is $g$. Therefore both $\tau_+$ and $\tau_-$ are isometries of the Calabi-Yau metric $g$. It follows that there is a Ricci-flat Kähler metric on the Enriques surface $Y_+$ for which $\tau_-$ is an antiholomorphic isometry and a Ricci-flat metric on $Z = Y_+/\tau_-$. Note that $\pi_1(Z) = \mathbb{Z}_2 \times \mathbb{Z}_2$ and $b_2^+(Z) = 0$. Thus $Z$ cannot be Kähler or symplectic. The Kähler form $\omega_{g_+}$ on $Y_+$ is taken to $-\omega_{g_+}$ by $\tau_-$ and so does not descend to $Z$. However the density $|\omega|$, where $\omega = \omega_{g_+}$, is well-defined on $Z$ and we can define the notion of a lagrangian in $Z$. We say a map $\ell : \Sigma \to Z$ is lagrangian if $\ell^*(|\omega|) = 0$. Note that the lift of a lagrangian to $Y_+$ is lagrangian in the usual sense. Most of the existence and regularity results of this paper easily extend to this generalized notion of lagrangian.

Note that $\tau_+ : X \to X$ is the covering involution for the covering map $\pi_+ : X \to Y_+$. The canonical bundle $K_X$ on $X$ is trivial. Let $\xi$ denote a parallel section of $K_X$ of unit length. Then $\tau_+^*(\xi)$ is also a parallel section of unit length. Since $\tau_+$ is a holomorphic involution $\tau_+^*(\xi) = \pm\xi$. However if $\tau_+^*(\xi) = \xi$ then $\xi$ descends to a parallel section of $K_Y$ the canonical bundle on $Y$. Thus $K_Y$ is trivial. This contradiction implies that $\tau_+^*(\xi) = -\xi$.

Let $f : T^2 \to Z$ be an incompressible immersion. The immersion $f$ is double covered by an immersion $\tilde f : T^2 \to Y_+$. Clearly

$$\int_{T^2} \tilde f^* \omega = 0.$$

It follows that in a tubular neighborhood of $\tilde f(T^2)$ we have $\omega = d\eta$, for a one-form $\eta$. Note that $\tau_-^*\eta = -\eta$, so $\eta$ does not descend to $Z$, though the density $|\eta|$ is well-defined in a neighborhood of $f(T^2)$. Suppose that $f(T^2)$ has a fine triangulation such that a neighborhood of $f(T^2)$ is covered by open Darboux balls with the properties that each closed 2-simplex on $f(T^2)$ lies in at least one ball of the cover and such that if a 1-simplex intersects a ball in the cover then it lies entirely in the ball. As in the proof of Proposition 2.1 we can perturb each 1-simplex, $\lambda$, keeping its endpoints fixed, so that on the perturbed 1-simplex, $\lambda'$, we have $\int_{\pm\lambda'} \eta = 0$. Then the integral of $\pm\eta$ around the boundary of each perturbed 2-simplex is zero and the 2-simplex can be replaced by a lagrangian simplex. The resulting piecewise $C^1$ lagrangian map $\ell_0 : T^2 \to Z$ is a perturbation of $f$ and so is incompressible.

We next minimize area among lagrangian maps $\ell : T^2 \to Z$ with the same induced map on $\pi_1$ as that of $\ell_0$. By the existence and regularity theory developed in the previous sections a minimizer $\ell_Z$ exists. The map $\ell_Z : T^2 \to Z$ is a Lipschitz lagrangian map that is smooth except (perhaps) at finitely many points. These points are $(p, p+1)$ cone-type singularities described precisely above. The induced map $\ell_{Z*}$ on $\pi_1$ is the same as that of $\ell_0$. We suppose, by way of contradiction, that $\ell_Z$ is smooth everywhere. Then $\ell_Z$ is double covered by a smooth lagrangian map $\ell_Y : T^2 \to Y_+$ which is equivariant with respect to an involution on $T^2$ and the antiholomorphic involution $\tau_-$. Clearly $\ell_Y$ minimizes area among such lagrangian maps. In particular, since $Y_+$ is Kähler-Einstein and $\ell_Y$ is smooth the mean curvature $H$ is an admissible variation. Using $H$ in the first variation formula it follows that $\ell_Y$ is a minimal ($H = 0$) branched immersion. Since one generator of $\pi_1(T^2)$ is mapped to a non-trivial element of $\pi_1(Y_+) = \mathbb{Z}_2$, $\ell_Y$ is double covered by a minimal lagrangian map $\ell_X : T^2 \to X$. Denote the image of $\ell_X$ by $\Sigma$ and note that $\sigma(\Sigma) = \Sigma$. Since $\Sigma$ is minimal lagrangian it is calibrated by a parallel unit section $\xi$ of $K_X$. Thus,

$$\int_\Sigma \operatorname{Re}(\xi) = \int_\Sigma \mathrm{dvol}.$$

But $\sigma^*(\xi) = -\xi$ and therefore

$$\int_\Sigma \operatorname{Re}(\xi) = -\int_\Sigma \mathrm{dvol}.$$

We conclude that $\int_\Sigma \mathrm{dvol} = 0$ and therefore the map $\ell_X$ is a map to a point. Thus the map $\ell_Z$ is also a map to a point. But this is impossible since it is incompressible. This contradiction implies that $\ell_Z$ cannot be regular everywhere and hence that there are singular points on $\ell_Z$. Taking coverings this implies that in the K3 surface $X$ there is a lagrangian stationary torus that contains singular points.

## 8. Main Results

We can combine the existence Theorem 5.21, the regularity Theorem 4.10 and the instability of the hamiltonian stationary lagrangian $(p, q)$ cones for $|p - q| > 1$ to prove:

**Theorem 8.1.** *Let $(N, \omega, J)$ be a compact symplectic 4-manifold with compatible metric $g$. Then the lagrangian homology is generated by classes that can be represented by lagrangian Lipschitz maps $\ell$ that are branched immersions except at finite number of singular points. The local Maslov index of each singular point is either $1$ or $-1$. Twice the sum of the local Maslov indices equals the first Chern class of $N$ paired with the homology class of the surface. The maps $\ell$ are minimizers of area constrained to lagrangian cycles. The mean curvature $H$ of each surface satisfies a first order elliptic system of "Hodge-type".*

**Lemma 8.2.** *A branched lagrangian immersion in a Kähler-Einstein surface is stationary if and only if it is lagrangian stationary.*

*Proof.* Since the ambient manifold is Kähler-Einstein $\sigma_H = H \lrcorner \, \omega$ is closed. To prove the lemma, we need only show that the mean curvature generates a smooth deformation through lagrangian submanifolds. We first consider the immersion case and then extend the argument to include branch points. Let $\ell : \Sigma \to N$ be a lagrangian immersion. We can extend $\ell$ to be an immersion $L$ from a neighborhood $O$ of the zero section in the normal bundle of $\Sigma$ where we identify the zero section with $\Sigma$. The map $L$ is then a local diffeomorphism, and we can pull back both $g$ and $J$ to $O$ using $L$. It suffices to construct a lagrangian variation of the zero section whose initial derivative is $H$. To accomplish this, let $\Pi : O \to \Sigma$ denote the natural projection map, and let $\sigma = -\Pi^*(\sigma_H)$. There is a unique vector field $V$ such that $\omega(H,X) = g(V,X)$, and since $\sigma$ is closed, $V$ is locally a hamiltonian vector field. Also $V = H$ on $\Sigma$. We let $F_t$ denote the flow defined by $V$ near $\Sigma$. The transformations $F_t$ are then symplectic, and $L(F_t(\Sigma))$ gives a lagrangian variation. If we assume that $\ell$ is lagrangian stationary, then the first variation of area is zero for this variation. The first variation of area is given by minus the $L^2$ inner product of the mean curvature vector field with the variation vector field. Therefore, it is strictly negative unless $H = 0$. Thus $\ell$ is stationary. This proves the lemma for immersions.

Let $p_1, \ldots, p_k \in \Sigma$ be the branch points of $\ell$. Choose $r > 0$ such that the discs, $D_r(p_i)$, of radius $r$ centered at $p_i$ are disjoint. On $D_r(p_i)$ there is a smooth function, $\beta_i$, with $\beta_i(p_i) = 0$ such that $\sigma_H = d\beta_i$. Choose local coordinates $\{x_1, x_2\}$ on $D_r(p_i)$ centered at $p_i$. For $\varepsilon < r$, define a function $(h_i)_\varepsilon$ on $D_r(p_i)$ by,

$$(h_i)_\varepsilon(x) = \begin{cases} 0, & |x| < \varepsilon^2 \\ \frac{\log \frac{|x|}{\varepsilon^2}}{\log \frac{1}{\varepsilon}}, & \varepsilon^2 \leq |x| \leq \varepsilon \\ 1, & |x| > \varepsilon. \end{cases}$$

Then a simple computation yields $\lim_{\varepsilon \to 0} \int |\nabla (h_i)_\varepsilon|^2 d\mu = 0$. Define $\sigma_\varepsilon = d((h_i)_\varepsilon \beta_i)$ on $D_r(p_i)$ and $\sigma_\varepsilon = \sigma_H$ otherwise. The closed 1-form $\sigma_\varepsilon$ vanishes near the branch points. Set $V_\varepsilon \lrcorner \, \omega = \sigma_\varepsilon$. Since $\ell$ is an immersion on $\Sigma \setminus \cup_i D_{\varepsilon^2/2}(p_i)$ and $V_\varepsilon$ vanishes on $\cup_i D_{\varepsilon^2}(p_i)$, for each $\varepsilon > 0$ the above construction gives a lagrangian variation of $\ell$ that fixes $\ell$ near the branch points. Suppose that $\ell$ is lagrangian stationary. Then the first variation of area is zero for each such variation. Thus we have for each $\varepsilon > 0$,

$$\begin{aligned} 0 &= -\int \langle V_\varepsilon, H \rangle d\mu \\ &= -\int_{\Sigma \setminus \cup_{i=1}^k D_\varepsilon(p_i)} |H|^2 d\mu - \sum_i^k \int_{D_\varepsilon(p_i)} \langle V_\varepsilon, H \rangle d\mu \end{aligned} \qquad (8.1)$$

For each $i$:

$$\left| \int_{D_\varepsilon(p_i)} \langle V_\varepsilon, H \rangle d\mu \right|$$

$$\leq \int_{D_\varepsilon(p_i)} (h_i)_\varepsilon |\nabla \beta_i|^2 d\mu + \int_{D_\varepsilon(p_i)} |\beta_i| |\langle \nabla(h_i)_\varepsilon, \nabla \beta_i \rangle| d\mu$$

$$\leq \int_{D_\varepsilon(p_i)} (h_i)_\varepsilon |\nabla \beta_i|^2 d\mu + \left( \int_{D_\varepsilon(p_i)} |\nabla(h_i)_\varepsilon|^2 d\mu \right)^{\frac{1}{2}} \left( \int_{D_\varepsilon(p_i)} \beta_i^2 |\nabla \beta_i|^2 d\mu \right)^{\frac{1}{2}}$$

Thus $\int_{D_\varepsilon(p_i)} \langle V_\varepsilon, H \rangle d\mu \to 0$ as $\varepsilon \to 0$. It follows from (8.1) that $0 = \int_\Sigma |H|^2 d\mu$ and therefore that $\ell$ is stationary. □

**Proposition 8.3.** *Let $(N, \omega, J)$ be a compact Kähler-Einstein surface. Suppose the maps $\ell$ in Theorem 8.1 are all branched immersions. Then the lagrangian homology is generated by classes that can be represented by branched lagrangian immersions that are classical minimal surfaces.*

*Proof.* Let $\ell : \Sigma \to N$ be a minimizer. Since the ambient manifold is Kähler-Einstein the 1-form $\sigma_H$ on $\Sigma$ is harmonic. Thus if $\Sigma = S^2$, $\sigma_H = 0$ and therefore $H = 0$. If genus$(\Sigma) \geq 1$ then $\Sigma$ is non-collapsible. If $\gamma$ is a separating simple closed curve then $[\gamma]$ is zero in homology and therefore $[\sigma_H]([\gamma]) = 0$. In the construction of $\ell$ the period condition is only used for classes that are represented by separating curves on $\Sigma$. It follows that variation by the mean curvature preserves the non-collapsibility condition and therefore $H$ is an admissible variation. The Proposition then follows from the lemma. □

Recall that a homotopy class $\alpha \in \pi_2(N)$ is called *lagrangian* if there is a weakly lagrangian map $\ell \in W^{1,2}(S^2, N)$ that represents $\alpha$. Applying the existence and regularity results we can represent a lagrangian homotopy class $\alpha$ by a sum of Lipschitz lagrangian maps $\ell_\lambda : S^2 \to N$ satisfying the same properties as the maps in Theorem 8.1.

Suppose that $\alpha \in \pi_2(N)$ is a lagrangian homotopy class. We define the *lagrangian area* of $\alpha$:

$$\mathrm{Larea}(\alpha) = \inf\{\mathrm{area}(\ell) \ : \ell \in W^{1,2}_L(S^2, N) \text{ and } \ell \text{ represents } \alpha\}.$$

If $m$ is an integer greater than 1, we say a lagrangian homotopy class $\alpha$ is *m-stable* if:

$$\mathrm{Larea}(m\alpha) = m \, \mathrm{Larea}(\alpha).$$

**Theorem 8.4.** *Suppose that $\alpha \in \pi_2(N)$ is a lagrangian homotopy class that can be represented by a lagrangian minimizer $\ell : S^2 \to N$. If $\alpha$ is m-stable then $\ell$ is a branched immersion (there are no singular points).*

*Proof.* Suppose not. Then $\ell$ contains at least two singular points, one with Maslov index $+1$, one with Maslov index $-1$. Consider an $m$-covering of $S^2$, $\phi : S^2 \to S^2$, branched at the two singular points. The map $\ell \circ \phi$ is lagrangian, represents $m\alpha$ and has $A(\ell \circ \phi) = mA(\ell) = m\, \text{Larea}(\alpha)$. It is therefore an area minimizer. But $\ell \circ \phi$ has two singular points with multiply covered tangent cones and is therefore unstable by Proposition 7.5. □

**Corollary 8.5.** *Under the assumptions of the theorem if, in addition, $N$ is a Kähler-Einstein surface then the minimizer $\ell$ is a branched minimal lagrangian immersion.*

## Appendix A. The Geometric Equations

Suppose that $N$ is a $2n$-manifold with symplectic form $\omega$, almost complex structure $J$ and compatible metric $g$.

**Proposition A.1.** *Let $\ell : \Sigma \to N$ be a lagrangian immersion. Denote by $H$ the mean curvature vector field of $\Sigma$ in $V$ and by Ric, the Ricci 2-form of $g$. Then*
$$d(H \lrcorner\, \omega) = d * \tau + \ell^* Ric.$$
*Here $\tau$ is a one form on $\Sigma$ determined by the torsion of a connection on $N$. The connection is uniquely determined by $g$ and $\omega$.*

*Proof.* Exactly as in the Kähler case, the almost complex structure $J$ determines an orthogonal splitting of the complexified tangent space, $TN \otimes \mathbb{C}$, into $+i$ and $-i$ eigenspaces denoted, $T^{(1,0)}N$ and $T^{(0,1)}N$, respectively. Let,
$$\{e_1, e_2, f_1, f_2\} \tag{A.1}$$
be an oriented orthonormal frame adapted to $\Sigma$ in the sense that,

 (i) $\{e_1, e_2\}$ is an oriented frame of the lagrangian plane $T\Sigma$
 (ii) $\{f_1, f_2\}$ is an oriented frame of the lagrangian plane $T\Sigma^\perp$.
 (iii) $Je_j = f_j$ $j = 1, 2$.

The vectors $\{u_j = e_j - if_j\}$ form a unitary frame of $T^{(1,0)}N$ and the vectors $\{\bar{u}_j = e_j + if_j\}$ form a unitary frame of $T^{(0,1)}N$. Thus we have,
$$g(u_j, u_k) = 0, \quad g(u_j, \bar{u}_k) = 2\delta_{jk}. \tag{A.2}$$
Let $\{\theta_1, \theta_2, \eta_1, \eta_2\}$ be the orthonormal coframe dual to the orthonormal frame (A.1). From (i) it follows that,

$\{\theta_1, \theta_2\}$ is an orthonormal coframe on $L$ for the induced metric. (A.3)

From (ii) it follows that,
$$\eta_1 = \eta_2 = 0 \text{ on } \Sigma. \tag{A.4}$$

From (iii) it follows that,

$$J\eta_j = \theta_j, \ j = 1, 2. \tag{A.5}$$

The 1-forms

$$\omega_j = \theta_j + i\eta_j, \ j = 1, 2, \tag{A.6}$$

form a unitary coframe adapted to $\Sigma$. They are dual to the unitary frame $\{u_1, u_2\}$.

Let $\nabla$ denote a metric compatible connection for $g$. Then,

$$\nabla u_j = \sum_k \omega_{jk} \otimes u_k + \sum_k \tau_{jk} \otimes \bar{u}_k, \ j = 1, 2. \tag{A.7}$$

From (A.2) it follows that,

$$\omega_{jk} + \bar{\omega}_{kj} = 0, \quad \tau_{jk} + \tau_{kj} = 0. \tag{A.8}$$

Thus the unitary coframe (A.6) satisfies the structure equations:

$$d\omega_j = \sum_k \omega_{jk} \wedge \omega_k + \sum_k \tau_{jk} \wedge \bar{\omega}_k. \tag{A.9}$$

In particular $(\omega_{jk})$ is the connection one-form and $(\tau_{jk})$ is the torsion form with respect to the coframe $\{\omega_1, \omega_2\}$. Note that the only non-zero component of the torsion is, $\tau_{12} = -\tau_{21}$. The symplectic form,

$$\omega = \frac{i}{2} \sum_j \omega_j \wedge \bar{\omega}_j.$$

satisfies $d\omega = 0$. This is equivalent to the condition:

$$\tau_{12} \text{ is a form of type } (0, 1). \tag{A.10}$$

We remark that the condition (A.10) implies that the connection defined by (A.7) is unique. Thus the symplectic form and the metric, together, determine a unique connection with torsion $\tau_{12}$.

The connection form $(\omega_{jk})$ can be written

$$\omega_{jk} = \theta_{jk} + i\eta_{jk} \tag{A.11}$$

where $\theta_{jk} = -\theta_{kj}$ and $\eta_{jk} = \eta_{kj}$. The torsion form $(\tau_{jk})$ can be written

$$\tau_{jk} = \sigma_{jk} + i\rho_{jk} \tag{A.12}$$

where $\sigma_{jk} = -\sigma_{kj}$ and $\rho_{jk} = -\rho_{kj}$. The structure equations (A.9) become:

$$d\theta_j = \sum_k (\theta_{jk} + \sigma_{jk}) \wedge \theta_k - \sum_k (\eta_{jk} - \rho_{jk}) \wedge \eta_k, \tag{A.13}$$

$$d\eta_j = \sum_k (\eta_{jk} + \rho_{jk}) \wedge \theta_k + \sum_k (\theta_{jk} - \sigma_{jk}) \wedge \eta_k. \tag{A.14}$$

Since $\eta_j = 0, \ j = 1, 2$, along $\Sigma$ these equations imply that on $\Sigma$:

$$\sum_k (\eta_{jk} + \rho_{jk}) \wedge \theta_k = 0. \tag{A.15}$$

Thus,
$$\eta_{jk} + \rho_{jk} = \sum_l h_{jkl}\theta_l \tag{A.16}$$

where $h_{jkl} = h_{jlk}$. The $\{h_{jkl}\}$ are the components of the second fundamental form of $\Sigma$ in $N$. The mean curvature vector is:
$$H = \sum_{k,j} h_{jkk} f_j. \tag{A.17}$$

Then,
$$J(H) = -\sum_{k,j} h_{jkk} e_j.$$

Dualizing we have,
$$J(H)^\# = -\sum h_{jkk}\theta_j. \tag{A.18}$$

Let $*$ denote the Hodge star operator on $\Lambda^1(\Sigma)$. This operator depends only on the conformal structure on $\Sigma$. From (A.16) we have,
$$J(H)^\# = 2*\rho_{12} - \sum_k \eta_{kk}. \tag{A.19}$$

where $\rho_{12} = \operatorname{Im}\tau_{12}$. By (A.11), $i\eta_{kk} = \omega_{kk}$ and thus we have,
$$H \lrcorner\, \omega = *\tau + i\sum_k \omega_{kk}, \tag{A.20}$$

where we have set $\tau = 2\rho_{12}$. Taking the exterior derivative the result follows. □

Thus,
$$\eta_{jk} + \rho_{jk} = \sum_l h_{jkl}\theta_l \tag{A.16}$$

where $h_{jkl} = h_{jlk}$. The $\{h_{jkl}\}$ are the components of the second fundamental form of $\Sigma$ in $N$. The mean curvature vector is:
$$H = \sum_{k,j} h_{jkk} f_j. \tag{A.17}$$

Then,
$$J(H) = -\sum_{k,j} h_{jkk} e_j.$$

Dualizing we have,
$$J(H)^\# = -\sum h_{jkk}\theta_j. \tag{A.18}$$

Let $*$ denote the Hodge star operator on $\Lambda^1(\Sigma)$. This operator depends only on the conformal structure on $\Sigma$. From (A.16) we have,
$$J(H)^\# = 2*\rho_{12} - \sum_k \eta_{kk}. \tag{A.19}$$

where $\rho_{12} = \operatorname{Im}\tau_{12}$. By (A.11), $i\eta_{kk} = \omega_{kk}$ and thus we have,
$$H \lrcorner\, \omega = *\tau + i\sum_k \omega_{kk}, \tag{A.20}$$

where we have set $\tau = 2\rho_{12}$. Taking the exterior derivative the result follows. □

70 R. SCHOEN AND J. WOLFSON